\documentclass[12pt]{article}
\usepackage{graphicx}
\usepackage{float}
\usepackage{amssymb,amsmath,amsfonts,relsize,mathrsfs}
\usepackage[usenames]{color}
\usepackage{cite}
\usepackage{esvect}
\usepackage{txfonts}

\makeatletter
\renewcommand{\fnum@figure}{Fig. \thefigure}
\makeatother
\newtheorem{theorem}{Theorem}[section]
\newtheorem{definition}[theorem]{Definition}

\newtheorem{proposition}[theorem]{Proposition}
\newtheorem{remark}[theorem]{Remark}
\newtheorem{lemma}[theorem]{Lemma}

\newtheorem{problem}[theorem]{Problem}

\newtheorem{claim}[theorem]{Claim}

\newtheorem{pra}[theorem]{The Projection Algorithm}
\newtheorem{alg}[theorem]{Algorithm}
\newtheorem{statement}[theorem]{Statement}
\newcommand {\Ac}      {{\mathcal A}}
\newcommand {\Bc}      {{\mathcal B}}
\newcommand {\Cc}      {{\mathcal C}}

\newcommand {\Gc}      {{\mathcal G}}
\newcommand {\Hc}      {{\mathcal H}}
\newcommand {\Ic}      {{\mathcal I}}

\newcommand {\Kc}      {{\mathcal K}}

\newcommand {\Mc}      {{\mathcal M}}

\newcommand {\Pc}      {{\mathcal P}}

\newcommand {\Rc}      {{\mathcal R}}
\newcommand {\Sc}      {{\mathcal S}}
\newcommand {\Tc}      {{\mathcal T}}

\newcommand {\Wc}      {{\mathcal W}}

\newcommand {\Hcr}      {{\mathscr H}}
\newcommand {\R}       {{\bf R}}

\newcommand {\RN}      {\R^n}
\newcommand {\RD}      {\R^D}
\newcommand {\RM}      {\R^m}
\newcommand {\RT}      {\R^2}

\newcommand {\lh}      {{\bf h}}

\newcommand {\tF}      {\widetilde{F}}

\newcommand {\tH}      {\widetilde{H}}
\newcommand {\tI}      {\widetilde{I}}

\newcommand {\tK}      {\widetilde{K}}

\newcommand {\tQ}      {\widetilde{Q}}

\newcommand {\tMc}     {\widetilde{{\mathcal M}}}

\newcommand {\tMf}     {\widetilde{\mathfrak M}}

\newcommand {\tf}      {\tilde{f}}

\newcommand {\tu}      {\tilde{u}}
\newcommand {\tv}      {\tilde{v}}

\newcommand {\tx}      {\tilde{x}}

\newcommand {\tz}      {\tilde{z}}

\newcommand {\tlm}     {\tilde{\lambda}}
\newcommand {\trh}     {\tilde{\rho}}

\newcommand {\tell}    {\tilde{\ell}}

\newcommand {\hI}      {\widehat{I}}
\newcommand {\hx}      {\hat{x}}
\newcommand {\Cf}      {{\mathfrak C}}
\newcommand {\Df}      {{\mathfrak D}}

\newcommand {\Mf}      {{\mathfrak M}}

\newcommand {\Rf}      {{\mathfrak R}}

\newcommand {\Tf}      {{\mathfrak T}}

\newcommand {\ve}      {\varepsilon}

\newcommand {\MR}      {(\Mc,\rho)}
\newcommand {\MS}      {\Mf}
\newcommand {\lmv}     {\vec{\lambda}}
\newcommand {\vf}      {\varphi}
\newcommand {\emp}     {\emptyset}

\newcommand {\CRT}     {\Conv(\RT)}
\newcommand {\RCT}     {\Rf(\RT)}

\newcommand {\HR}      {\Hc}

\newcommand {\Tfw}     {\widetilde{\Tf}}

\newcommand {\cbg}    {\,\mathlarger{\mathlarger{\cap}}\,}

\newcommand {\blbig}    {\mathlarger{\mathlarger{\bullet\hspace{0.2mm}}}}

\newcommand {\BXR}     {I_0}

\newcommand {\LTI}     {\ell^2_\infty}

\newcommand {\CNV}     {\Conv(\RT)}

\newcommand {\HPL}     {\Hc\Pc(\RT)}

\newcommand {\ip}[1]   {\langle{#1}\rangle}

\newcommand {\FM}      {|F|_{\Mf}}

\newcommand {\PLG}     {\Pc\Gc_L(\RT)}
\newcommand {\PLGA}    {\Pc\Gc(\RT)}
\newcommand {\DI}      {\rho_W}

\newcommand {\s}       {s^{[1]}}
\newcommand {\st}      {t^{[1]}}

\newcommand {\VST}     {}

\newcommand {\Lip}     {\operatorname{Lip}}

\newcommand {\dhf}     {\operatorname{d_H}}

\newcommand {\Prj}     {\operatorname{\text{\sc Pr}}}
\newcommand {\Prm}     {\operatorname{\text{\bf Pr}}}

\newcommand {\dist}    {\operatorname{dist}}

\newcommand {\Conv}    {\operatorname{Conv}}
\newcommand {\sign}    {\operatorname{sign}}
\newcommand {\cent}    {\operatorname{center}}

\newcommand {\clos}    {\operatorname{closure}}
\newcommand {\lend}    {\operatorname{\it L}}
\newcommand {\rend}    {\operatorname{\it R}}
\newcommand {\SN}      {\operatorname{Sn}}
\newcommand {\smsk}    {\smallskip}
\newcommand {\msk}     {\medskip}
\newcommand {\bsk}     {\bigskip}
\newcommand {\bx}      {\hspace{10mm}$\blacksquare$}
\newcommand {\rbx}     {\hspace{10mm}$\blacktriangleleft$}

\newcommand {\nn}      {\nonumber}
\newcommand {\rf}[1]    {(\ref{#1})}      
\newcommand {\reff}[1] {\ref{#1}}         
\newcommand{\lbl}[1]      {\label{#1}}       
\newcommand{\be}          {\begin{eqnarray}}
\newcommand{\bel}[1]      {\begin{eqnarray} \label{#1}}
\newcommand{\ee}           {\end{eqnarray}}
\newcommand {\SECT}[2] {\section*{\centerline{\normalsize
{\bf #1}}} \setcounter{section}{#2}
\setcounter{theorem}{0}\setcounter{equation}{0}}
\begin{document}
\parindent 1em
\parskip 0mm
\medskip
\centerline{{\bf Efficient Algorithms for Lipschitz Selections of Set-Valued Mappings in $\R^2$: long version}}
\vspace*{5mm}
\centerline{By~ {\sc Pavel Shvartsman}}\vspace*{3 mm}
\centerline {\it Department of Mathematics, Technion - Israel Institute of Technology,}\vspace*{1 mm}
\centerline{\it 32000 Haifa, Israel}\vspace*{1 mm}
\centerline{\it e-mail: pshv@technion.ac.il}
\vspace*{5mm}
\renewcommand{\thefootnote}{ }
\footnotetext[1]{{\it\hspace{-6mm}Math Subject
Classification:} 46E35
\smsk
\par This research was supported by the ISRAEL SCIENCE FOUNDATION (grant No. 520/22).}
\begin{abstract} Let $F$ be a set-valued mapping from an $N$-element metric space $(\Mc,\rho)$ into the family of all closed half-planes in $\RT$. In this paper, we provide an efficient algorithm for a Lipschitz selection of $F$, i.e., a Lipschitz mapping $f:\Mc\to\RT$ such that $f(x)\in F(x)$ for all $x\in\Mc$. Given a constant $\lambda>0$,
this algorithm produces the following two outcomes: (1) The algorithm guarantees that there is no Lipschitz selection of $F$ with Lipschitz constant at most $\lambda$; (2) The algorithm returns a Lipschitz selection of $F$ with Lipschitz constant at most $3\lambda$.
\par The total work and storage required by this selection algorithm are at most $CN^2$ where $C$ is an absolute constant.
\end{abstract}
\vspace*{-15mm}
\renewcommand{\contentsname}{ }
\tableofcontents
\addtocontents{toc}{{\centerline{\sc{Contents}}}\par}
\SECT{1. Introduction.}{1}
\addtocontents{toc}{\hspace*{3.2mm} 1. Introduction.\hfill \thepage\par\VST}
\indent

\par {\bf 1.1 Main results.}
\addtocontents{toc}{~~~~1.1 Main results.\hfill \thepage\par\VST}
\msk
\par In this paper we study some aspects of the following selection problem formulated by C. Fefferman \cite{F-2019} (see also \cite[Section 8.7]{FI-2020}).
\begin{problem}\lbl{LS-PRB}  {\em Let $(\Mc,\rho)$ be an $N$-point metric space. For each $x\in\Mc$, let $F(x) \subset \R^D$ be a convex polytope.
\par {\it How can one compute a map $f: \Mc\to\R^D$ such that $f(x)\in F(x)$ for all $x\in \Mc$, with Lipschitz norm as small as possible up to a factor $C(D)$?}
\par This is a big ill-conditioned linear programming problem.
\par {\it Can we do better than just applying general-purpose linear programming? How does the work of an optimal algorithm scale with the number of points $N$?}}
\end{problem}
\par For this problem, we are looking for an {\it efficient algorithm} that requires minimal use of computer resources.
\par C. Fefferman and B. G. Pegueroles \cite{FP-2019}
provided such an algorithm which tackles a general version of Problem \reff{LS-PRB} related to ``{\it approximate}'' $C^m$-smooth selections. Their approach involves an $N$-point subset $E\subset\RN$ as the metric space $\MR$, where the value $f(x)$ for $x \in E$ is allowed to be in an enlarged version of $F(x)$, specifically
$(1 + \tau) \diamond F(x)$.
Here $\tau>0$ and $(1 + \tau) \diamond F(x)$ is the convex set obtained by dilating $F(x)$ about its center of mass by a factor of $(1+\tau)$.
\par In this paper, we focus on a special case of this result related to $C^1$-selections. By applying the standard transitions between $C^1$-functions and Lipschitz functions, we can reformulate the main result from \cite{FP-2019} for the case $m=1$ into the following equivalent statement.
\begin{theorem}\lbl{FP-LSA} (An Approximate Selection Algorithm \cite{FP-2019}) Let $n,D,L,N$ be positive integers and let $M,\tau>0$. Let $E\subset\RN$ be an $N$-point set, and let $F$ be a mapping on $E$ which to every $x\in E$ assigns a compact, convex polytope $F(x)\subset\RD$. We suppose that each $F(x)$ is
defined by at most $L$ linear constraints. Under these conditions, we produce one of the following two outcomes:
\smsk
\par \underline{\sc Outcome 1} {\bf(``No-Go'')}: We guarantee that there does not exist a Lipschitz function $f:E\to\RD$ with Lipschitz constant at most $M$ such that $f(x)\in F(x)$ for all $x\in E$.
\smsk
\par \underline{\sc Outcome 2} {\bf(``Success'')}: We return a Lipschitz function  $f:E\to\RD$, with Lipschitz constant at most $\gamma(n,D) M$ such that $f(x)\in(1+\tau)\lozenge F(x)$ for all $x\in E$.
\smsk
\par This Selection Algorithm requires at most
$CN\log N$ computer operations and at most $CN$ units of computer memory. Here, the constant $C$ depends only on $n,D,L$ and $\tau$.
\end{theorem}
\par Our paper is inspired by this striking result. As the authors of \cite{FP-2019} point out, the approach developed in that paper does not work for the case $\tau=0$, i.e., for ``exact'' (rather than ``approximate'') Lipschitz selections.
\par In our paper, we study the Lipschitz selection Problem \reff{LS-PRB} for the case $\tau=0$ and $D=2$, i.e., for the special case of set-valued mappings defined on finite metric spaces and taking values in the family of all {\it convex polygons} in $\RT$.
\smsk
\par In what follows, it will be more convenient for us to consider the problem in a somewhat broader context of {\it pseudometric} spaces $(\Mc,\rho)$, that is, we assume that the ``distance function'' $\rho:\Mc\times\Mc\to [0,+\infty]$ satisfies
$$
\rho(x,x)=0,~ \rho(x,y)=\rho(y,x),~~~\text{and}~~~\rho(x,y)\le \rho(x,z)+\rho(z,y)
$$
for all $x,y,z\in\Mc$. Note that $\rho(x,y)=0$ may hold with $x\ne y$, and $\rho(x,y)$ may be $+\infty$. We call a pseudometric space $(\Mc,\rho)$ finite if $\Mc$ is finite, but we say that the pseudometric $\rho$ is finite if $\rho(x,y)$ is finite for every $x,y\in\Mc$.
\par Also note that each polygon is the intersection of a finite number of half-planes. This allows us, instead of Problem \reff{LS-PRB}, to study its equivalent version, in which the metric space and the family of convex polygons are replaced by a pseudometric space and the family $\HPL$ of {\it all closed half-planes in $\RT$}.
\smsk
\par Let us formulate our main result. By $\Lip(\Mc)$ we denote the space of all Lipschitz mappings from $\Mc$ into $\RT$ equipped with Lipschitz seminorm
$$
\|f\|_{\Lip(\Mc)}=\inf\{\,\lambda\ge 0:\|f(x)-f(y)\|
\le\lambda\,\rho(x,y)~~\text{for all}~~x,y\in\Mc\,\}.
$$
\par Hereafter $\|\cdot\|$ denotes the uniform norm in $\RT$, i.e., $\|x\|=\max\{|x_1|,|x_2|\}$ for $x=(x_1,x_2)\in\RT$.
\smsk
\par  Let $F$ be a set-valued mapping which to each element $x\in\Mc$ assigns a nonempty convex closed set $F(x)\subset\RT$. A {\it selection} of $F$ is a map $f:\Mc\to \RT$ such that $f(x)\in F(x)$ for all $x\in\Mc$.
A selection $f$ is said to be Lipschitz if $f\in\Lip(\Mc)$.
\par We introduce the quantity $\FM$ by letting
\bel{FM}
\FM=\inf\{\,\|f\|_{\Lip(\Mc)}: f~~\text{is a Lipschitz selection of}~~F\}
\ee
whenever $F$ has a Lipschitz selection, and we set $\FM=+\infty$ otherwise.
\smsk
\par Thus, we are given a positive parameter $\lambda$, an $N$-element pseudometric space $\Mf=\MR$ and a set-valued mapping $F:\Mc\to\HPL$. We want to either compute a Lipschitz function $f:\Mc\to\RT$ such that $f(x)\in F(x)$ on $\Mc$, with $\|f\|_{\Lip(\Mc)}\le \gamma\,\lambda$, or show that there is no Lipschitz selection of $F$ with Lipschitz constant at most $\lambda$. (Here $\gamma>0$ is an absolute constant.)
\par Our objective is to develop an efficient algorithm, suitable for implementation on an idealized computer with a standard von Neumann architecture and the ability to handle exact real numbers, to solve the problem outlined above. By {\it efficient}, we mean requiring a small number of basic computer operations. We consider an ``operation'' to be either an arithmetic operation (addition, subtraction, multiplication, division, comparison of two real numbers) or accessing a value in RAM. The ``running time'' or ``work'' of our algorithm will be assessed by the total count of these operations.
\smsk
\par In Sections 1.2 and Sections 3-6 we introduce and analyze an algorithm that meets all these requirements. We call it the {\it Projection Algorithm}. This straightforward, five-step geometric procedure is the main object of our study. Its construction relies on fundamental geometric operations, including intersections of half-planes in $\RT$, the rectangular hulls of these intersections, and metric projections onto closed convex sets (in the uniform norm in $\RT$). See Section 1.2 and Figures 2-5.
\smsk
\par Let us formulate the main result of the paper.
\begin{theorem}\lbl{PA-PG-I} The Projection Algorithm
receives as input a real number $\lambda>0$, an $N$-element pseudometric space $\MR$ and a half-plane $F(x)\subset\RT$ for each $x\in\Mc$.
\par Given the above input, we produce one of the two following outcomes:
\smsk
\par {$\blbig$} \underline{\sc Outcome 1} {\bf(``No Go''):} We guarantee that there does not exist $f\in\Lip(\Mc)$ with Lipschitz seminorm $\|f\|_{\Lip(\Mc)}\le \lambda$ such that $f(x)\in F(x)$ for all $x\in\Mc$.
\smsk
\par {$\blbig$} \underline{\sc Outcome 2} {\bf(``Success''):} The algorithm produces a mapping $f:\Mc\to\RT$ with Lipschitz seminorm
\bel{G-3}
\|f\|_{\Lip(\Mc)}\le 3\lambda
\ee
satisfying $f(x)\in F(x)$ for all $x\in\Mc$.
\smsk
\par This algorithm requires at most $CN^2$ computer operations and at most $CN^2$ units of computer memory. Here $C>0$ is an absolute constant.
\end{theorem}
\par Theorem \reff{PA-PG-I} implies an analog of Theorem
\reff{FP-LSA} (for the case $\tau=0$ and $D=2$) for set-valued mappings defined on a pseudometric space and taking values in the family of all convex (not necessarily bounded) polygons in $\RT$. This is Theorem \reff{PLG} which we prove in Section 6.6.

\begin{remark} {\em Efficiency estimates of the Projection Algorithm given in Theorem \reff{PA-PG-I} lead naturally to the following interesting (open) problem: Given a specific, non-trivial $N$-element (pseudo)metric space $\MR$, under what conditions can we construct a Lipschitz selection algorithm with a signi\-fi\-cantly lower computational cost than the $O(N^2)$ complexity of the Projection Algorithm presented in Theorem \reff{PA-PG-I}? For instance, let $\MR$ be an $N$-point subset of $\RN$ equipped with the standard Euclidean metric. {\it Can we replace the Projection Algorithm with an algorithm requiring at most $O(N\log N)$ operations and $O(N)$ storage, similar to the Approximate Selection Algorithm described in \cite{FP-2019}?} (A question posed by C. Fefferman.)
\smsk
\par We note that the answer to this question remains unclear even in the seemingly simpler one dimensional setting $D=1$, i.e., for set-valued mappings whose images are closed intervals in $\R$.\rbx}
\end{remark}

\begin{remark} {\em We now restate Theorem \reff{PA-PG-I} in terms of the quantity $|F|_{\Mf}$ defined in \rf{FM}. Given a constant $\lambda>0$, the Projection Algorithm yields one of two outcomes:
\par {\bf (``No Go''):} The algorithm confirms that
$|F|_{\Mf}\ge\lambda$, indicating no Lipschitz selection exists with seminorm strictly less than $\lambda$.
\par {\bf (``Success''):} The algorithm produces a Lipschitz selection $\tf$ of $F$ with seminorm
$$
|F|_{\Mf}\le \|\tf\|_{\Lip(\Mc)}\le 3\lambda.
$$
\par Thus, Theorem \reff{PA-PG-I} provides an efficient method to obtain lower and upper bounds for $|F|_{\Mf}$, cor\-responding to the {\bf (``No Go'')} and {\bf (``Success'')} outcomes, respectively.
\par However, this result does not fully address Problem \ref{LS-PRB}, which seeks an efficient algorithm that, given a finite pseudometric space $\Mf=\MR$ and a set-valued mapping $F:\Mc\to\HPL$, computes a Lipschitz selection $f$ of $F$ with seminorm at most
$\gamma |F|_{\Mf}$, where $\gamma>0$ is an absolute constant.
\par If we could compute $|F|_{\Mf}$ up to an absolute constant $\gamma$, and then apply Theorem \reff{PA-PG-I} with $$\lambda=\gamma\,|F|_{\Mf},$$
Problem \ref{LS-PRB} would be resolved. In this case, the Projection Algorithm would yield the outcome {\bf (``Success'')} and produce a Lipschitz selection $\tf$ of $F$ with seminorm $$\|\tf\|_{\Lip(\Mc)}\le 3\gamma\,|F|_{\Mf},$$ as required.
\par We explore the problem of efficient algorithms for computing $|F|_{\Mf}$ and related issues in the forthcoming paper \cite{S-2025}. Preliminary results are presented in
\cite{S-2023}, Section 6.3, where we introduce an algorithm that, given an $N$-element pseudometric space $\Mf=\MR$ and a set-valued mapping $F:\Mc\to\HPL$, computes $|F|_{\Mf}$ (up to an absolute constant factor) with a computational complexity (work and storage) of at most $C\,N^3$ for some universal constant $C$.\rbx}
\end{remark}
\smsk
\par Problem \reff{LS-PRB} is a special case of the general Lipschitz selection problem, which studies the existence and properties of Lipschitz selections of set-valued mappings from (pseudo)metric spaces into various families of convex subsets of Banach spaces. This problem can be viewed as a search for a Lipschitz mapping that approximately agrees with the data.
\smsk
\par There is an extensive literature devoted to various aspects of the Lipschitz selection and related smooth selection problems. For a comprehensive overview of known results and methods, we refer the reader to \cite{AF-1990,BL-2000,FIL-2016,FP-2019,FS-2018,
JLO-2022-1,JLO-2022-2,JLLL-2023,
PR-1992,PY-1989,PY-1995,S-2001,S-2002,S-2004,S-2008} and references therein.
\smsk
\par The Lipschitz selection problem is of great interest in recent years, mainly because of its deep connection with the classical {\it Whitney Extension Problem} \cite{Wh-1934}:
\smsk

\par {\it Given a positive integer $m$ and a function $f$ deﬁned on a closed subset of $\RN$, how can one determine whether $f$ extends to a $C^m$-function on all of $\RN$?}
\smsk
\par Over the years (since 1934) this problem has attracted a lot of attention, and there is an extensive literature devoted to this problem and its analogues for various spaces of smooth functions. For a detailed account of the history of extension and restriction problems for $m$-smooth functions, as well as various references related to this topic, we refer the reader to
\cite{BS-1994,BS-2001,F-2005,F-2009,FK-2009,F-2013,
FI-2020,JL-2021,FJL-2023,S-2008} and references therein.
\smsk
\par As an example, let us illustrate the connection between the Lipschitz selection problem and the Whitney problem for the space $C^2(\RN)$.
In \cite{S-1987,BS-2001,S-2002} we showed that the Whitney problem for the restrictions of $C^2$-functions to {\it finite} subsets of $\RN$ can be reduced to a certain Lipschitz selection problem for {\it affine-set valued} mappings. The solution to this special Lipschitz selection problem, presented in \cite{S-2001,S-2002,S-2004}, revealed an intriguing characteristic of $C^2$-function restrictions: the {\it Finiteness Principle}. This principle, named by C. Fefferman \cite{F-2009} for the general $C^m$-case, allows us to reduce the Whitney problem for $C^2(\RN)$-restrictions to {\it arbitrary finite subsets} of $\RN$ to a similar problem but for {\it finite sets containing at most $k^{\#}=3\cdot 2^{n-1}$ points} (see \cite{S-1987}). C. Fefferman later established in \cite{F-2005} that the Finiteness Principle holds for $C^m(\RN)$ for any $m, n \ge 1$, with a constant $k^{\#}=k^{\#}(m,n)$ depending only on $m$ and $n$.
\smsk
\par Furthermore, in \cite{S-2002} we solved the Lipschitz selection problem for the special case of the {\it line-set valued} mappings in $\RT$ and showed how constructive geometric criteria for Lipschitz selections of such mappings transform into purely analytical descriptions of the restrictions of $C^2$-functions to finite subsets of the plane.
\msk
\par These close connections between these two problems, the geometric Lipschitz selection problem and the analytic Whitney extension problem, allow us to transform efficient Lipschitz selection algorithms into efficient algorithms for various extension problems for spaces of smooth functions. Note also that efficient algorithms for solving  $C^m$-interpolation problems were developed by C. Fefferman and B. Klartag in \cite{F-2009-2,FK-2009,FK-2009-2}.

\bsk\msk
\par {\bf 1.2 The Projection Algorithm: main steps and main properties.}
\addtocontents{toc}{~~~~1.2 The Projection Algorithm: main steps and main properties.\hfill \thepage\par\VST}
\msk
\par In this section we present a geometric description of the general version of the Projection Algorithm for the family $\CRT$ of all convex closed subsets of $\RT$.
\par The input data of this algorithm are:
\par ($\blbig 1$) Two real numbers $\lambda_1,\lambda_2\ge 0$;\smsk
\par ($\blbig 2$) An $N$-point pseudometric space $\MR$;
\smsk
\par ($\blbig 3$) A set-valued mapping  $F:\Mc\to\CRT$.
\msk
\par Let us prepare the ingredients that are needed for this geometric description.
\begin{definition}\lbl{M-REF} {\em Given a constant $\lambda\ge 0$ we define a set-valued mapping $F^{[1]}[\cdot:\lambda;\rho]$ on $\Mc$ by letting
\bel{F-1}
F^{[1]}[x:\lambda;\rho]=\bigcap_{y\in \Mc}\,
\left[F(y)+\lambda\,\rho(x,y)\,Q_0\right],~~~~x\in\Mc.
\ee
\par We refer to the mapping $F^{[1]}[\cdot:\lambda;\rho]$ as {\it the $(\lambda;\rho)$-metric refinement of $F$}.
If the pseudometric $\rho$ is clear from the context, we omit $\rho$ in this notation and call the mapping in \rf{F-1} the $\lambda$-metric refinement of $F$.}
\end{definition}
\par Given a convex set $S\subset\RT$, we let $\HR[S]$ denote  {\it the smallest closed rectangle (possibly unbounded) with sides parallel to the coordinate axes, containing} $S$. We refer to $\HR[S]$ as {\it the rectangular hull} of the set $S$. See Fig. 1.

\begin{figure}[H]
\hspace*{27mm}
\includegraphics[scale=0.65]{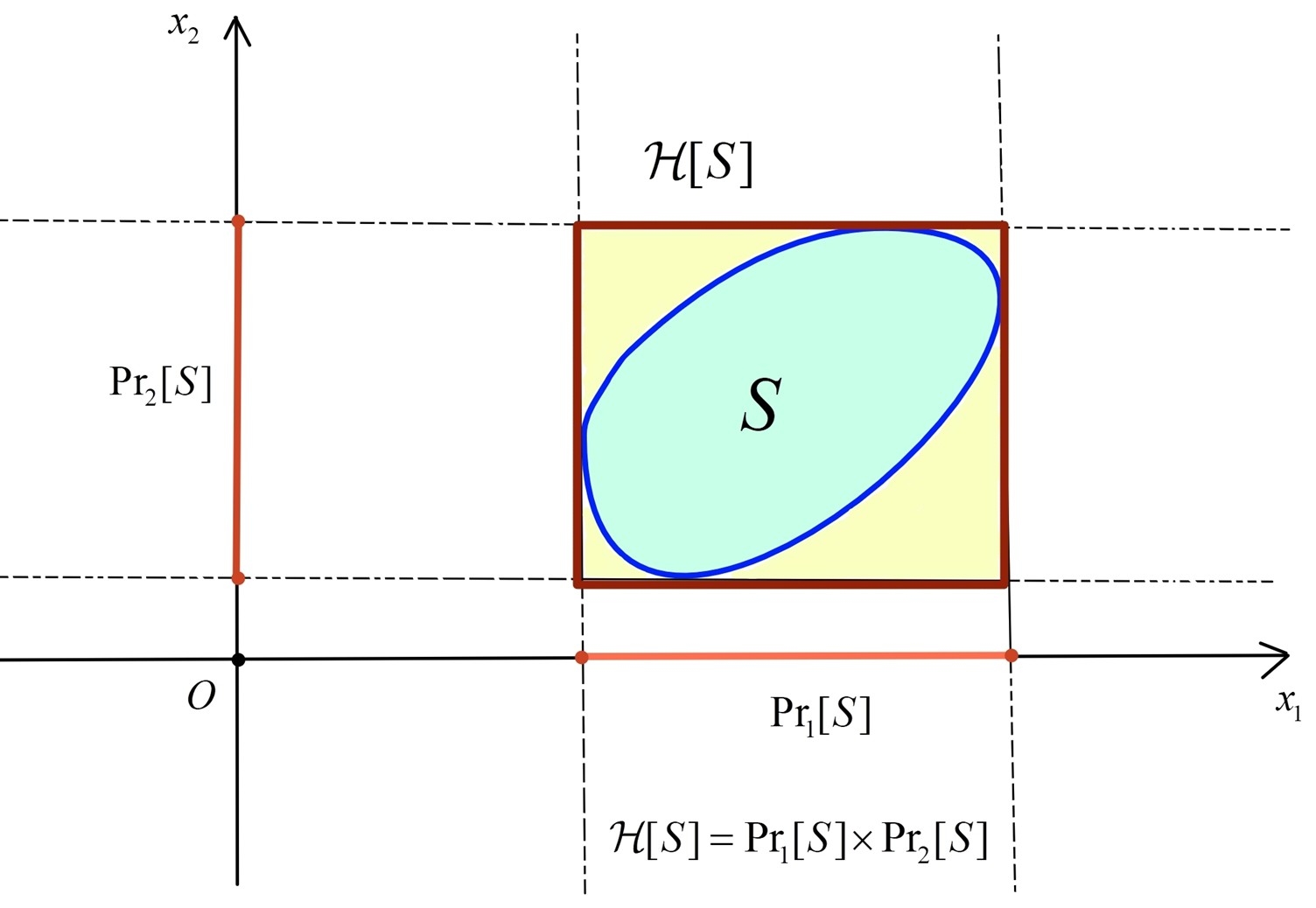}
\caption{The rectangular hull of a set $S$.}
\end{figure}

\par By $\Prm(\cdot,S)$ we denote the operator of metric projection onto a closed convex subset $S\subset\RT$. Finally, we let $\cent\,(\cdot)$ denote the center of a centrally symmetric bounded set in $\RT$.
\smsk
\begin{pra} {\em Given a vector $\lmv=(\lambda_1,\lambda_2)$ with non-negative coordinates $\lambda_1$ and $\lambda_2$, a finite pseudometric space $\Mf=\MR$, and a set-valued mapping $F:\Mc\to\CRT$, the Projection Algorithm which we define below either produces a selection $f_{\lmv;F}$ of $F$ (the outcome {\bf ``Success''}) or terminates (the outcome {\bf ``No go''}). This procedure includes the following five main steps.
\bsk
\par {\bf STEP 1.} At this step we construct the $\lambda_1$-metric refinement of $F$, i.e., the set-valued mapping
\bel{F-LM1}
F^{[1]}[x:\lambda_1]=\bigcap_{y\in \Mc}\,
\left[F(y)+\lambda_1\,\rho(x,y)\,Q_0\right].
\ee
\par If $F^{[1]}[x:\lambda_1]=\emp$ for some $x\in\Mc$, the algorithm produces the outcome {\bf ``No go''} and terminates.
\bsk
\par {\bf STEP 2.} Let us assume that the above condition does not hold, i.e., {\it for every element $x\in\Mc$ the $\lambda_1$-metric refinement $F^{[1]}[x:\lambda_1]$ is nonempty}.
\smsk
\par In this case, for each $x\in\Mc$, we construct {\it the rectangular hull} of $F^{[1]}[x:\lambda_1]$, the set
\bel{TC-DF}
\Tc_{F,\lambda_1}(x)=\Hc[F^{[1]}[x:\lambda_1]].
\ee
\par Thus, $\Tc_{F,\lambda_1}$ maps $\Mc$ into the family $\RCT$ of all rectangles in $\RT$.
\par See Fig. 2.
\bsk
\begin{figure}[h!]
\hspace{27mm}
\includegraphics[scale=0.7]{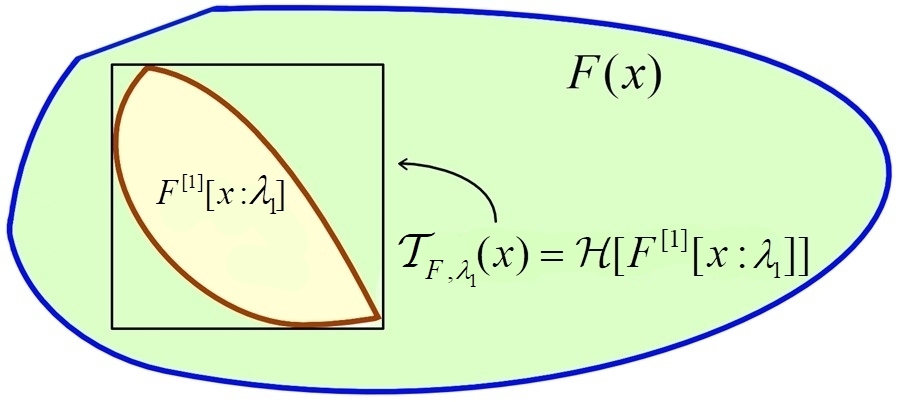}
\caption{The second step of the Projection Algorithm.}
\end{figure}
%
\par {\bf STEP 3.} For every $x\in\Mc$, we construct the $\lambda_2$-metric refinement of the mapping $\Tc_{F,\lambda_1}$, i.e., the rectangle $\Tc^{[1]}_{F,\lambda_1}[x:\lambda_2]$ defined by
\bel{RT-1}
\Tc^{[1]}_{F,\lambda_1}[x:\lambda_2]=\bigcap_{y\in \Mc}\,
\left[\Tc_{F,\lambda_1}(y)+\lambda_2\,\rho(x,y)\,Q_0\right].
\ee
See Fig. 3.

\begin{figure}[H]
\hspace{28mm}
\includegraphics[scale=0.72]{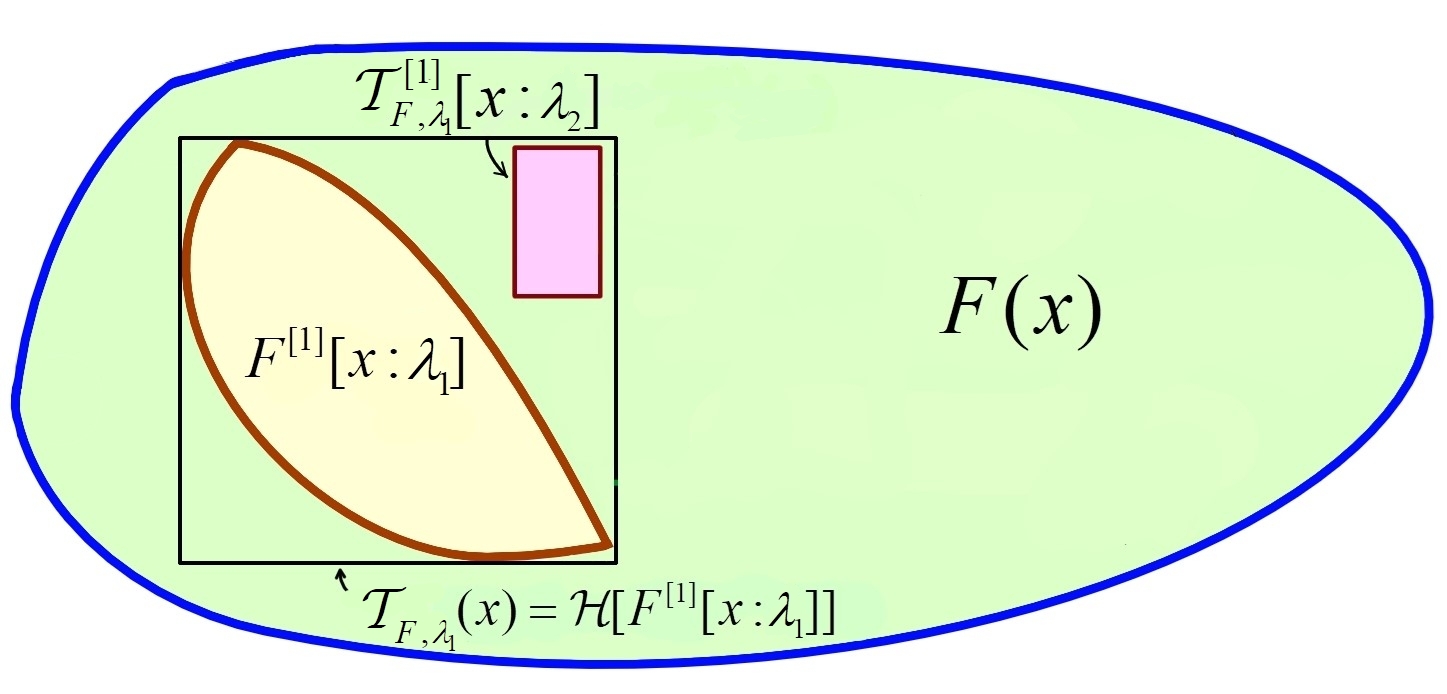}
\vspace*{3mm}
\caption{The third step of the Projection Algorithm.}
\end{figure}
\par If
$$
\Tc^{[1]}_{F,\lambda_1}[x:\lambda_2]=\emp~~~~\text{for some}~~~x\in\Mc,
$$
the algorithm produces the outcome {\bf ``No go''} and terminates.
\bsk\smsk
\par {\bf STEP 4.} At this step, we assume that {\it for each $x\in\Mc$ the rectangle $\Tc^{[1]}_{F,\lambda_1}[x:\lambda_2]\ne\emp$}. Let $O=(0,0)$ be the origin.
\smsk
\par We define a mapping $g_F:\Mc\to\RT$ by letting
\smsk
\bel{G-CNT}
g_F(x)=\cent\left
(\Prm\left(O,\Tc^{[1]}_{F,\lambda_1}[x:\lambda_2]\right) \right), ~~~~~~x\in\Mc.
\hspace{18mm}
\ee
\par See Fig. 4.
\vspace*{5mm}
\begin{figure}[h!]
\hspace{25mm}
\includegraphics[scale=0.76]{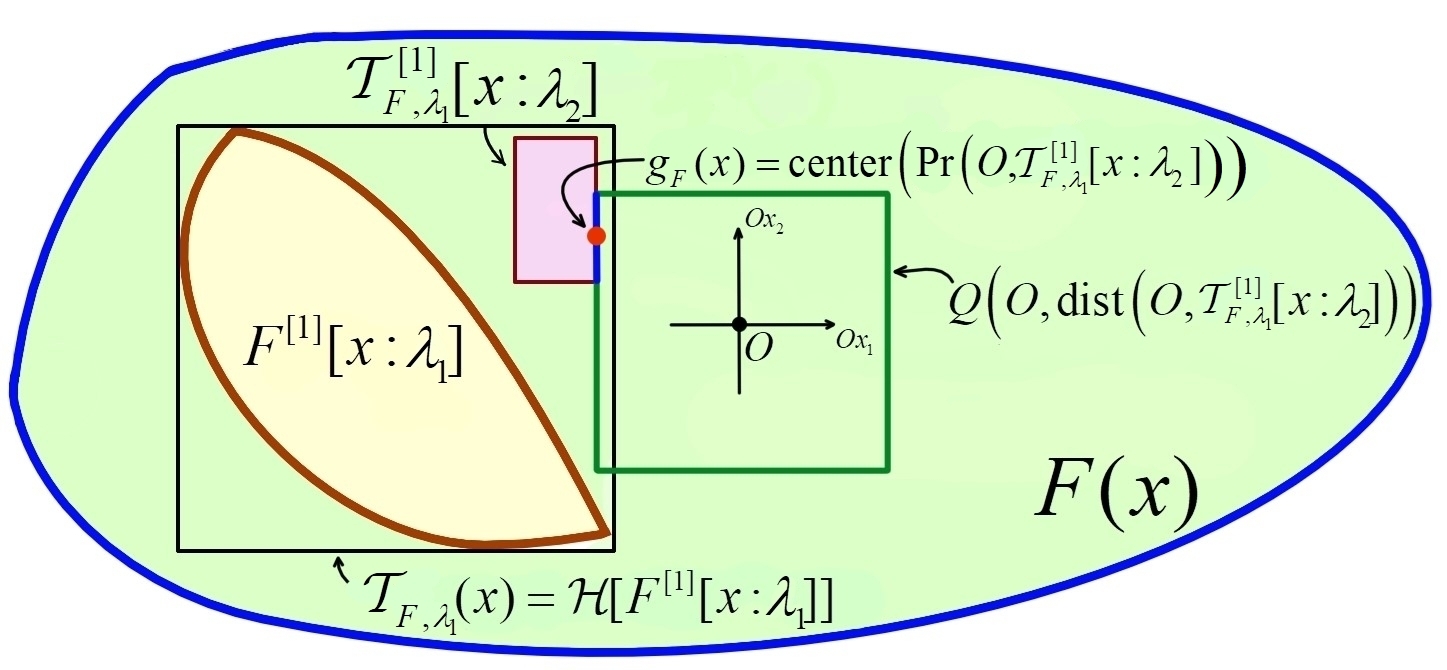}
\vspace*{3mm}
\caption{The mapping $g_F:\Mc\to\RT$, a selection of the set-valued mapping $\Tc_{F,\lambda_1}$.}
\end{figure}
\msk\bsk
\par {\bf STEP 5.} We define the mapping $f_{\lmv;F}:\Mc\to\RT$ by letting\msk
\bel{F-LF}
f_{\lmv;F}(x)=\Prm(g_F(x),F^{[1]}[x:\lambda_1]), ~~~~~~x\in\Mc.
\ee
\par See Fig. 5.
\begin{figure}[h]
\hspace{25mm}
\includegraphics[scale=1.18]{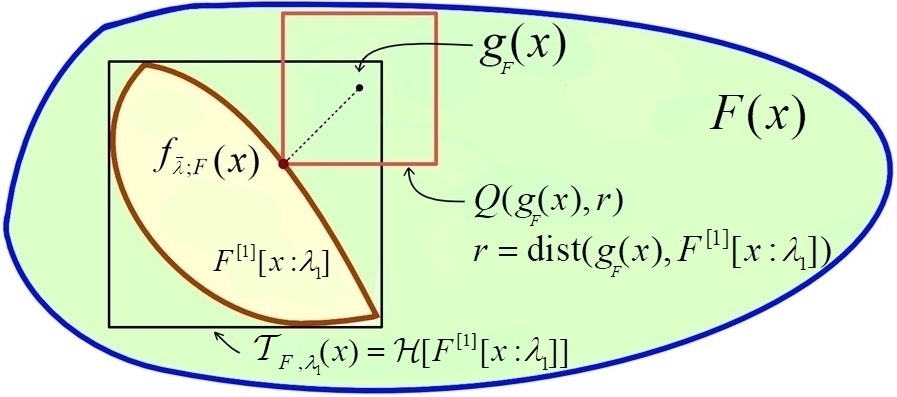}
\caption{The final step of the Projection Algorithm.}
\end{figure}
\smsk
\par At this stage, the algorithm produces the outcome
{\bf ``Success''} and terminates.\bx}
\end{pra}
\par To specify the dependence on the parameters $\lambda_1$ and $\lambda_2$, we call the above algorithm
$$
\text{{\it the $\lmv$-Projection Algorithm}}~~~~~~
\text{($\lmv$-PA for short)}
$$
or $(\lambda_1,\lambda_2)$-Projection Algorithm ($(\lambda_1,\lambda_2)$-PA for short).

\begin{remark} {\em Let us note the following important properties of the $(\lambda_1,\lambda_2)$-Projection Algorithm which are immediate from the construction of this algorithm presented above:
\smsk
\par (i) The $(\lambda_1,\lambda_2)$-PA produces the outcome {\bf ``No Go''} if and only if
\bel{NG-O}
\text{there exists}~~ x\in\Mc~~\text{such that either} ~~F^{[1]}[x:\lambda_1]=\emp~~\text{or}~~
\Tc^{[1]}_{F,\lambda_1}[x:\lambda_2]=\emp;
\ee
\par (ii) The $(\lambda_1,\lambda_2)$-PA produces the outcome {\bf ``Success''} if and only if
\bel{SC-P}
F^{[1]}[x:\lambda_1]\neq\emp~~~\text{and}~~~
\Tc^{[1]}_{F,\lambda_1}[x:\lambda_2]\neq\emp~~~~
\text{for every}~~~x\in\Mc.
\ee
}
\end{remark}
\par The next theorem describes the main properties of
the $\lmv$-Projection Algorithm.
\begin{theorem}\lbl{ALG-T} Let $\lambda_1,\lambda_2\ge 0$, and let $\lmv=(\lambda_1,\lambda_2)$. Let $\Mf=\MR$ be a finite pseudometric space, and let $F:\Mc\to\CRT$ be a set-valued mapping.
\smsk
\par (i) If $\lmv$-Projection Algorithm produces the outcome {\bf ``No go''}, then we gua\-ran\-tee that there does not exist a Lipschitz selection of $F$ with Lipschitz seminorm at most $\min\{\lambda_1,\lambda_2\}$.
\smsk
\par  (ii) Suppose that the $\lmv$-Projection Algorithm produces the outcome {\bf ``Success''}. In this case, this algorithm returns the mapping $f_{\lmv;F}:\Mc\to\RT$, see \rf{F-LF}, with the following properties:
\smsk
\par ($\bigstar A$) The mapping $f_{\lmv;F}$ is well defined. This means the following:
\par (a) For every $x\in\Mc$, the set $F^{[1]}[x:\lambda_1]$ and the rectangle $\Tc^{[1]}_{F,\lambda_1}[x:\lambda_2]$ are nonempty;
\par (b) The mapping $g_F$ (see \rf{G-CNT}) is well defined;
\par (c) The metric projection defined by the right hand side of \rf{F-LF} is a {\it singleton}.
\smsk
\par ($\bigstar B$) $f_{\lmv;F}$ is a Lipschitz selection of $F$ with Lipschitz seminorm
$$
\|f_{\lmv;F}\|_{\Lip(\Mc)}\le \lambda_1+2\lambda_2.
$$
\end{theorem}
\par In Section 5 we show that these properties of
the $\lmv$-Projection Algorithm are immediate from a stronger result, Theorem \reff{W-CR}, which we formulate below.
\par Given $\lambda\ge 0$ and $x,x',x''\in\Mc$, we let   $\Wc_F[x,x',x'':\lambda]$ denote a subset of $\RT$ defined by
\msk
$$
\Wc_F[x,x',x'':\lambda]=
\HR[\{F(x')+\lambda\,\rho(x',x)\,Q_0\}
\cap \{F(x'')+\lambda\,\rho(x'',x)\,Q_0\}].
$$
See Fig. 6.
\smsk
\begin{figure}[H]
\hspace{4mm}
\includegraphics[scale=1.03]{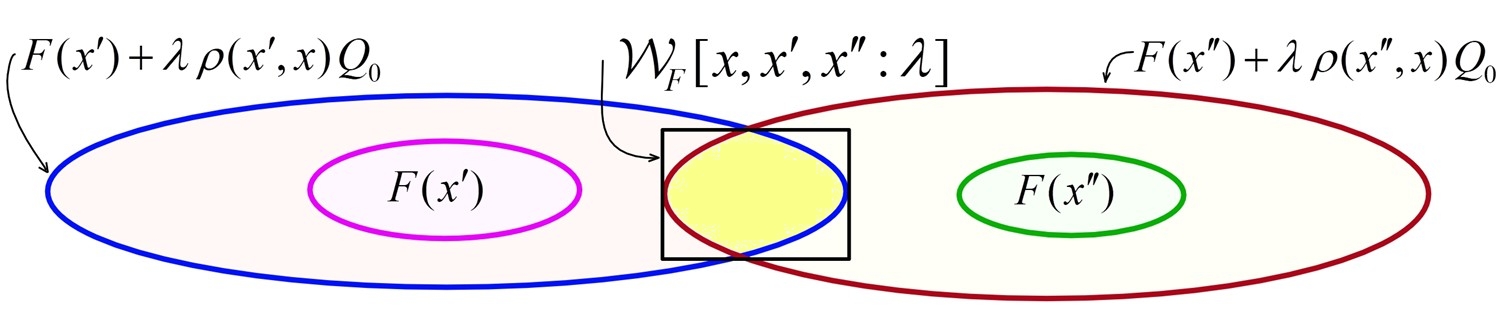}
\caption{The set ${\mathcal W}_F[x,x',x'':\lambda]$.}
\end{figure}
\smsk
\par We recall that by $\HR[\cdot]$ we denote the rectangular hull of a set.

\begin{theorem}\lbl{W-CR} Let $\Mf=\MR$ be a finite pseudometric space, and let $F:\Mc\to\CRT$ be a set-valued mapping. Given non-negative constants $\tlm$ and $\lambda$, let us assume that for every $x,x',x'',y,y',y''\in\Mc$ the following condition
\bel{WNEW}
\Wc_F[x,x',x'':\tlm]\cap \{\Wc_F[y,y',y'':\tlm]+\lambda\,\rho(x,y)\,Q_0\}\ne\emp
\ee
holds. Then $F$ has a Lipschitz selection with Lipschitz seminorm at most $2\lambda+\tlm$.
\end{theorem}

\begin{remark} {\em Note that condition \rf{WNEW} is equivalent to the following property:
\par $\Wc_F[x,x',x'':\tlm]\ne\emp$ and
$\Wc_F[y,y',y'':\tlm]\ne\emp$, and
$$
\dist\left(\Wc_F[x,x',x'':\tlm], \Wc_F[y,y',y'':\tlm]\right)\le\lambda\,\rho(x,y)
$$
for every $x,x',x'',y,y',y''\in\Mc$. See Fig. 7.~ \rbx}
\end{remark}

\begin{figure}[H]
\hspace{-5mm}
\includegraphics[scale=1.05]{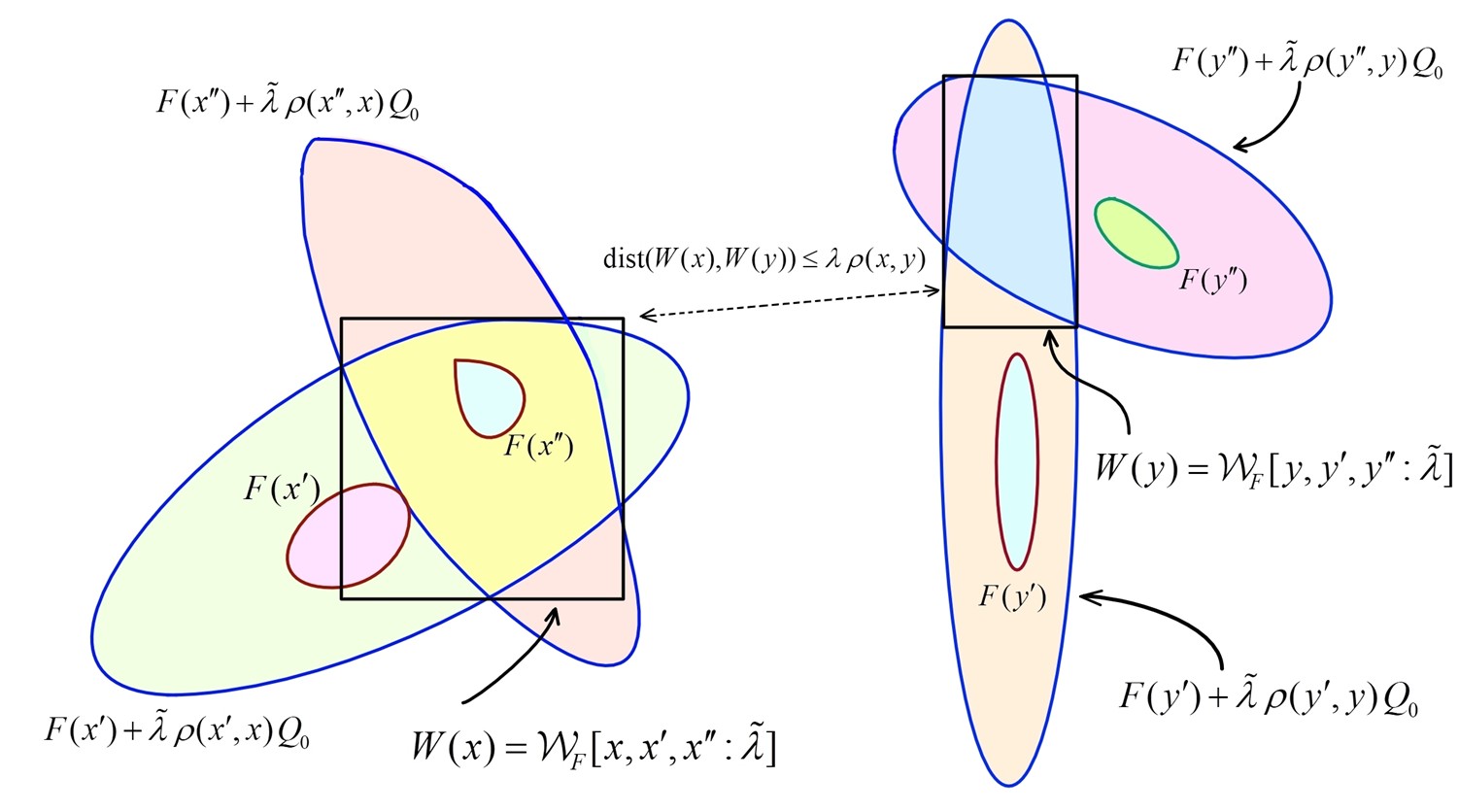}
\caption{The hypothesis of Theorem 1.9.}
\end{figure}

\par We refer to Theorem \reff{W-CR} as {\it the key theorem.} This theorem is the most technically difficult part of the present paper. We prove Theorem \reff{W-CR} in Sections 3 and 4, and Theorem \reff{ALG-T} in Section 5.
\par In Section 5 we also show that these theorems imply the following result.
\smsk
\begin{theorem}\lbl{PA-W} Let $\Mf=\MR$ be a finite pseudometric space, and let $F:\Mc\to\CRT$ be a set-valued mapping. Given $\tlm,\lambda\ge 0$, let us assume that  condition \rf{WNEW} holds for every $x,x',x'',y,y',y''\in\Mc$. Then the $\lmv$-PA with $\lmv=(\tlm,\lambda)$ produces the outcome {\bf ``Success''} and returns the mapping $f_{\lmv;F}:\Mc\to\RT$ which has the following properties:
\bel{FL-LS}
f_{\lmv;F}~~~\text{is a Lipschitz selection of}~~~F~~\text{with}~~~\|f_{\lmv;F}\|_{\Lip(\Mc)}\le 2\lambda+\tlm.
\ee
\end{theorem}
\par In other words, Theorem \reff{PA-W} states that the mapping $f_{\lmv;F}$ provides a selection of the set-valued mapping $F$ mentioned in the statement of Theorem \reff{W-CR}.
\smsk

\par We return to the $(\lambda_1,\lambda_2)$-Projection Algorithm for the family $\HPL$ of all closed half-planes in $\RT$. {\bf STEPS 1-5} contain a geometric description of this algorithm. However, this description is incomplete because we do not specify the algorithms for constructing several geometric objects that occur in {\bf STEPS 2-5}. Given $x\in\Mc$, such objects include:
\msk
\par {\bf STEP 1}: The set $F^{[1]}[x:\lambda_1]$, i.e., the $\lambda_1$-metric refinement of $F$. See \rf{F-LM1}.
\smsk
\par {\bf STEP 2}: The rectangle $\Tc_{F,\lambda_1}(x)$, see \rf{TC-DF}, i.e., the rectangular hull of the set $F^{[1]}[x:\lambda_1]$.
\smsk
\par {\bf STEP 3}: The rectangle $\Tc^{[1]}_{F,\lambda_1}[x:\lambda_2]$, see \rf{RT-1}, i.e., the $\lambda_2$-metric refinement of $\Tc_{F,\lambda_1}$.
\smsk
\par {\bf STEP 4}: The point $g_F(x)$, the center of the metric projection of the origin onto  $\Tc^{[1]}_{F,\lambda_1}[x:\lambda_2]$.
\par See \rf{G-CNT}.
\smsk
\par {\bf STEP 5}: The point $f_{\lmv;F}(x)$ defined by formula \rf{F-LF}, i.e., the metric projection of $g_F(x)$ onto the set $F^{[1]}[x:\lambda_1]$.
\smsk
\par The corresponding computational algorithms for these objects are presented in Sections 6.1 - 6.5. In these  sections, we also show that for each $x\in\Mc$,
the work and storage required to execute each of these algorithm are linear in $N$ (up to an absolute constant factor).
\par Algorithms of this kind are more or less elementary for the objects in {\bf STEPS 3-5}. However, at {\bf STEP 1} and {\bf STEP 2} we encounter classical problems of linear optimization for finite families of half-planes in $\RT$. In particular, at {\bf STEP 2} this is the problem of an efficient algorithm for computing the rectangle hull of a polygon determined by $N$ linear constraints. In this case, to solve these problems, we use classical results on low-dimensional linear programming by N. Megiddo \cite{M-1983-2} and M. E. Dyer \cite{Dy-1984}, which provide the required efficient algorithms with $O(N)$ running time and $O(N)$ units of computer memory.
\par Finally, by successively applying all these algorithms to our problem, we obtain the required algorithm for Lipschitz selection satisfying the conditions of Theorem \reff{PA-PG-I}.
\smsk
\par We complete the paper with Section 7 devoted to some generalizations of Theorem \reff{PA-PG-I}. First of them relates to wider families of convex subsets in $\RT$ than the family $\HPL$ of all closed half-planes. We also discuss efficiency of a polynomial algorithm for Lipschitz selection which provides the same outcomes as in Theorem \reff{PA-PG-I}, but with the estimate $\|f\|_{\Lip(\Mc)}\le \lambda$ instead of inequality \rf{G-3}. Our last generalization relates to a variant of Problem
\reff{LS-PRB} where the metric $\rho$ is replaced with an arbitrary non-negative function defined on $\Mc\times\Mc$.
\smsk
\smsk
\par {\bf Acknowledgements.} I am very grateful to Charles Fefferman for stimulating discussions and valuable advice.

\SECT{2. Notation and preliminaries.}{2}
\addtocontents{toc}{2. Notation and preliminaries.\hfill \thepage\par\VST}

\indent\par {\bf 2.1 Background notation.}

\addtocontents{toc}{~~~~2.1 Background notation.\hfill \thepage\par\VST}

\msk
\par Let $A$ and $B$ be nonempty subsets of $\RT$. We let
$A+B=\{a+b: a\in A, b\in B\}$ denote the Minkowski sum of these sets. Given $\lambda\ge 0$, we set
$\lambda\,A =\{\lambda\,a: a\in A\}$.
\par We write
$$
\dist(A,B)=\inf\{\|a-b\|:a\in A,~b\in B\}
$$
to denote the distance between $A$ and $B$. For $x\in\RT$, we set $\dist(x,A)=\dist(\{x\},A)$. We put $\dist(\emp,A)=0$ provided $A$ is an arbitrary (possibly empty) subset of $\RT$.
\smsk
\par We let $\dhf(A,B)$ denote the Hausdorff distance between $A$ and $B$ in the Banach space $\LTI=(\RT,\|\cdot\|)$:
\bel{HD-DF}
\dhf(A,B)=\inf\{r>0: A+rQ_0\supset B,~B+rQ_0\supset A\}
\ee
where $Q_0=[-1,1]^2$.
\par Given $a,b\in\RT$, by $[a,b]$ we denote the line segment (or a point if $a=b$) in $\RT$ with the ends in $a$ and $b$. Thus,
\bel{LSEG}
[a,b]=\{x\in\RT:x=(1-t)\,a+t\,b, 0\le t\le 1\}.
\ee
\par We write $[x]_+=\max\{x,0\}$ for the positive part of the real $x$. We also use the natural convention that
\bel{INF-S}
\frac{0}{0}=0,~~\frac{a}{0}=+\infty~~\text{for}~~a>0,
~~a-b=0~~~\text{if}~~a=b=\pm\infty,~~\text{and}~~
(\pm\infty)-(\mp\infty)=\pm\infty.
\ee
\par If $S$ is a finite set, by $\#S$ we denote the number of elements of $S$.
\par We let $\ip{\cdot,\cdot}$ denote the standard inner product in $\RT$. By
$$
Ox_1=\{x=(t,0):t\in\R\}~~~\text{and}~~~  Ox_2=\{x=(0,t):t\in\R\}
$$
we denote the coordinate axes in $\RT$.
\par Let $\Ic(Ox_i)$, $i=1,2$, be the family of all nonempty convex closed subsets of the coordinate axis $Ox_i$. In other words, $\Ic(Ox_i)$ is the family of all nonempty closed intervals (bounded or unbounded) and points lying on the $Ox_i$ axis. Given $I\in\Ic(Ox_1)$, let
\bel{LR1-END}
\lend(I)=\inf\{s\in\R:(s,0)\in I\}~~~\text{and}~~~
\rend(I)=\sup\{s\in\R:(s,0)\in I\}.
\ee
We call the numbers $\lend(I)$ and $\rend(I)$  {\it the left end and right end of $I$}, respectively.
\par In the same way we define the left and right ends of an interval $I\in\Ic(Ox_2)$:
\bel{LR2-END}
\lend(I)=\inf\{t\in\R:(0,t)\in I\}~~~\text{and}~~~
\rend(I)=\sup\{t\in\R:(0,t)\in I\}.
\ee
Note that the numbers $\lend(I)$ and $\rend(I)$ my take the values $-\infty$ or $+\infty$ provided $I$ is unbounded.
\smsk
\par By $\Prj_i$, $i=1,2$, we denote the operator of the orthogonal projection onto the axis $Ox_i$;
thus, given $x=(x_1,x_2)\in\RT$, we have $\Prj_1[x]=(x_1,0)$ and $\Prj_2[x]=(0,x_2)$. Given $i=1,2$, and a convex closed $S\subset\RT$, we let $\Prj_i[S]$, denote the orthogonal projection of $S$ onto the axis $Ox_i$. Thus,
$$
\Prj_1[S]=\{(s,0):(s,t)\in S~\text{for some}~t\in\R \} ~~\text{and}~~
\Prj_2[S]=\{(0,t):(s,t)\in S~\text{for some}~s\in\R \}.
$$


\begin{remark}\lbl{N-ST} {\em Note that the orthogonal projections of a convex closed set $S\subset\RT$ onto the coordinate axes are determined by the numbers
$$
\sigma_1(S)=\inf\{s:(s,t)\in S\},~~~
\sigma_2(S)=\sup\{s:(s,t)\in S\},
$$
and
$$
\tau_1(S)=\inf\{t:(s,t)\in S\},~~~
\tau_2(S)=\sup\{t:(s,t)\in S\}.
$$
\par This means that
$$
\sigma_1(S)=\lend(\Prj_1[S]),~~
\sigma_2(S)=\rend(\Prj_1[S])
$$
and
$$
\tau_1(S)=\lend(\Prj_2[S])~~\text{and}~~
\tau_2(S)=\rend(\Prj_2[S]),
$$
i.e., $\sigma_i(S)$, $i=1,2$, are the left and right ends of $\Prj_1[S]$, and $\tau_i(S)$, $i=1,2$, are the left and right ends of $\Prj_2[S]$. See \rf{LR1-END} and \rf{LR2-END}. More specifically,
$$
\Prj_1[S]=\{(s,0):\sigma_1(S)\le s\le \sigma_2(S)\}
$$
provided $S$ is bounded, and
$$
\Prj_1[S]=\{(s,0):-\infty<s\le \sigma_2(S)\}~~\text{or}~~
\Prj_1[S]=\{(s,0):\sigma_1(S)\le s<+\infty\}
$$
provided
$$
\sigma_1(S)=-\infty,~~\sigma_2(S)<\infty,
~~~\text{or}~~~
-\infty<\sigma_1(S),~~\sigma_2(S)=+\infty,
$$
respectively.
\par Finally,
$$
\Prj_1[S]=Ox_1~~~\text{if}~~~\sigma_1(S)=-\infty~~~
\text{and}~~~\sigma_2(S)=+\infty.
$$
\par In the same fashion the numbers $\tau_1(S)$ and $\tau_2(S)$ determine the set $\Prj_2[S]$.\rbx}
\end{remark}

\par Given sets $A_i\subset Ox_i$, $i=1,2$, we let $A_1\times A_2$ denote a subset of $\RT$ defined by
\bel{A1TA2}
A_1\times A_2=\{a=(a_1,a_2)\in\RT: (a_1,0)\in A_1, (0,a_2)\in A_2\}.
\ee
\par Given $a\in\RT$ and $r>0$, we let $Q(a,r)$ denote the square with center $a$ and length of side $2r$:
$$
Q(a,r)=\{y\in\RT:\|y-a\|\le r\}.
$$
In particular,
$$
Q_0=[-1,1]^2=Q(0,1)
$$
is the unit ball of the Banach space $\LTI=(\RT,\|\cdot\|)$.
\smsk
\par Let $S$ be a nonempty {\it convex} closed subset of $\RT$.  By $\Prm(\cdot,S)$ we denote the operator of metric projection onto $S$ in $\LTI$-norm. To each $a\in\RT$ this operator assigns the set of all points in $S$ that are nearest to $a$ on $S$ in the uniform norm. Thus,
\bel{MPR}
\Prm(a,S)=S\cap Q(a,\dist(a,S)).
\ee
Clearly, the set $\Prm(a,S)$ is either a singleton or a line segment in $\RT$ parallel to one of the coordinate axes.
\par If $S\subset\RT$ is convex bounded and centrally symmetric, by $\cent(S)$ we denote the center of $S$.
\smsk
\par Given non-zero vectors $u,v\in\RT$ we write $u\parallel v$ if $u$ and $v$ are collinear, and we write $u\nparallel v$ whenever these vectors are non-collinear.
\smsk
\par Let $\tu$ and $\tv$ be the {\it unit} vectors in the direction of the vectors $u$ and $v$ respectively. By $\theta(u,v)\in[0,2\pi)$ we denote {the angle of rotation} from $\tu$ to $\tv$ in {\it the counterclockwise direction}. (Thus, $\theta(v,u)=2\pi-\theta(u,v)$.) We refer to $\theta(u,v)$ as the angle between the vectors $u$ and $v$.
\smsk
\par Let $\ell_1$ and $\ell_2$ be two non-parallel straight lines in $\RT$; in this case, we write $\ell_1\nparallel\ell_2$. Let $V=\ell_1\cap \ell_2$. These two lines form two angles $\vf_1,\vf_2\in[0,\pi)$, $\vf_1+\vf_2=\pi$, with the vertex at the point $V$. Let
$$
\vf(\ell_1,\ell_2)=\min\{\vf_1,\vf_2\};~~~~~
\text{clearly,}~~~~ \vf(\ell_1,\ell_2)\in[0,\pi/2].
$$
\par We refer to $\vf(\ell_1,\ell_2)$ as {\it ``the angle between the straight lines $\ell_1$ and $\ell_2$''.} If $\ell_1\parallel \ell_2$ (i.e., $\ell_1$ and $\ell_2$ are parallel), we set $\vf(\ell_1,\ell_2)=0$.
\bsk\bsk

\par {\bf 2.2 Rectangles and rectangular hulls.}

\addtocontents{toc}{~~~~2.2 Rectangles and rectangular hulls.\hfill \thepage\par\VST}

\msk
\par Recall that $\Ic(Ox_i)$, $i=1,2$, is the family of all nonempty convex closed subsets of the axis $Ox_i$. We set
\bel{RCT-D}
\RCT=\{\Pi=I_1\times I_2: I_1\in\Ic(Ox_1),I_2\in\Ic(Ox_2)\}.
\ee
We refer to every member of the family $\RCT$ as a {\it ``rectangle''}. Furthermore, throughout the paper, the word ``rectangle'' will mean an element of $\RCT$, i.e., a closed rectangle (possibly unbounded) with sides parallel to the coordinate axes.
\par Clearly, thanks to definition \rf{RCT-D},
\bel{PR-RW}
\Pi=\Prj_1[\Pi]\times \Prj_2[\Pi]~~~~\text{for every rectangle}~~~\Pi\in\RCT.
\ee
\par Because $\Prj_i$, $i=1,2$, is a continuous operator, for every {\it convex} closed set $S\subset\RT$ its orthogonal projection $\Prj_i[S]$ onto the axis $Ox_i$ is a {\it convex} closed subset of $Ox_i$, i.e., $\Prj_i[S]\in\Ic(Ox_i)$. Thus, thanks to this property and \rf{A1TA2}, we have the following:
\bel{RC-PD}
\text{a convex closed set}~~\Pi~~\text{belongs to}~~\RCT~~\text{if and only if}~~   \Pi=\Prj_1[\Pi]\times \Prj_2[\Pi].
\ee
\par Let us also note the following elementary properties of rectangles.
\begin{claim}\lbl{CLH-R} Let $\Rc\subset\RCT$ be a finite collection of rectangles in $\RT$.
\par (i) Suppose that the rectangles from the family $\Rc$ have a common point. Then the following equality holds:
\bel{RC-IP}
\bigcap_{\Pi\in\Rc} \Pi=
\left\{\bigcap_{\Pi\in\Rc}\Prj_1[\Pi]\right\}
\times \left\{\bigcap_{\Pi\in\Rc}\Prj_2[\Pi]\right\}.
\ee
\par In particular, for every $i=1,2$, we have
$$
\Prj_i\left[\bigcap_{\Pi\in\Rc} \Pi\right]=
\bigcap_{\Pi\in\Rc}\Prj_i[\Pi].
$$
\par (ii) The intersection of all rectangles of the family $\Rc$ is nonempty if and only if
$$
\text{both}~~~\bigcap_{\Pi\in\Rc}\Prj_1[\Pi]\ne\emp~~~
\text{and}~~~ \bigcap_{\Pi\in\Rc}\Prj_2[\Pi]\ne\emp.
$$
\end{claim}
\par The proof of this claim is obvious.
\smsk
\par We also note that any closed interval in $\RT$ lying on a line parallel to a coordinate axis is a ``rectangle''. In particular, every closed interval on the axis $Ox_1$ or $Ox_2$ belongs to the family $\RCT$.
\par Finally, given a bounded rectangle $\Pi\in\RCT$ we let $\cent(\Pi)$ denote the center of $\Pi$.
\smsk
\par We recall that $\CRT$ denotes the family of all nonempty convex closed subsets of $\RT$. Given a set $S\in\CRT$, we let $\HR[S]$ denote the {\it ``rectangular hull``} of $S$, i.e., the smallest (with respect to inclusion) rectangle containing $S$. Thus,
\bel{HRS}
\HR[S]=\cap\{\Pi: \Pi\in\RCT, \Pi\supset S\}.
\ee
\par Combining this definition with \rf{RC-IP}, we conclude that
\bel{P1XP2}
\HR[S]=\Prj_1[S]\times\Prj_2[S].
\ee
\par Thus, given $S\in\CRT$, its rectangular hull $\HR[S]$ is the only rectangle $\Pi$ for which
\bel{HS-U}
\Prj_1[\Pi]=\Prj_1[S]~~~~\text{and}~~~~~
\Prj_2[\Pi]=\Prj_2[S].
\ee

\par Let us also note the following elementary property of rectangles: for every closed convex set $S\subset\RT$ and every $r\ge 0$ we have $\HR[S+r Q_0]=\HR[S]+r Q_0.$
\msk
\par In this section we present two important auxiliary results. The first of them is a variant of the classical Helly's intersection theorem for rectangles.
\begin{lemma}\lbl{H-R} Let $\Kc\subset\RCT$ be a finite collection of rectangles in $\RT$. Suppose that the intersection of every two rectangles from $\Kc$ is nonempty. Then there exists a point in $\RT$ common to all of the family $\Kc$.
\end{lemma}
\par {\it Proof.} Representation \rf{PR-RW} reduces the problem to the one dimensional case. In this case the statement of the lemma is a variant of Helly's theorem in $\R$. See, e.g. \cite{DGK-1963}.\bx
\smsk
\par The second auxiliary result is a Helly-type theorem formulated in terms of the orthogonal projections onto the coordinate axes.
\begin{proposition}\lbl{INT-RE} Let $\Cf$ be a finite family of convex closed subsets in $\RT$ such that
\bel{P1-C}
\Prj_1[C_1\cap C_1']\,\cbg \Prj_1[C_2\cap C_2'] \ne\emp
\ee
for every $C_1,C_1',C_2,C_2'\in\Cf$. Then
\bel{IP-C}
\cap\{C:C\in\Cf\}\ne\emp.
\ee
\par Furthermore, in this case
\bel{PH-C}
\HR\left[\cap\{C:C\in\Cf\}\right]=
\cap\{\HR[C\cap C']: C,C'\in\Cf\}\,.
\ee
\end{proposition}
\par {\it Proof.}  Condition \rf{P1-C} tells us that for every $C,C'\in\Cf$ the set $C\cap C'$ is nonempty. Because $C\cap C'$ is a convex subset of $\RT$, its projection onto $Ox_1$, the set $\Prj_1[C\cap C']\subset Ox_1$, is convex as well, i.e., this set is a closed interval in $Ox_1$.
\par Let
$$
\Wc=\{\Prj_1[C\cap C']:C,C'\in\Cf\}.
$$
Then $\Wc$ is a finite family of intervals, and, thanks to \rf{P1-C}, every two members of this family have a common point. Helly's theorem tells us that in this case there exists a point in $Ox_1$ common to all of the family $\Wc$. See Lemma \reff{H-R}. Thus,
\bel{V-L}
V=\bigcap_{C,C'\in\,\Cf}\Prj_1[C\cap C']\ne\emp\,.
\ee
\par Fix a point $v\in V$. Then, thanks to \rf{V-L},
\bel{SV-LB}
v\in \Prj_1[C\cap C']~~~~\text{for every}~~~~C,C'\in\Cf.
\ee
\par Let
\bel{L-DFB}
L=\{w\in\RT:\Prj_1[w]=v\}
\ee
be the straight line through $v$ orthogonal to the axis $Ox_1$. See Fig. 8.

\begin{figure}[H]
\includegraphics[scale=0.75]{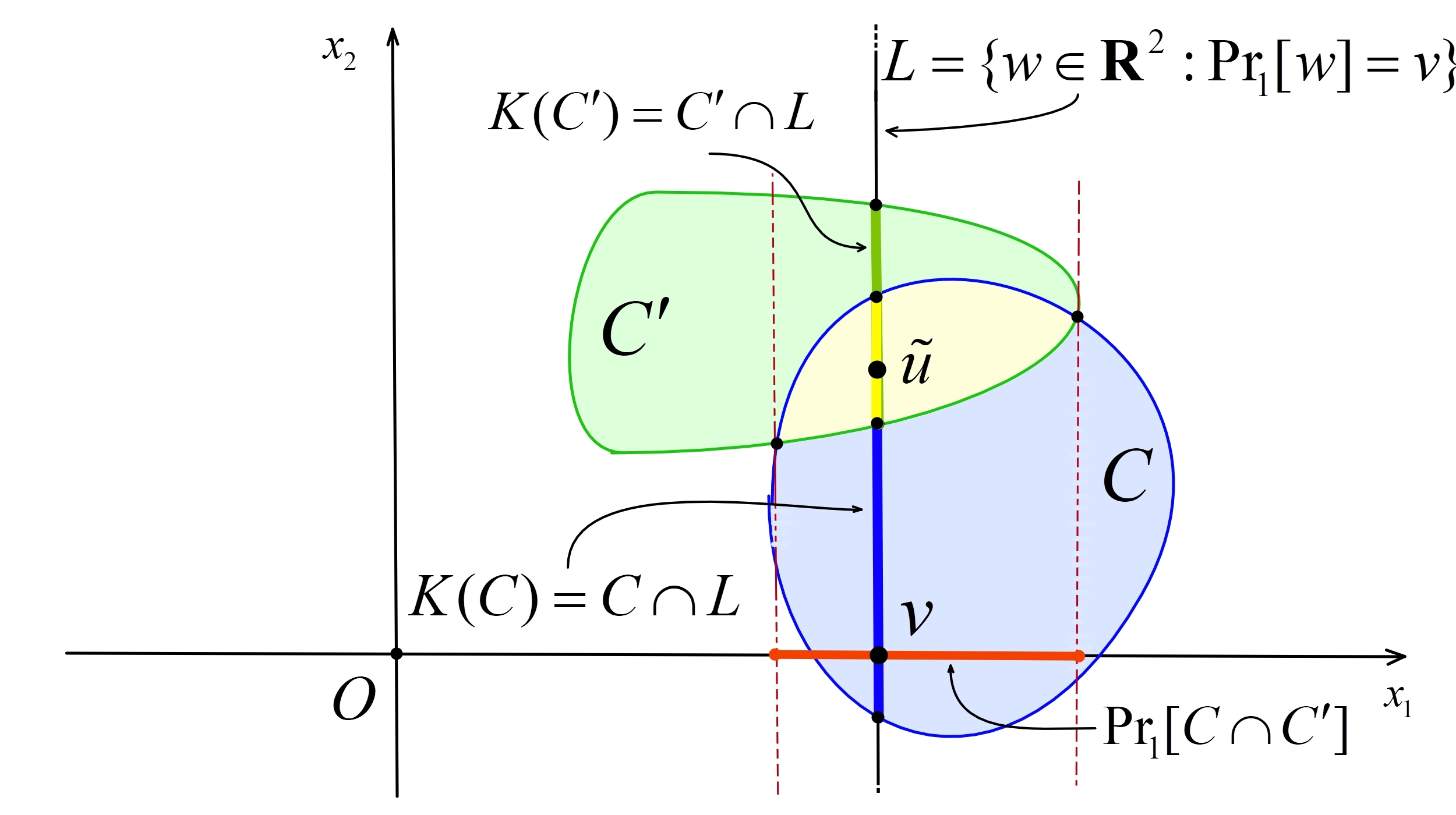}
\caption{Proposition 2.4.}
\end{figure}


\par Given $C\in\Cf$, we set $K(C)=C\cap L$. Thanks to \rf{SV-LB}, $v\in\Prj_1[C\cap C]=\Prj_1[C]$ so that there exists $u_C\in C$ such that $\Prj_1[u_C]=v$. From this and \rf{L-DFB}, we have $u_C\in C\cap L$ proving that $K(C)\ne\emp$ for every $C\in\Cf$. Clearly, each $K(C)$ is a closed interval lying on the straight line $L$. Let us show that there exists a point in $L$ common to all these intervals.
\smsk
\par Property \rf{SV-LB} tells us that for every $C,C'\in\Cf$ there exists a point $\tu\in C\cap C'$ such that $\Prj_1[\tu]=v$. Hence, thanks to \rf{L-DFB},
$\tu\in L$ so that
$$
\tu\in L\cap (C\cap C')=(L\cap C)\cap (L\cap C')=K(C)\cap K(C').
$$
\par This proves that any two members of the family $\Kc=\{K(C): C\in\Cf\}$ have a common point. Furthermore,  $\Kc$ is a {\it finite} (because $\Cf$ is finite) family of intervals lying in $L$. Helly's theorem tells us that in this case $\cap\{K(C):C\in\Cf\}\ne\emp$. Thus,
$$
\cap\{K(C):C\in\Cf\}=\cap\{L\cap C:C\in\Cf\}
=L\cap\left(\cap\{C:C\in\Cf\}\right)\ne\emp
$$
proving \rf{IP-C}.
\smsk
\par Let us prove \rf{PH-C}. Clearly, the left hand side of the equality \rf{PH-C} is contained in its right hand side.
Let us prove the converse statement.
\par Fix a point
\bel{U-FS}
u=(u_1,u_2)\in\cap \{\HR[C\cap C']: C,C'\in\Cf\}
\ee
and prove that $u\in\HR\left[\cap\{C:C\in\Cf\}\right]$.
Thanks to \rf{P1XP2}, property \rf{U-FS} is equivalent to the following one:
$u_1\in\Prj_1[\Tc]$ and $u_2\in\Prj_2[\Tc]$ where $\Tc=\cap\{C:C\in\Cf\}$.
\smsk
\par Prove that $u_1\in\Prj_1[\Tc]$. We let $L_1$ denote the straight line in $\RT$ through the point $u=(u_1,u_2)$ orthogonal to the axis $Ox_1$. Thus,
\bel{L1-OR}
L_1=\{w\in\RT:\Prj_1[w]=u_1\}.
\ee
See Fig. 9.

\begin{figure}[H]
\includegraphics[scale=0.75]{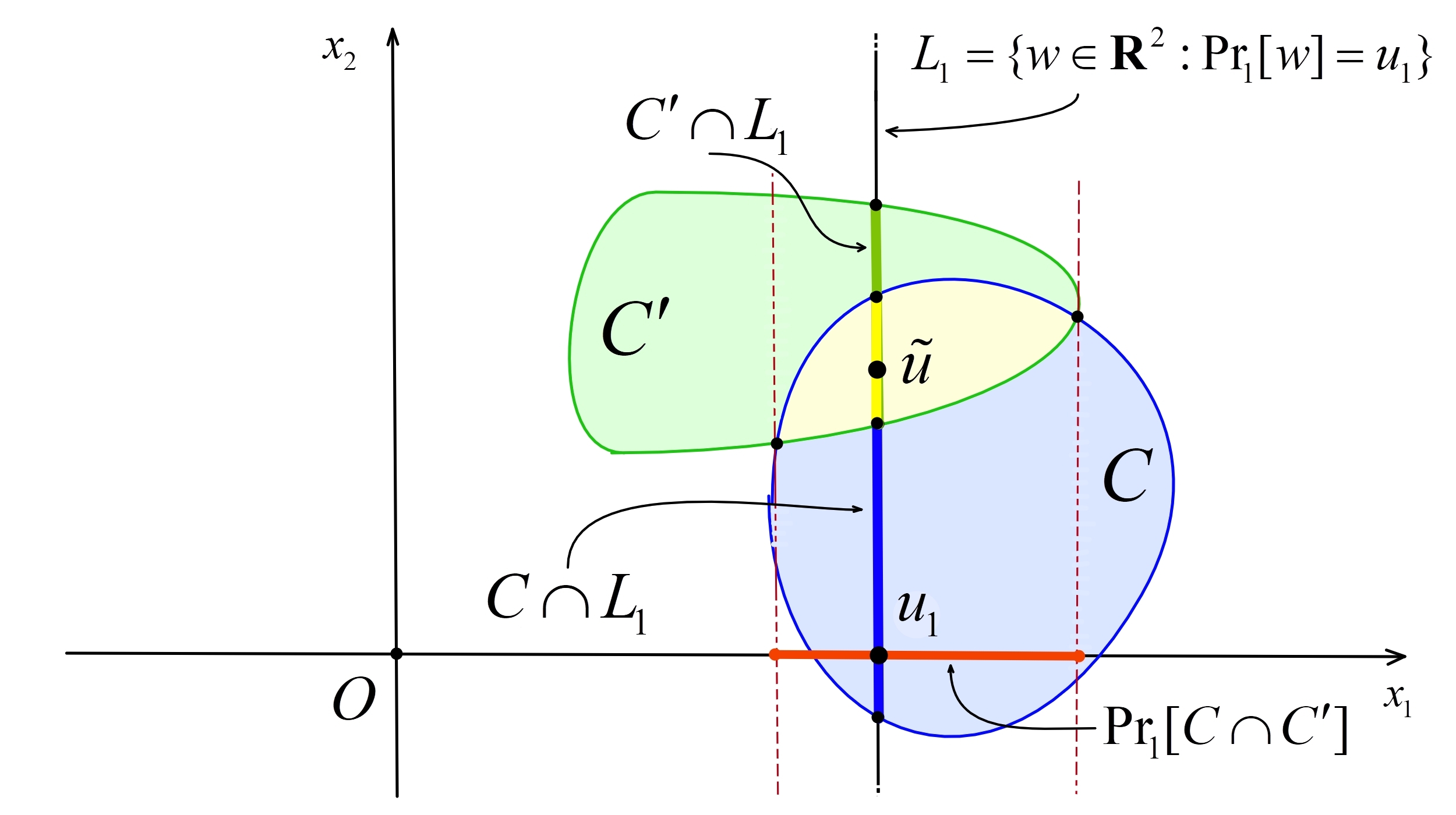}
\caption{Proposition 2.4.}
\end{figure}

\par Let us show that $L_1\cap\Tc\ne\emp$. Indeed, thanks to \rf{U-FS}, $u\in\HR[C\cap C']$ provided  $C,C'\in\Cf$ so that, thanks to \rf{P1XP2},
$u_1\in\Prj_1[C\cap C']$ for every $C,C'\in\Cf$.
Combining this property with definition \rf{L1-OR}, we conclude that $L_1\cap C\cap C'\ne\emp$ for all $C,C'\in\Cf$. Thus, every two members of the (finite) family $\Kc_1=\{L_1\cap C:C\in\Cf\}$ have a common point.
Therefore, thanks to Helly's theorem, there exists a point $\tu\in\RT$ such that
$$
\tu\in \cap\{L_1\cap C:C\in\Cf\}=L_1\cap\{C:C\in\Cf\}=L_1\cap\Tc.
$$
\par Thus, $\tu\in L_1$ so that, thanks to \rf{L1-OR}, $u_1=\Prj[\tu]$. Furthermore, $\tu\in \Tc=\cap\{C:C\in\Cf\}$ so that $u_1\in \Prj_1[\Tc]$.
\par In the same way we prove that $u_2\in\Prj_2[\Tc]$ completing the proof of the proposition.\bx
\smsk
\par We conclude the section with the following claim, which we use in Section 6.4 below.
\begin{claim}\lbl{MPR-C} Let $I_i\in\Ic(Ox_i)$, be a closed interval on the axis $Ox_i$, $i=1,2$, and let $\Tc=I_1\times I_2$ be a rectangle in $\RT$. Let $a_i=L(I_i)$ and $b_i=R(I_i)$, $i=1,2$, be the left and right ends of the interval $I_i$ respectively. (See \rf{LR1-END} and \rf{LR2-END}.)
\par Then the distance from the origin to $\Tc$ can be calculated as follows:
\bel{DS-Q}
\dist(O,\Tc)=\max\{[a_1]_+,[-b_1]_+,[a_2]_+,[-b_2]_+\}.
\ee
\par Furthermore, let
$$
g=(g_1,g_2)=\cent\left(\,\Prm\left(O,\Tc\right)\right)
$$
be the center of the metric projection of the origin onto the rectangle $\Tc$, see \rf{MPR}. Then
\bel{G-2T}
g_1=(L_1+R_1)/2~~~~\text{and}~~~~g_2=(L_2+R_2)/2
\ee
where
\bel{C-MP}
L_i=\max\{-\dist(O,\Tc),a_i\}~~~\text{and}~~~
R_i=\min\{\dist(O,\Tc),b_i\},~~~i=1,2.
\ee
\end{claim}
\par {\it Proof.} Clearly,
$$
\dist(O,\Tc)=\dist(O,I_1\times I_2)=
\max\{\dist(O,I_1),\dist(O,I_2)\}.
$$
But
$$
\dist(O,I_i)=\min\{|x_i|:a_i\le x\le b_i\}
=\max\{[a_i]_+,[-b_i]_+\}, ~~i=1,2,
$$
proving \rf{DS-Q}.
\par Let us prove \rf{G-2T}. To do this, it suffices to show that
\bel{FE-5}
\Prm\left(O,\Tc\right)=[L_1,R_1]\times [L_2,R_2].
\ee
which is equivalent to the following property:
$$
\Prj_i[\Prm\left(O,\Tc\right)]=[L_i,R_i],~~i=1,2.
$$
(Clearly, thanks to \rf{C-MP}, $-\infty<L_i\le R_i<+\infty$, $i=1,2$.)
\smsk
\par Recall that $\Prm\left(O,\Tc\right)=Q(O,r)\cap (I_1\times I_2)$ where $r=\dist(O,\Tc)$. (See \rf{MPR}.) Therefore,
$$
\Prj_i[\Prm\left(O,\Tc\right)]=
\Prj_i[Q(O,r)\cap \Tc]=\Prj_i[Q(O,r)]\cap \Prj_i[I_1\times I_2]=[-r,r]\cap I_i,~~~i=1,2,
$$
so that
$$
\Prj_i[\Prm\left(O,\Tc\right)]=[-r,r]\cap [a_i,b_i]=[\max\{-r,a_i\},\min\{r,b_i\}]=[L_i,R_i],~~~i=1,2,
$$
proving \rf{FE-5} and the claim.\bx

\bsk

\indent\par {\bf 2.3 Rectangles: intersections, neighborhoods and selections.}

\addtocontents{toc}{~~~~2.3 Rectangles: intersections, neighborhoods and selections.\hfill \thepage\par\VST}

\msk

\par In this section we present several criteria and several constructive formulae for the optimal Lipschitz selections of set-valued mappings taking values in the
family $\RCT$ of all closed rectangles in $\RT$ with sides parallel to the coordinate axes. See \rf{RCT-D}.
\par Let $\Ic(\R)$ be the family of all closed intervals and all points in $\R$, and let $I_0=[-1,1]$. Given $a\in\R$ and $r\ge 0$, we set $rI_0=[-r,r]$. We also recall that, given a bounded interval $I\in\Ic(\R)$, by $\cent(I)$ we denote the center of $I$.
\begin{lemma}\lbl{H-RB} Let $\Kc\subset\RCT$ be a family of rectangles in $\RT$ with nonempty intersection. Then for every $r\ge 0$ the following equality
\bel{H-SM}
\left(\,\bigcap_{K\in\,\Kc} K\right) +r Q_0
=\bigcap_{K\in\,\Kc}\,\left\{\,K+rQ_0\,\right\}
\ee
holds.
\end{lemma}
\par {\it Proof.} Obviously, the right hand side of \rf{H-SM} contains its left hand side. Let us prove that
\bel{H-SM-LR}
\left(\,\bigcap_{K\in\,\Kc} K\right) +r Q_0
\supset\bigcap_{K\in\,\Kc}\,\left\{\,K+rQ_0\,\right\}.
\ee
\par This inclusion is based on the following simple claim: Let $\Ic$ be a family of convex closed subsets of $\R$ (intervals) with nonempty intersection. Let  $K=[a,b]\subset\R$, be a closed bounded interval such that $K\cap I\ne\emp$ for every $I\in\Ic$. Then there exists a point common to $K$ and all of the members of the family $\Ic$. (The proof is immediate from Helly's theorem in $\R$ applied to the family $\Ic\cup\{K\}$ of closed intervals.)
\par This claim implies the following one dimensional variant of inclusion \rf{H-SM-LR}: Let $\Ic$ be a family of intervals in $\R$ with nonempty intersection. Then
\bel{H-SM-R1}
\left(\,\bigcap_{I\in\,\Ic} I\right) +r I_0
\supset\bigcap_{I\in\,\Ic}\,\left\{\,I+rI_0\,\right\}
~~~~\text{where}~~~~I_0=[-1,1].
\ee
\par Indeed, if $u\in\cap\{I+rI_0:I\in\Ic\}$ then $[u-r,u+r]\cap I\ne\emp$ for every $I\in\Ic$. Therefore, thanks to the above claim, $[u-r,u+r]\cap(\cap\{I:I\in\Ic\})\ne\emp$ proving that $u$ belongs to the left hand side of \rf{H-SM-R1}.
\smsk
\par Now, let us prove \rf{H-SM-LR} using \rf{H-SM-R1} and properties \rf{RC-PD} and \rf{RC-IP} of rectangles. For every $i=1,2$, we have
$$
\Prj_i\left[\left(\,\bigcap_{K\in\,\Kc} K\right)
+r Q_0\right]=
\left(\,\bigcap_{K\in\,\Kc} \Prj_i[K]\right)
+r\Prj_i[Q_0]=U_i.
$$
Furthermore,
$$
\Prj_i\left[
\bigcap_{K\in\,\Kc}\,\left\{\,K+rQ_0\,\right\}
\right]=
\bigcap_{K\in\,\Kc}\left\{\Prj_i[K]+r\Prj_i[Q_0]\right\}
=V_i.
$$
\par Thanks to inclusion \rf{H-SM-R1}, $U_i\supset V_i$, $i=1,2$, proving that the orthogonal projections onto the coordinate axes of the left hand side of \rf{H-SM-LR} contain the corresponding projections of its right hand side. Because the left and right hand sides of \rf{H-SM-LR} are {\it rectangles}, inclusion \rf{H-SM-LR} holds.
\par The proof of the lemma is complete.\bx

\smsk
\par We will also need the following simple claim.
\begin{claim}\lbl{TWO} (i) Let $A$ and $B$ be two closed intervals in $\R$. Then
\bel{AB-END}
\dhf(A,B)=\max\{|\inf A-\inf B|,|\sup A-\sup B|\}.
\ee
(See also our convention \rf{INF-S} for the cases of $\inf A,\inf B=-\infty$ and $\sup A,\sup B=+\infty$.)
\smsk
\par (ii) Let $\Ac,\Bc\in\RCT$ be two bounded rectangles in $\RT$. Then
\bel{CAB-H}
\|\cent(\Ac)-\cent(\Bc)\|\le\dhf(\Ac,\Bc).
\ee
\end{claim}
\par {\it Proof.} (i) Let
$$
r=\dhf(A,B)~~~\text{and}~~~
\delta=\max\{|\inf A-\inf B|,|\sup A-\sup B|\}.
$$
Then, thanks to \rf{HD-DF}, $A+rI_0\supset B$ proving that $\sup A+r\ge \sup B$ and $\inf A-r\le \inf B$. By interchanging the roles of $A$ and $B$ we obtain also $\sup B+r\ge \sup A$ and $\inf B-r\le \inf A$ proving that $\delta\le r$.
\par Let us prove that $r\le \delta$. Suppose that both $A$ and $B$ are bounded, i.e., $A=[\inf A,\sup A]$ and
$B=[\inf B,\sup B]$. Let $\alpha\in[0,1]$, and let
$$
a_\alpha=\alpha\inf A+(1-\alpha)\sup A~~~\text{and}~~~
b_\alpha=\alpha\inf B+(1-\alpha)\sup B.
$$
Then $a_\alpha\in A$, $b_\alpha\in B$, and $|a_\alpha-b_\alpha|\le \delta$ proving that $\dist(a,B)\le r$ for every $a\in A$ and $\dist(b,A)\le r$ for every $b\in B$. Hence, $A+rI_0\supset B$ and $B+rI_0\supset A$, so that, thanks to \rf{HD-DF}, $r\le \delta$. In a similar way, we prove this inequality whenever one of the intervals is unbounded. We leave the details to the interested reader as an easy exercise.
\par (ii) By orthogonal projecting to the coordinate axes, we can reduce the problem to the one dimensional case. In this case, given bounded intervals $A,B\in\Ic(\R)$, we have
$$
\cent(A)=(\inf A+\sup A)/2~~~\text{and}~~~
\cent(B)=(\inf B+\sup B)/2.
$$
This and inequality \rf{AB-END} imply the required inequality
$$
|\cent(A)-\cent(B)|\le\dhf(A,B)
$$
proving the claim.\bx
\msk
\par Let $\MR$ be a finite pseudometric space and let $\Tc:\Mc\to \RCT$ be a set-valued mapping.
\par Given $\eta\ge 0$, let
\bel{TAU-1}
\Tc^{[1]}[x:\eta]=\bigcap_{z\in\Mc}\,
\left[\Tc(z)+\eta\,\rho(x,z)\,Q_0\right],~~~x\in\Mc,
\ee
be the $\eta$-metric refinement of $\Tc$. See Definition \reff{M-REF}.

\begin{proposition}\lbl{X2-C} 
\par (a) The set $\Tc^{[1]}[x:\eta]$ is not empty for every $x\in\Mc$ if and only if
\bel{TD-DR}
\Tc(x)\cap\{\Tc(y)+\eta\,\rho(x,y)Q_0\}\ne\emp~~~
\text{for all}~~~ x,y\in\Mc.
\ee
\par (b) If
\bel{TAU-NE}
\Tc^{[1]}[x:\eta]\ne\emp~~~\text{for every}~~~x\in\Mc,
\ee
then
\bel{TAM-1}
\dhf\left(\Tc^{[1]}[x:\eta],\Tc^{[1]}[y:\eta]\right)
\le \eta\,\rho(x,y)~~~\text{for all}~~~x,y\in\Mc.
\ee
\par Furthermore, if \rf{TAU-NE} holds and the set $\Tc^{[1]}[x:\eta]$ is bounded for every $x\in\Mc$, then the mapping
$$
\tau(x)=\cent\left(\Tc^{[1]}[x:\eta]\right),~~~ x\in\Mc,
$$
is a Lipschitz selection of $\Tc$ with $\|\tau\|_{\Lip(\Mc)}\le\eta$.
\end{proposition}
\par {\it Proof.} {\it (a)} Thanks to definition \rf{TAU-1}, for every $x,y\in\Mc$,
$$
\Tc^{[1]}[x:\eta]\subset\Tc(x)\cap
\left[\Tc(y)+\eta\,\rho(x,y)\,Q_0\right]
$$
so that if $\Tc^{[1]}[x:\eta]\ne\emp$ on $\Mc$ then condition \rf{TD-DR} holds.
\smsk
\par Now, let us assume that \rf{TD-DR} holds and prove that $\Tc^{[1]}[x:\eta]\ne\emp$ for every $x\in\Mc$.
\par First, given $x\in\Mc$, let us prove that
\bel{FFP}
\{\Tc(z)+\eta\,\rho(x,z)\,Q_0\}
\cap\{\Tc(z')+\eta\,\rho(x,z')\,Q_0\}
\ne\emp~~~~~\text{for every}~~~z,z'\in\Mc.
\ee
\par Thanks to \rf{TD-DR}, there exist points $a\in \Tc(z)$ and $a'\in \Tc(z')$ such that $\|a-a'\|\le \eta\,\rho(z,z')$. Therefore, thanks to the triangle inequality,
$$
\|a-a'\|\le \eta\,\rho(z,z')\le \eta\,\rho(z,x)+\eta\,\rho(x,z').
$$
This implies the existence of a point $w\in\RT$ such that
$$
\|a-w\|\le \eta\,\rho(z,x)~~~\text{and}~~~ \|a'-w\|\le \eta\,\rho(z',x).
$$
\par But $a\in \Tc(z)$ and $a'\in \Tc(z')$ so that $w$ belongs to the left hand side of \rf{FFP} proving this property.
\smsk
\par Property \rf{FFP} tells us that every two rectangles from the family $\{\Tc(z)+\eta\,\rho(x,z):z\in\Mc\}$ have a common point. Therefore, thanks to Lemma \reff{H-R}, the intersection of these rectangles, i.e., the set $\Tc^{[1]}[x:\eta]$ (see \rf{TAU-1}), is nonempty.
\msk
\par {\it (b)} We know that $\Tc^{[1]}[x:\eta]\ne\emp$, $x\in\Mc$, so that, thanks to Lemma \reff{H-RB} and definition \rf{TAU-1}, we have
\be
\Tc^{[1]}[x:\eta]+\eta\,\rho(x,y)Q_0&=&
\left\{\bigcap_{z\in \Mc}
\left[\Tc(z)+\eta\,\rho(x,z)Q_0\right]\right\}
+\eta\,\rho(x,y)Q_0
\nn\\
&=&\bigcap_{z\in \Mc}
\left[\Tc(z)+(\eta\,\rho(x,z)+\eta\,\rho(x,y))\,Q_0\right].
\nn
\ee
From this and the triangle inequality, we have
$$
\Tc^{[1]}[x:\eta]+\eta\,\rho(x,y)\,Q_0\supset
\bigcap_{z\in \Mc}\,
\left[\Tc(z)+\eta\,\rho(y,z)\,Q_0\right]=\Tc^{[1]}[y:\eta].
$$
\par By interchanging the roles of $x$ and $y$ we also obtain
$$
\Tc^{[1]}[y:\eta]+\eta\,\rho(x,y)\,Q_0\supset \Tc^{[1]}[x:\eta].
$$
These two inclusions imply the required inequality \rf{TAM-1}.
\smsk
\par Finally, inequalities \rf{CAB-H} and \rf{TAM-1} imply the required inequality $\|\tau\|_{\Lip(\Mc)}\le\eta$ completing the proof of the proposition.\bx
\msk
\par Let us represent the set-valued mapping
$\Tc:\Mc\to\RCT$ in the form
\bel{PR-F}
\Tc(x)=I_1(x)\times I_2(x)~~~
\text{where}~~~
I_i(x)=\Prj_i[\Tc(x)],~~~x\in\Mc,~i=1,2.
\ee
\par Let $s_1(x)$ and $s_2(x)$ be the left and right ends of $I_1(x)$, see \rf{LR1-END}, i.e.,
\bel{E-S12}
s_1(x)=\inf\{s\in\R:(s,0)\in I_1(x)\}~~~\text{and}~~~
s_2(x)=\sup\{s\in\R:(s,0)\in I_1(x)\},
\ee
and let $t_1(x)$ and $t_2(x)$ be the left and right ends of $I_2(x)$, see \rf{LR2-END}, i.e.,
\bel{T12-F}
t_1(x)=\inf\{t\in\R:(0,t)\in I_2(x)\}~~~\text{and}~~~
t_2(x)=\sup\{t\in\R:(0,t)\in I_2(x)\}.
\ee
\par Given $\eta\ge 0$, let us introduce similar objects for the set-valued mapping $\Tc^{[1]}[\cdot:\eta]:\Mc\to\RCT$, see \rf{TAU-1}, provided condition \rf{TAU-NE} holds. We represent this mapping in the form
\bel{PRR-1}
\Tc^{[1]}[x:\eta]=\tI_1(x)\times \tI_2(x)~~~\text{where}~~~
\tI_i(x)=\Prj_i[\Tc^{[1]}[x:\eta]],~~~x\in\Mc,~i=1,2.
\ee
\par Let $\s_1(x)$ and $\s_2(x)$ be the left and right ends of $\tI_1(x)$, see \rf{LR1-END}, and  let $\st_1(x)$ and $\st_2(x)$ be the left and right ends of $\tI_2(x)$, see \rf{LR2-END}.
\begin{lemma}\lbl{RF-ENDS} Suppose that
$\Tc^{[1]}[x:\eta]\ne\emp$ for every $x\in\Mc$. Then for every $x\in\Mc$, the following equalities hold:
\bel{S12-L}
\s_1(x)=\max\{s_1(y)-\eta\,\rho(x,y):y\in\Mc\},
~~~
\s_2(x)=\min\{s_2(y)+\eta\,\rho(x,y):y\in\Mc\},
\ee
and
\bel{T12-L}
\st_1(x)=\max\{t_1(y)-\eta\,\rho(x,y):y\in\Mc\},
~~~
\st_2(x)=\min\{t_2(y)+\eta\,\rho(x,y):y\in\Mc\}.
\ee
\end{lemma}
\par {\it Proof.} Let us prove the first equality in \rf{S12-L}. Thanks to \rf{TAU-1}, \rf{PR-F} and \rf{PRR-1},
\be
\tI_1(x)&=&
\Prj_1[\Tc^{[1]}[x:\eta]]=
\Prj_1\left[\bigcap_{y\in \Mc}\,
\left[\Tc(y)+\eta\,\rho(x,y)\,Q_0\right]
\right]=\bigcap_{y\in \Mc}\,
\Prj_1\left[\Tc(y)+\eta\,\rho(x,y)\,Q_0\right]
\nn\\
&=&
\bigcap_{y\in \Mc}\,
\left\{\Prj_1\left[\Tc(y)\right]+
\eta\,\rho(x,y)\,\Prj_1[Q_0]\right\}
=\bigcap_{y\in \Mc}\,
\left\{I_1(x)+\eta\,\rho(x,y)\,J_1\right\}
\nn
\ee
where $J_1=\Prj_1[Q_0]=[-e_1,e_1]$ and $e_1=(1,0)$.
\par Thanks to \rf{E-S12}, for every $x,y\in\Mc$, the set
$$
\Ic(x,y)=I_1(x)+\eta\,\rho(x,y)\,J_1
$$
is a line segment in $Ox_1$ with
$$
\text{the left end}~~~L(\Ic(x,y))=s_1(x)-\eta\,\rho(x,y)
~~~
\text{and the right end} ~~~R(\Ic(x,y))=s_2(x)+\eta\,\rho(x,y).
$$
See \rf{LR1-END}. Therefore,
$$
\s_1(x)=L\left(\tI_1(x)\right)=
L\left(\bigcap_{y\in \Mc}\,\Ic(x,y)\right)=
\max_{y\in \Mc} L\left(\Ic(x,y)\right)=
\max_{y\in \Mc}\{s_1(x)-\eta\,\rho(x,y)\}
$$
proving the first equality in \rf{S12-L}.
\par In a similar way, we prove the second equality in \rf{S12-L} and equalities of \rf{T12-L} completing the proof of the lemma.\bx
\msk
\par  Let us formulate one more useful statement that we will need in Section 6.3.
\begin{claim}\lbl{CL-M} Given $x,y\in\Mc$, we have
\bel{TD-XY}
\Tc(x)\cap\{\Tc(y)+\eta\,\rho(x,y)Q_0\}\ne\emp
\ee
if and only if
\bel{CL-REN}
\max\{s_1(x)-s_2(y),s_1(y)-s_2(x),
t_1(x)-t_2(y),t_1(y)-t_2(x)\}\,\le\,\eta\,\rho(x,y).
\ee
\end{claim}
\par {\it Proof.} We recall that, given $u\in\Mc$, we  represent the set-valued mapping $\Tc:\Mc\to\RCT$ in the form
$$
\Tc(x)=I_1(x)\times I_2(x)~~~
\text{where}~~~
I_i(x)=\Prj_i[\Tc(x)],~~~x\in\Mc,~i=1,2.
$$
See \rf{PR-F}.
\par We also recall that $s_1(x)$ and $s_2(x)$ are the left and right ends of $I_1(x)$, see \rf{E-S12}, and $t_1(x)$ and $t_2(x)$ are the left and right ends of $I_2(x)$, see \rf{T12-F}.
\smsk
\par Let
$$
W(y)=\Tc(y)+\eta\,\rho(x,y)\,Q_0
$$
and let
$$
\hI_1(y)=\Prj_1[W(y)]~~~\text{and}~~~
\hI_2(y)=\Prj_2[W(y)].
$$
Thus $W(y)=\hI_1(y)\times \hI_2(y)$. Furthermore,
$$
\hI_i(y)=\Prj_i[W(y)]=
\Prj_i[\Tc(y)+\eta\,\rho(x,y)\,Q_0]=
\Prj_i[\Tc(y)]+
\eta\,\rho(x,y)\,\Prj_i[Q_0].
$$
Therefore,
$$
\hI_i(y)=I_i(y)+\eta\,\rho(x,y)\,J_i
$$
where $J_i=\Prj_i[Q_0]=[-e_i,e_i]$, $i=1,2$. (Recall that $e_1=(1,0)$ and $e_2=(0,1)$.)
\par Because $s_1(x)$ and $s_2(x)$ are the left and right ends of $I_1(x)$, and $t_1(x)$ and $t_2(x)$ are the left and right ends of $I_2(x)$, the set $\hI_1(y)$ is a line segment on $Ox_1$ with the left end
$$
L(\hI_1(y))=s_1(x)-\eta\,\rho(x,y)~~~\text{and right end}~~~R(\hI_1(y))=s_2(x)+\eta\,\rho(x,y).
$$
In turn, $\hI_2(y)$ is a line segment on $Ox_2$ with the left end
$$
L(\hI_2(y))=t_1(x)-\eta\,\rho(x,y)~~~
\text{and right end}~~~
R(\hI_t(y))=t_2(x)+\eta\,\rho(x,y).
$$
\par Part {\it (ii)} of Claim \reff{CLH-R} tells us that condition \rf{TD-XY} holds (i.e., $\Tc(x)\cap W(y)\ne\emp$) if and only if
$$
I_1(x)\cap\hI_1(y)\ne\emp~~~\text{and}~~~ I_2(x)\cap\hI_2(y)\ne\emp.
$$
\par Clearly, two line segments on the axis $Ox_1$ have a common point if and only if the left end of the first of them is located to the left of the right end of the second, and vice versa, the left end of the second is located to the left of the right end of the first.
Therefore, $I_1(x)\cap\hI_1(y)\ne\emp$ if and only if
$$
s_2(y)+\eta\,\rho(x,y)\ge s_1(x)~~~
\text{and}~~~
s_2(x)\ge s_1(y)-\eta\,\rho(x,y).
$$
which is equivalent to the inequality
$$
\max\{s_1(x)-s_2(y),s_1(y)-s_2(x)\}\,\le
\,\eta\,\rho(x,y).
$$
\par In a similar way, we prove that $I_2(x)\cap\hI_2(y)\ne\emp$ if and only if
$$
\max\{t_1(x)-t_2(y),t_1(y)-t_2(x)\}\,
\le\,\eta\,\rho(x,y).
$$
\par This proves that property \rf{TD-XY} holds if and only if inequality \rf{CL-REN} is satisfied, and the proof of the claim is complete.\bx

\msk
\par Let us give several explicit formulae for Lipschitz selections of set-valued mappings in the one dimensional case. Let $G:\Mc\to \Ic(\R)$ be a set-valued mapping. We set
\bel{AB-11}
a_G(x)=\inf G(x)~~~~\text{and}~~~~b_G(x)=\sup G(x).
\ee
Thus, $a_G$ and $b_G$ are two functions on $\Mc$ such that
$$
a_G:\Mc\to\R\cup\{-\infty\},~~~
b_G:\Mc\to\R\cup\{+\infty\}~~~~\text{and}~~~~a_G(x)\le b_G(x)~~~\text{for all}~~~x\in\Mc.
$$
\par Clearly,
\bel{DR-XY}
\dist(G(x),G(y))=\max\{[a_G(x)-b_G(y)]_+,[a_G(y)-b_G(x)]_+\}.
\ee
(See our convention \rf{INF-S} for the case of $a_G(x)=-\infty$, $b_G(x)=+\infty$.)
\par Given $\eta\ge 0$, we introduce the following functions on $\Mc$:
\bel{FP-D1}
a^{[1]}_G[x:\eta]=\sup_{y\in\Mc}
\,\left\{a_G(y)-\eta\,\rho(x,y)\right\},~~~~~
b^{[1]}_G[x:\eta]=
\inf_{y\in\Mc}\,\left\{b_G(y)+\eta\,\rho(x,y)\right\}
\ee
and
\bel{FS-D1}
c_G[x:\eta]=\left(a^{[1]}_G[x:\eta]
+b^{[1]}_G[x:\eta]\right)/2.
\ee
\par Let $G^{[1]}[\cdot:\eta]$ be the $\eta$-metric refinement of $G$, i.e., a set-valued mapping on $\Mc$ defined by
\bel{F1-PR1}
G^{[1]}[x:\eta]=
\bigcap_{z\in\Mc}\,\left[G(z)+
\eta\,\rho(x,z)\,\BXR\right],
~~~x\in\Mc.
\ee
\par Comparing this definition with definitions \rf{FP-D1} and \rf{FS-D1}, we conclude that for every $x\in\Mc$,
\bel{F1-ENDS}
a^{[1]}_G[x:\eta]=\inf G^{[1]}[x:\eta],~~~~~
b^{[1]}_G[x:\eta]=\sup G^{[1]}[x:\eta],
\ee
provided $G^{[1]}[x:\eta]\ne\emp$ and
\bel{F1-CNT}
c_G[x:\eta]=\cent\left(G^{[1]}[x:\eta]\right)
\ee
provided $G^{[1]}[x:\eta]$ is nonempty and {\it bounded}.
\smsk
\begin{remark}\lbl{RM-ABF} {\em We note that the function
$g^+=b^{[1]}_G[\cdot:\eta]$ maps $\Mc$ into $\R$
if and only if $g^+\nequiv+\infty$, i.e., $g^+(x^+)<\infty$ for some $x^+\in\Mc$. Analogously, the mapping
$g^-=a^{[1]}_G[\cdot:\eta]$ maps $\Mc$ into $\R$
if and only if $g^-\nequiv-\infty$. Finally, the mapping
$c^{[1]}_G[\cdot:\eta]=(g^++g^-)/2$, see \rf{FS-D1} and \rf{F1-CNT}, is well defined if and only if
both $g^+\nequiv+\infty$ and $g^-\nequiv-\infty$. \rbx}
\end{remark}

\begin{proposition} Let $\Mc$ be finite and let
$G:\Mc\to \Ic(\R)$ be a set-valued mapping.
\par Given $x\in\Mc$ and $\eta\ge 0$, the set $G^{[1]}[x:\eta]$ is nonempty if and only if
\bel{AB-C}
a^{[1]}_G[x:\eta]\le b^{[1]}_G[x:\eta].
\ee
See \rf{FP-D1}. Furthermore, if \rf{AB-C} holds, then
\bel{E-F1}
a^{[1]}_G[x:\eta]=\inf G^{[1]}[x:\eta]~~~\text{and}~~~
b^{[1]}_G[x:\eta]=\sup G^{[1]}[x:\eta].
\ee
\end{proposition}
\par {\it Proof.} Suppose that $G^{[1]}[x:\eta]\ne\emp$. Then for every $y_1,y_2\in\Mc$, we have
$$
(G(y_1)+\eta\,\rho(x,y_1)Q_0)\cap
(G(y_2)+\eta\,\rho(x,y_2)Q_0)\ne\emp,
$$
so that
$$
\inf(G(y_1)+\eta\,\rho(x,y_1)Q_0)\le
\sup(G(y_2)+\eta\,\rho(x,y_2)Q_0).
$$
\par Therefore, thanks to \rf{AB-11},
\bel{MA-B}
a_G(y_1)-\eta\,\rho(x,y_1)\le b_G(y_1)+\eta\,\rho(x,y_2)
~~~\text{for all}~~~y_1,y_2\in\Mc,
\ee
so that, thanks to \rf{FP-D1}, inequality \rf{AB-C} holds.
\smsk
\par Let us now assume that \rf{AB-C} holds. Then, for every $y_1,y_2\in\Mc$ inequality \rf{MA-B} holds as well so
that
$$
a_G(y_2)-b_G(y_1)\le\eta\,\rho(x,y_1)+\eta\,\rho(x,y_2).
$$
By interchanging the roles of $y_1$ and $y_2$, we obtain also
$$
a_G(y_1)-b_G(y_2)\le\eta\,\rho(x,y_1)+\eta\,\rho(x,y_2).
$$
\par From these inequalities and formula \rf{DR-XY}, we have
$$
\dist(G(y_1),G(y_2))\le \eta\,\rho(x,y_1)+\eta\,\rho(x,y_2).
$$
\par But because $G(y_1)$ and $G(y_2)$ are two closed intervals in $\R$, there exist points $u_1\in G(y_1)$ and $u_2\in G(y_2)$ such that
$$
\|u_1-u_2\|\le \eta\,\rho(x,y_1)+\eta\,\rho(x,y_2).
$$
Therefore, there exists a point $v\in\R$ such that
$$
\|u_1-v\|\le \eta\,\rho(x,y_1)~~~\text{and}~~~
\|u_2-v\|\le \eta\,\rho(x,y_2).
$$
Clearly,
$$
v\in(G(y_1)+\eta\,\rho(x,y_1)Q_0)\cap
(G(y_2)+\eta\,\rho(x,y_2)Q_0).
$$
\par This proves that any two closed intervals from the family
$$
\Gc=\{G(y)+\eta\,\rho(x,y)Q_0:y\in\Mc\}
$$
have a common point. Therefore, thanks to Helly's theorem in $\R$, there exists a point $w\in\R$ common to all of the intervals from $\Gc$.
\par Thus, $w$ belongs to the intersection of all these intervals, and therefore, thanks to definition \rf{F1-PR1}, $w\in G^{[1]}[x:\eta]$ proving that $G^{[1]}[x:\eta]\ne\emp$.
\par Finally, it remains to note that property \rf{E-F1} coincides with \rf{F1-ENDS}, and the proof of the proposition is complete.\bx
\begin{proposition}\lbl{FP-R1} (The Finiteness Principle for Lipschitz selections in $\R$.) Let $\Mc$ be finite and let $G:\Mc\to \Ic(\R)$ be a set-valued mapping. Let $\eta\ge 0$.
\par Suppose that for every $x,y\in\Mc$ the restriction $G|_{\{x,y\}}$ of $G$ to $\{x,y\}$ has a Lipschitz selection $g_{\{x,y\}}$ with $\|g_{\{x,y\}}\|_{\Lip(\{x,y\},\R)}\le \eta$. Then $G$ has a Lipschitz selection $g:\Mc\to\R$ with Lipschitz seminorm $\|g\|_{\Lip(\Mc,\R)}\le \eta$.
\par Furthermore, one can set $g=c^{[1]}_G[\cdot:\eta]$ provided there exist $x^+,x^-\in\Mc$ such that $\inf G(x^-)>-\infty$ and $\sup G(x^+)<\infty$. Also, one can set $g=b^{[1]}_G[\cdot:\eta]$ if $G(x^+)$ is bounded from above for some $x^+\in\Mc$, or $g=a^{[1]}_G[\cdot:\eta]$ if $G(x^-)$ is bounded from below for some $x^-\in\Mc$.
\end{proposition}
\par {\it Proof.} The hypothesis of the proposition tells us that for every $x,y\in\Mc$ there exist points $a\in G(x)$ and $b\in G(x_2)$ such that $\|a-b\|\le \eta\,\rho(x,y)$. Hence,
$$
G(x)\cap\{G(y)+\eta\,\rho(x,y)I_0\}\ne\emp~~~\text{for all}~~~ x,y\in\Mc.
$$
Here $I_0=[-1,1]$. Proposition \reff{X2-C} tells us that in this case $G^{[1]}[x:\eta]\ne\emp$ for every $x\in\Mc$, and
$$
\dhf\left(G^{[1]}[x:\eta],G^{[1]}[y:\eta]\right)
\le \eta\,\rho(x,y)~~~\text{for all}~~~x,y\in\Mc.
$$
From this, part (i) of Claim \reff{TWO}, and definitions \rf{F1-ENDS} it follows that the inequality
\bel{AB-HDG}
\max\left\{\,|a^{[1]}_G[x:\eta]-a^{[1]}_G[y:\eta]|,\,
|b^{[1]}_G[x:\eta]-b^{[1]}_G[y:\eta]|\,\right\}
\le \eta\,\rho(x,y)
\ee
holds for all $x,y\in\Mc$.
\par Clearly, if $G\equiv \R$ then the constant mapping $g\equiv \{0\}$ on $\Mc$ is a Lipschitz selection of $G$ (with $\|g\|_{\Lip(\Mc,\R)}=0$). Otherwise, either $g^+=b^{[1]}_G[\cdot:\eta]\nequiv+\infty$ or $g^-=a^{[1]}_G[\cdot:\eta]\nequiv-\infty$. Therefore, thanks to Remark \reff{RM-ABF}, either $g^+:\Mc\to\R$
or $g^-:\Mc\to\R$.
\par Let us note that if at least one of the sets $G(x)$ is bounded, then the set $G^{[1]}[x:\eta]$ is bounded as well. In this case, the points
$$
a^{[1]}_G[x:\eta],b^{[1]}_G[x:\eta]\in G^{[1]}[x:\eta]\subset G(x).
$$
See definition \rf{F1-ENDS}. Therefore, in this case we can set either $g=g^+$ or $g=g^-$. Then, thanks to \rf{AB-HDG}, in both cases the mapping $g:\Mc\to\R$ will be a Lipschitz selection of $G$ with $\|g\|_{\Lip(\Mc,\R)}\le\eta$. Also from this it follows that the mapping $g=c_G[\cdot:\eta]=(g^+ + g^-)/2$, see \rf{F1-CNT}, has the same properties.
\par Note that, if   $g^+=b^{[1]}_G[\cdot:\eta]\nequiv+\infty$ but  $g^-=a^{[1]}_G[\cdot:\eta]\equiv-\infty$, then each interval $G^{[1]}[x:\eta]$ is {\it unbounded from below}. Because $\Mc$ is finite, all these intervals have a common point, say $A$. Then the constant mapping $g\equiv \{A\}$ is a Lipschitz selection of $G$ (with $\|g\|_{\Lip(\Mc,\R)}=0$). Analogously, if $g^-\nequiv-\infty$ but $g^+\equiv+\infty$, there is a constant mapping which provides a Lipschitz selection of $G$.
\par Let us suppose that both $g^+=b^{[1]}_G[\cdot:\eta]\nequiv+\infty$ and $g^-=a^{[1]}_G[\cdot:\eta]\nequiv-\infty$. In this case, thanks to Remark \reff{RM-ABF}, the mapping
$c^{[1]}_G[\cdot:\eta]=(g^++g^-)/2$, see \rf{FS-D1} and \rf{F1-CNT}, is well defined, i.e., each set $G^{[1]}[x:\eta]$ is nonempty and bounded. Clearly,
$$
g(x)=c^{[1]}_G[x:\eta]\in G^{[1]}[x:\eta]\subset G(x)
~~~~\text{for every}~~~~x\in\Mc,
$$
proving that $g$ is a selection of $G$. Thanks to \rf{AB-HDG}, its Lipschitz seminorm  $\|g\|_{\Lip(\Mc,\R)}\le\eta$, and the proof of the proposition is complete.\bx

\smsk
\par Proposition \reff{FP-R1} implies the following Finiteness Principle for rectangles in $\RT$.
\begin{proposition}\lbl{FP-RC} Let $\eta\ge 0$. Let $\MR$ be a finite pseudometric space, and let $\Tc:\Mc\to \RCT$ be a set-valued mapping. Suppose that for every $x,y\in\Mc$ the restriction $\Tc|_{\{x,y\}}$ of $\Tc$ to $\{x,y\}$ has a Lipschitz selection $\tau_{\{x,y\}}$ with $\|\tau_{\{x,y\}}\|_{\Lip(\{x,y\})}\le \eta$.
\par Then $\Tc$ has a Lipschitz selection $\tau:\Mc\to\RT$ with Lipschitz seminorm $\|\tau\|_{\Lip(\Mc)}\le \eta$.
\end{proposition}
\par {\it Proof.} By orthogonal projecting onto the coordinate axes, we reduce the problem to the Finiteness Principle for Lipschitz selections in $\R$ proven in Proposition \reff{FP-R1}.\bx
\smsk

\SECT{3. The key theorem: Lipschitz selections and rectangular hulls.}{3}

\addtocontents{toc}{3. The key theorem: Lipschitz selections and rectangular hulls.\hfill \thepage\par\VST}


\indent
\par We proceed to the proof of Theorem \reff{W-CR}. This proof is based on an essential refinement of the approach developed in the author's paper \cite[Section 3.1]{S-2002}.
\par Let $\MS=(\Mc,\rho)$ be a finite pseudometric space, and let $F:\Mc\to \CNV$ be a set-valued mapping. We recall that, given $\lambda\ge 0$ and $x,x',x''\in\Mc$, by   $\Wc_F[x,x',x'':\lambda]$ we denote a (possibly empty) subset of $\RT$ defined by
\bel{WC-DF}
\Wc_F[x,x',x'':\lambda]=
\HR[\{F(x')+\lambda\,\rho(x',x)\,Q_0\}
\cap \{F(x'')+\lambda\,\rho(x'',x)\,Q_0\}].
\ee
\par (Recall also that $\HR[\cdot]$ denotes the rectangular hull of a set. See \rf{HRS}.)

\msk

\par We begin the proof of the key Theorem \reff{W-CR} in this section and conclude it in the next section.
\msk
\par \underline{\sc Proof of Theorem \reff{W-CR}.} Suppose that for every $x,x',x'',y,y',y''\in\Mc$ condition \rf{WNEW} holds.
\msk

\par Recall \rf{WNEW}: it states that, given non-negative constants $\tlm$ and $\lambda$, the following condition
$$
\Wc_F[x,x',x'':\tlm]\cap \{\Wc_F[y,y',y'':\tlm]+\lambda\,\rho(x,y)\,Q_0\}\ne\emp
$$
holds for every $x,x',x'',y,y',y''\in\Mc$. In turn, this condition is equivalent the following one: for every $x,x',x'',y,y',y''\in\Mc$ there exist points $u,v\in\RT$ such that
$$
u\in\HR[\{F(x')+\lambda\,\rho(x',x)\,Q_0\}
\cap \{F(x'')+\lambda\,\rho(x'',x)\,Q_0\}],
$$
\smsk
$$
v\in\HR[\{F(y')+\lambda\,\rho(y',y)\,Q_0\}
\cap \{F(y'')+\lambda\,\rho(y'',y)\,Q_0\}],
$$
and
$$
\|u-v\|\le\lambda\,\rho(x,y)
$$
See Fig. 10.
\msk

\begin{figure}[H]
\includegraphics[scale=1.05]{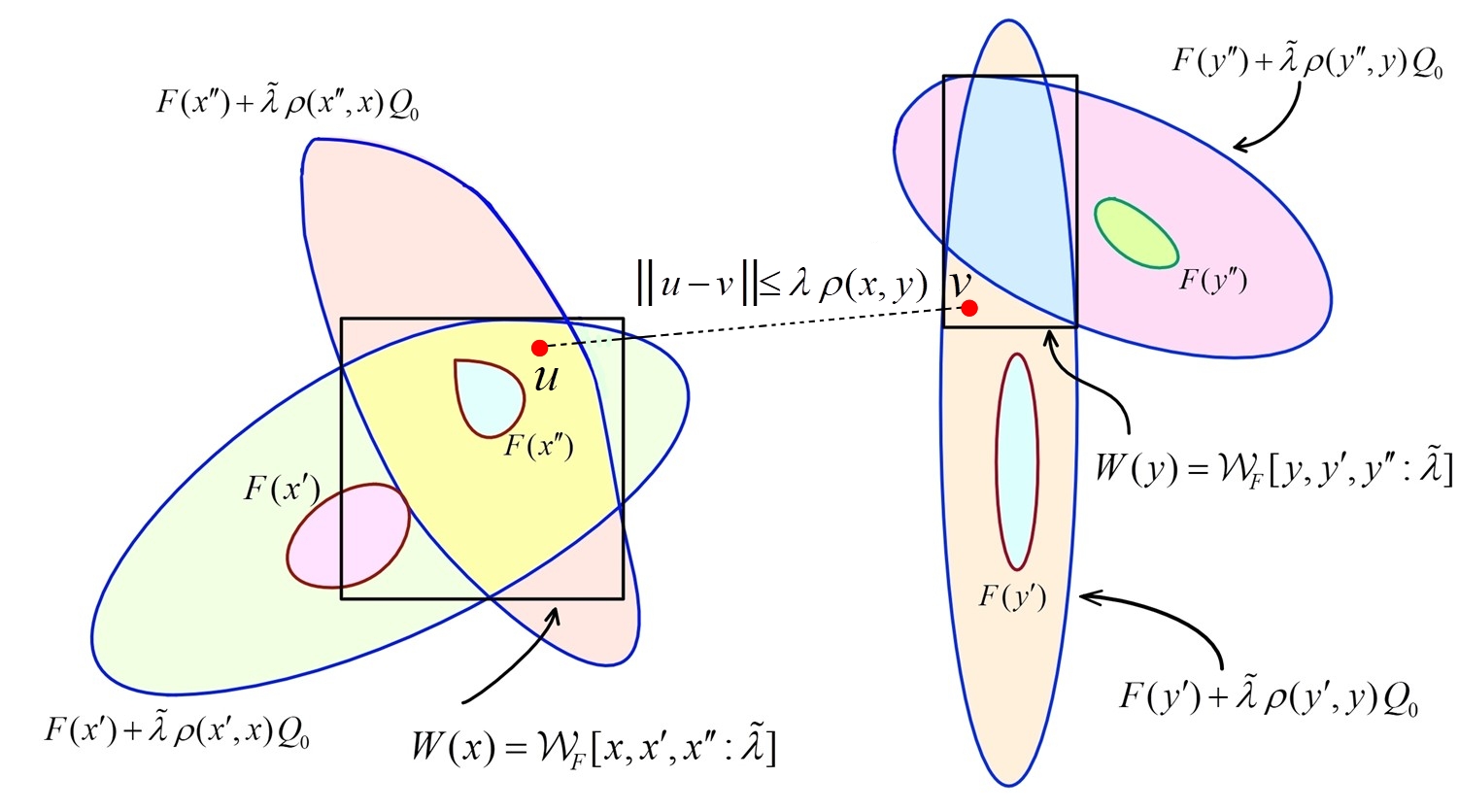}
\caption{The hypothesis of Theorem 1.9.}
\end{figure}

Let us construct a Lipschitz selection $f:\Mc\to\RT$ of $F$ with Lipschitz seminorm $\|f\|_{\Lip(\Mc)}\le 2\lambda+\tlm$. We will do this in three steps.
\bsk
\par \underline{\sc The First Step.} We introduce the $\tlm$-metric refinement of $F$, see \rf{F-1}, i.e., a set-valued mapping on $\Mc$ defined by the formula
\bel{WF1-3L}
F^{[1]}[x:\tlm]=
\bigcap_{y\in \Mc}\,\left[F(y)+\tlm\,\rho(x,y)\,Q_0\right],
~~~~x\in\Mc.
\ee
\begin{lemma}\lbl{WF1-NE} For each $x\in\Mc$, the set
$F^{[1]}[x:\tlm]$ is a nonempty closed convex subset of $\RT$. Moreover, for every $x\in\Mc$ the following representation holds:
\bel{WH-F1}
\HR[F^{[1]}[x:\tlm]]=\cap\{\Wc_F[x,y,y':\tlm]: y,y'\in\Mc\}.
\ee
\end{lemma}
\par {\it Proof.} Let us prove that
\bel{F1-IN}
F^{[1]}[x:\tlm]\ne\emp~~~~\text{for every}~~~~x\in\Mc.
\ee
\par Given $x\in\Mc$, we set
$$
\Cf_x=\{F(y)+\tlm\,\rho(x,y)\,Q_0:y\in\Mc\}.
$$
Then $F^{[1]}[x:\tlm]=\cap\{C:C\in\Cf_x\}$. See \rf{WF1-3L}.
\smsk
\par Let us prove that for every $y_1,y_1',y_2,y_2'\in\Mc$ the sets
\bel{C-PT}
C_i=F(y_i)+\tlm\,\rho(x,y_i)Q_0~~~~\text{and}~~~~ C'_i=F(y'_i)+\tlm\,\rho(x,y'_i)Q_0,~~ i=1,2,
\ee
satisfy property \rf{P1-C}.
\par First, let us note that, thanks to \rf{WNEW}, the set
$\Wc_F[x,y,z:\tlm]\ne\emp$ for all $x,y,z\in\Mc$.
\par In particular, from this and definition \rf{WC-DF}, it follows that
$$
\{F(y)+\tlm\,\rho(x,y)\,Q_0\}\cap
\{F(z)+\tlm\,\rho(x,z)\,Q_0\}\ne\emp~~~~
\text{for every}~~~~y,z\in\Mc
$$
proving that {\it any two elements of $\Cf_x$ have a common point}.
\par Property \rf{WNEW} tells us that
$$
\Wc_F[x,y_1,y'_1:\tlm]\cap\Wc_F[x,y_2,y'_2:\tlm]
\ne\emp~~~~\text{for every}~~~~y_1,y'_1,y_2,y'_2\in\Mc.
$$
Thanks to \rf{WC-DF} and \rf{C-PT},
$$
\Wc_F[x,y_1,y'_1:\tlm]=\Hc[C_1\cap C_1']~~~\text{and}~~~ \Wc_F[x,y_2,y'_2:\tlm]=\Hc[C_2\cap C_2']
$$
proving that
$$
\Hc[C_1\cap C_1']\,\cbg \Hc[C_2\cap C_2'] \ne\emp~~~~~\text{for every}~~~~~ C_1,C_1',C_2,C_2'\in\Cf_x.
$$
Hence,
$$
\Prj_1[\Hc[C_1\cap C_1']]\,\cbg \Prj_1[\Hc[C_2\cap C_2']] \ne\emp
$$
so that, thanks to \rf{HS-U},
$$
\Prj_1[C_1\cap C_1']\,\cbg \Prj_1[C_2\cap C_2'] \ne\emp
$$
proving \rf{P1-C}.
\smsk
\par Because the set $\Mc$ is {\it finite}, the family $\Cf_x$ is finite as well. Therefore, thanks to Proposition \reff{INT-RE}, $\cap\{C:C\in\Cf_x\}\ne\emp$ proving \rf{F1-IN}.
\smsk
\par Finally, \rf{PH-C} and \rf{WC-DF} imply formula \rf{WH-F1} completing the proof of the lemma.\bx
\msk
\par Following {\bf STEP 2} of the Projection Algorithm, given $x\in\Mc$, we let $\Tc_{F,\tlm}(x)$ denote the rectangular hull of the set $F^{[1]}[x:\tlm]$. Cf. \rf{TC-DF}. Thus, $\Tc_{F,\tlm}$ is a set-valued mapping from $\Mc$ into $\RCT$ defined by
\bel{WTH-F1}
\Tc_{F,\tlm}(x)=\HR[F^{[1]}[x:\tlm]]=
\HR\left[\bigcap_{y\in \Mc}
\left\{F(y)+\tlm\,\rho(x,y)\,Q_0\right\}\right], ~~~~~x\in\Mc.
\ee
Let us note that formula \rf{WH-F1} provides the following representation of the mapping $\Tc_{F,\tlm}$:
\bel{WR-TC}
\Tc_{F,\tlm}(x)=\bigcap_{x',x''\in\Mc}\Wc_F[x,x',x'':\tlm],
~~~~x\in\Mc.
\ee

\par Let us note that, thanks to Lemma \reff{WF1-NE}, $F^{[1]}[x:\tlm]\ne\emp$ for every $x\in\Mc$. Therefore, thanks to \rf{WTH-F1},
\bel{TF-NE}
\Tc_{F,\tlm}(x)\ne\emp~~~\text{for all}~~~x\in\Mc.
\ee
\smsk

\par \underline{\sc The Second Step.} At this step we prove the existence of a Lipschitz selection of the set-valued mapping $\Tc_{F,\tlm}$.
\par First, following {\bf STEP 3} of the Projection Algorithm, we let $\Tc_{F,\tlm}^{[1]}[\cdot:\lambda]$ denote the $\lambda$-metric refinement of the mapping $\Tc_{F,\tlm}$, see \rf{RT-1}. Thus,
$$
\Tc_{F,\tlm}^{[1]}[x:\lambda]=
\bigcap_{z\in\Mc}\,\left[\Tc_{F,\tlm}(z)+
\lambda\,\rho(x,z)\,Q_0\right],~~~x\in\Mc.
$$

\begin{proposition}\lbl{WLS-T} (i) The set-valued mapping $\Tc_{F,\tlm}:\Mc\to\RCT$ has a Lipschitz selection $g:\Mc\to\RT$ with Lipschitz seminorm $\|g\|_{\Lip(\Mc)}\le \lambda$;
\smsk
\par (ii) For every $x\in\Mc$, the following property
\bel{TC-NEMP}
\Tc_{F,\tlm}^{[1]}[x:\lambda]\ne\emp
\ee
holds. Furthermore, for every $x,y\in\Mc$, we have
\bel{DRT-1}
\dhf\left(\Tc_{F,\tlm}^{[1]}[x:\lambda],
\Tc_{F,\tlm}^{[1]}[y:\lambda]\right)\le \lambda\,\rho(x,y).
\ee
(Recall that $\dhf$ denotes the Hausdorff distance between sets.)
\smsk
\par (iii) If each rectangle $\Tc_{F,\tlm}^{[1]}[x:\lambda]$, $x\in\Mc$, is bounded, then the mapping
$$
g_F(x)=\cent\left(\Tc_{F,\tlm}^{[1]}[x:\lambda]\right),
~~~~x\in\Mc,
$$
is a Lipschitz selection of~ $\Tc_{F,\tlm}$ with  $\|g_F\|_{\Lip(\Mc)}\le\lambda$.
\end{proposition}
\par {\it Proof.} {\it (i)} Proposition \reff{FP-RC} tells us that the required Lipschitz selection $g$ exists provided for every $x,y\in\Mc$ the restriction $\Tc_{F,\tlm}|_{\{x,y\}}$ of $\Tc_{F,\tlm}$ to $\{x,y\}$ has a Lipschitz selection $g_{\{x,y\}}$ with Lipschitz  seminorm $\|g_{\{x,y\}}\|_{\Lip(\{x,y\})}\le \lambda$. Clearly, this requirement is equivalent to the following property:
\bel{WDTC}
\Tc_{F,\tlm}(x)\cap\{\Tc_{F,\tlm}(y)
+\lambda\,\rho(x,y)Q_0\}\ne\emp
~~~~\text{for every}~~~~x,y\in\Mc.
\ee
\par Let us prove that this property holds. Let $x,y\in\Mc$ and let
$$
\Tf_x=\{\Wc_F[x,x',x'':\tlm]: x',x''\in\Mc\} ~~~\text{and}~~~
\Tf_y=\{\Wc_F[y,y',y'':\tlm]: y',y''\in\Mc\}.
$$
Thanks to \rf{WR-TC},
\bel{TF-XY}
\Tc_{F,\tlm}(x)=\cap\{W:W\in\Tf_x\}~~~~\text{and}~~~~
\Tc_{F,\tlm}(y)=\cap\{W:W\in\Tf_y\}.
\ee
\par Since $\Mc$ is finite, $\Tf_x$ and $\Tf_y$ are {\it finite} families of rectangles. Thanks to \rf{TF-NE}, the set $\Tc_{F,\tlm}(z)\ne\emp$ for every $z\in\Mc$, so that each family has a nonempty intersection.
\par Let $r=\lambda\,\rho(x,y)$. Then, thanks to \rf{TF-XY} and Lemma \reff{H-RB},
$$
\Tc_{F,\tlm}(y)+\lambda\,\rho(x,y)Q_0=\cap\{W:W\in\Tf_y\}+rQ_0=
\cap\{W+rQ_0:W\in\Tf_y\}
$$
so that
\bel{WQ-NM}
\Tc_{F,\tlm}(x)\cap\{\Tc_{F,\tlm}(y)+\lambda\,\rho(x,y)Q_0\}=
[\cap\{W:W\in\Tf_x\}]\cap [\cap\{W+rQ_0:W\in\Tf_y\}].
\ee
\par Let $\Tfw=\Tf_x\cup \Tf^+_y$ where $\Tf^+_y=\{W+rQ_0:W\in\Tf_y\}$. Thanks to \rf{WQ-NM}, property \rf{WDTC} holds provided the family of rectangles
$\Tfw$ has a common point. Because $\Tf_x$ and $\Tf_y$ are {\it finite} families, the family $\Tfw$ is finite as well.
Therefore, thanks to Helly's intersection theorem for rectangles, see Lemma \reff{H-R}, there exists a point common to all of the family $\Tfw$ provided $W'\cap W''\ne\emp$ for every $W',W''\in\Tfw$.
\smsk
\par Clearly, $W'\cap W''\ne\emp$ if $W',W''\in\Tf_x$ or
$W',W''\in\Tf^+_y$ because both $\Tf_x$ and $\Tf^+_y$ has a nonempty intersection. Let
$$
W'=\Wc_F[x,x',x'':\tlm],~~~x',x''\in\Mc~~~
\text{and}~~~ W''=\Wc_F[y,y',y'':\tlm]+rQ_0,~~~y',y''\in\Mc,
$$
be two arbitrary members of $\Tf_x$ and  $\Tf^+_y$ respectively. Then, thanks to assumption \rf{WNEW} of Theorem \reff{W-CR}, $W'\cap W''\ne\emp$.
\par Thus, the hypothesis of Lemma \reff{H-R} holds for
$\Tfw$. Therefore, this family has a common point proving the required property \rf{WDTC}.
\smsk
\par Let us prove parts {\it (ii)} and {\it (iii)}. We note that, thanks to property \rf{WDTC}, the mapping $\Tc=\Tc_{F,\tlm}$ satisfies the conditions of the hypothesis of part (b) of Proposition \reff{X2-C} with $\eta=\lambda$. This proposition tells us that property \rf{TC-NEMP} and inequality \rf{DRT-1} hold proving part {\it (ii)}. Furthermore, part (b) of Proposition \reff{X2-C} proves part {\it (iii)}.
\par The proof of Proposition \reff{WLS-T} is complete.\bx
\bsk

\par \underline{\sc The Third Step.} At this step we construct a Lipschitz selection $f$ of the set-valued mapping $F$ with Lipschitz seminorm at most $2\lambda+\tlm$.
\msk
\par We recall that the set-valued mapping $F^{[1]}[\cdot:\tlm]$ and its rectangular hull, the set-valued mapping $\Tc_{F,\tlm}=\HR[F^{[1]}[\cdot:\tlm]]$, are defined by formulae \rf{WF1-3L} and \rf{WTH-F1} respectively.
\par Part {\it (i)} of Proposition \reff{WLS-T} tells us that $\Tc_{F,\tlm}$ has a Lipschitz selection with Lipschitz seminorm at most $\lambda$. In other words, there exists a mapping $g:\Mc\to\RT$ such that
\bel{WG-HF}
g(x)\in \Tc_{F,\tlm}(x)=\HR[F^{[1]}[x:\tlm]]~~~~~\text{for every}~~~~~x\in\Mc,
\ee
and
\bel{G-LIP}
\|g(x)-g(y)\|\le \lambda\,\rho(x,y)~~~~~\text{for all}~~~~~x,y\in\Mc.
\ee

\begin{proposition}\lbl{LSEL-F} Let $g:\Mc\to\RT$ be an arbitrary Lipschitz selection of the set-valued mapping $\Tc_{F,\tlm}:\Mc\to\RCT$ with Lipschitz seminorm at most $\lambda$, i.e., a mapping satisfying conditions \rf{WG-HF} and \rf{G-LIP}.
\smsk
\par We define a mapping $f:\Mc\to\RT$ by letting
$$
f(x)=\Prm\left(g(x),F^{[1]}[x:\tlm]\right),~~~~~~x\in\Mc.
$$
(Recall that $\Prm(\cdot,S)$ is the operator of metric projection onto a convex closed $S\subset\RT$. See \rf{MPR}.)
\smsk
\par Then the following properties hold:
\smsk
\par ($\bigstar 1$) The mapping $f$ is well defined, i.e., $f(x)$ is a singleton for every $x\in\Mc$. In this case
$$
f(x)=\Prm\left(g(x),F^{[1]}[x:\tlm]\right)\in F^{[1]}[x:\tlm]\subset F(x)~~~~\text{for every}~~~~x\in\Mc,
$$
so that $f$ is a {\it selection of $F$ on $\Mc$};
\smsk
\par ($\bigstar 2$) The mapping $f:\Mc\to\RT$ is Lipschitz with Lipschitz seminorm $\|f\|_{\Lip(\Mc)}\le 2\lambda+\tlm$.
\end{proposition}
\par The proof of this proposition is based on a number of auxiliary results. The first of these is the following lemma.
\begin{lemma}\lbl{WMP-S} Let $S\subset\RT$ be a nonempty convex closed set.
\par Then for every point $a\in\HR[S]$ the metric projection $\Prm(a,S)$ is a singleton.
\par Furthermore, $\Prm(a,S)$ coincides with a vertex of the square $Q(a,\dist(a,S))$.

See Fig. 11.
\msk

\begin{figure}[H]
\hspace{25mm}
\includegraphics[scale=0.77]{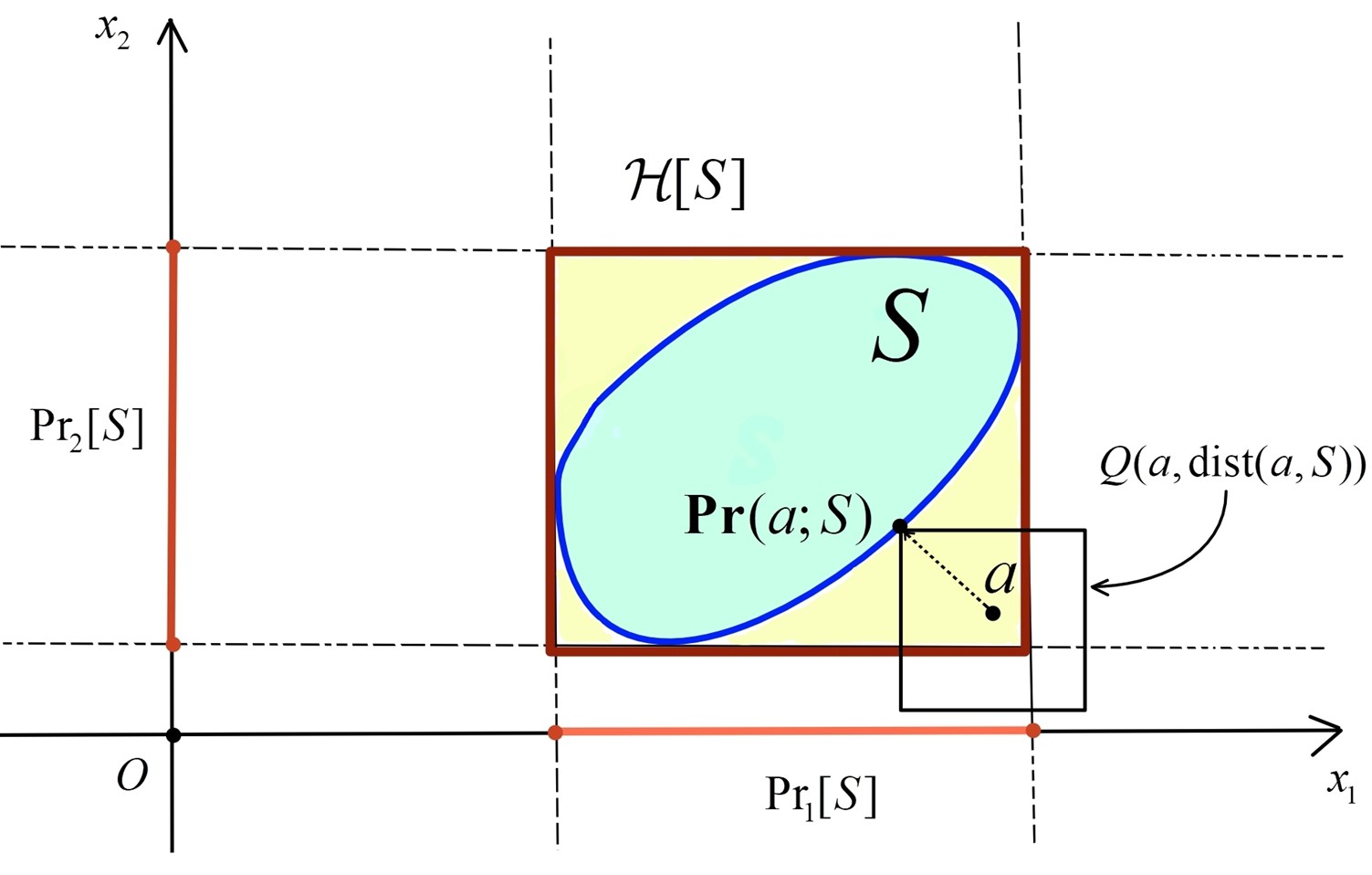}
\vspace*{3mm}
\caption{Lemma \reff{WMP-S}.}
\end{figure}

\end{lemma}

\par {\it Proof.} The proof of this lemma is given in \cite[p. 301]{S-2002}. See also \cite[p. 68]{S-2021-L}. For the convenience of the reader, we give it here.
\smsk
\par Clearly, if $a\in S$, nothing to prove. Suppose $a\notin S$ so that $r=\dist(a,S)>0$. Because $S$ is closed, $\Prm(a;S)\ne\emp$. Furthermore,
$$
\Prm(a;S)=S\cap Q=S\cap \partial Q
$$
where  $Q=Q(a,r)$.
\par Because $\Prm(a;S)$ is {\it a nonempty convex set} lying on the boundary of $Q$, it belongs to a certain side of the square $Q$. In other words, there exist two distinct vertices of $Q$, say $A$ and $B$, such that $\Prm(a;S)\subset [A,B]$.
\par Let us prove that
\bel{AB-V}
\text{either}~~~~\Prm(a;S)=\{A\}~~~~\text{or}~~~~ \Prm(a;S)=\{B\}.
\ee
\par Indeed, otherwise there exists a point
$$
p\in (A,B)\cap\Prm(a;S).
$$
\par Let $\ell$ be the straight line passing through $A$ and $B$. Clearly, $\ell$ is parallel to a coordinate axis. Let $H_1,H_2$ be the closed half-planes determined by $\ell$. (Thus $\ell=H_1\cap H_2$ and $H_1\cup H_2=\RT$.) Clearly, $Q$ is contained in one of these half-planes, say in $H_1$. Because $\dist(a,\ell)=r>0$, the point $a\in H^{int}_1$ where $H^{int}_1$ denotes the interior of $H_1$.

\smsk
\par Prove that in this case $S\subset H_2$, i.e., the straight line $\ell$ separates (not strictly) the square $Q$ and the set $S$. Indeed, suppose that there exists  a point $b\in S\cap H^{int}_1$. Then also $(p,b]\subset H^{int}_1$ because $p\in\partial H_1=\ell$. But
$$
p\in(A,B)~~~\text{so that}~~~(p,b]\cap Q^{int}\ne\emp.
$$
\par On the other hand, because $S$ is convex and $p\in\partial S$, the interval $(p,b]\subset S$ proving that $S\cap Q^{int}\ne\emp$. But $S\cap Q \subset \partial Q$, a contradiction. See Fig. 12.
\msk

\begin{figure}[H]
\hspace{35mm}
\includegraphics[scale=0.65]{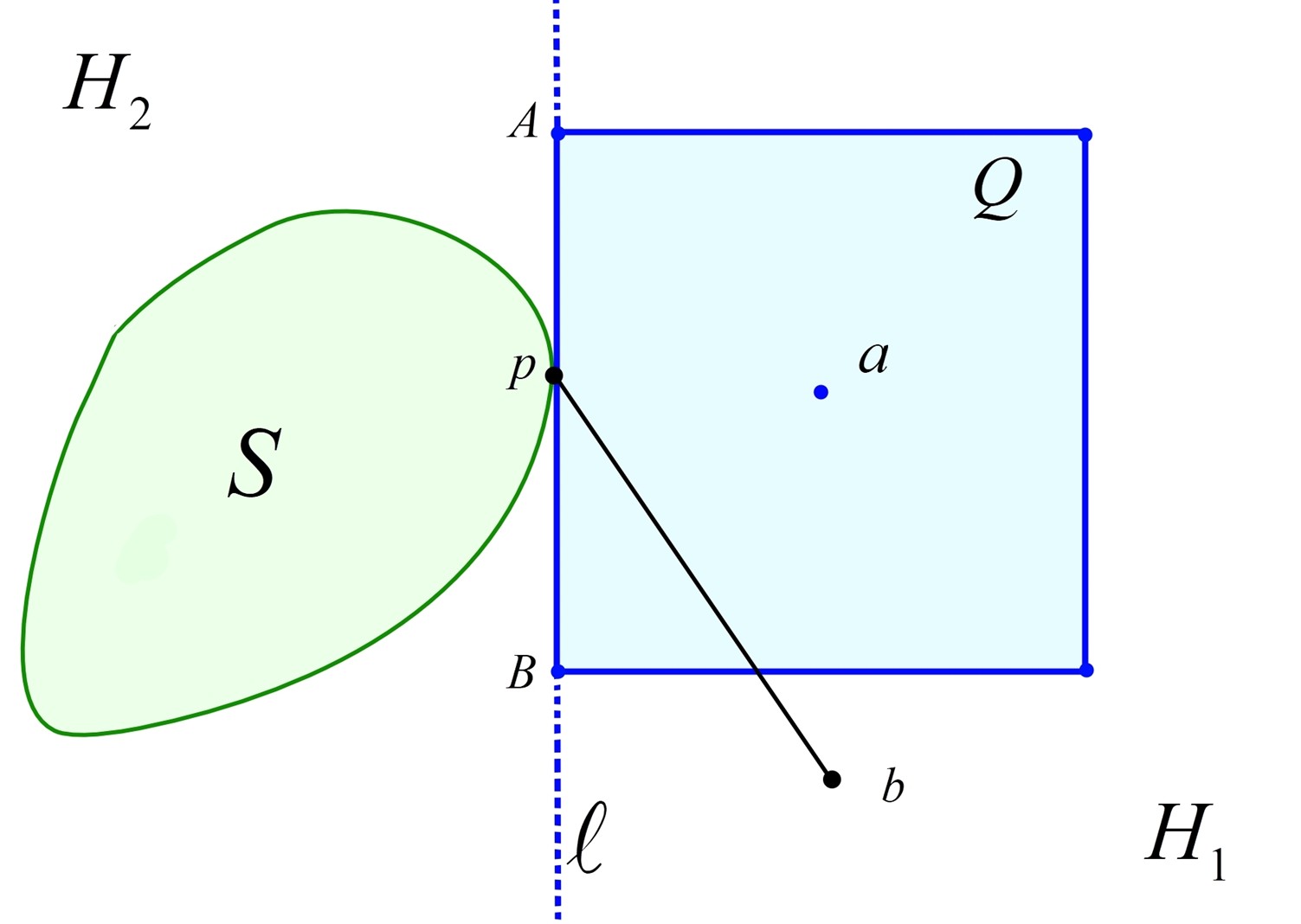}
\vspace*{3mm}
\caption{Lemma \reff{WMP-S}.}
\end{figure}

\par Thus, $S\subset H_2$ and $Q\subset H_1$. But $a\in H^{int}_1$ so that $a\notin H_2$. Clearly, $H_2\in\RCT$, i.e., $H_2$ is an (unbounded) rectangle. Therefore $\HR[S]\subset H_2$, see definition \rf{HRS}. Therefore, thanks to the lemma's hypothesis, $a\in\HR[S]\subset H_2$, a contradiction.
\par This contradiction implies \rf{AB-V} completing the proof of the lemma.\bx
\smsk
\par Clearly, this lemma implies the statement
($\bigstar 1$) of Proposition \reff{LSEL-F}.
\smsk
\par Let us prove the statement ($\bigstar 2$) which is equivalent to the inequality
\bel{LIP-FL}
\|f(x)-f(y)\|\le (2\lambda+\tlm)\,\rho(x,y)~~~
\text{for every}~~~x,y\in\Mc.
\ee
\par The proof of this inequality relies on a number of auxiliary results which we present in the next section.

\SECT{4. Proof of the key theorem: the final step.}{4}

\addtocontents{toc}{4. Proof of the key theorem: the final step. \hfill \thepage\par\VST}


\indent

\begin{lemma}\lbl{AB-PR} Let $A,B\subset \R^{2}$ be nonempty convex closed sets such that $A\subset B$, and let $a\in \HR[A]$.
\par Then $\Prm(a,A)$ and $\Prm(a,B)$ are singletons having the following properties:
\msk
\par (i) $\Prm(a,B)\in [\Prm(a,A),a]$. (See notation \rf{LSEG}.)
\smsk
\par  (ii) The following equality holds:
$$
\|\Prm(a,A)-\Prm(a,B)\|=\dist(a,A)-\dist(a,B).
$$

See Fig. 13.

\begin{figure}[H]
\hspace{28mm}
\includegraphics[scale=0.77]{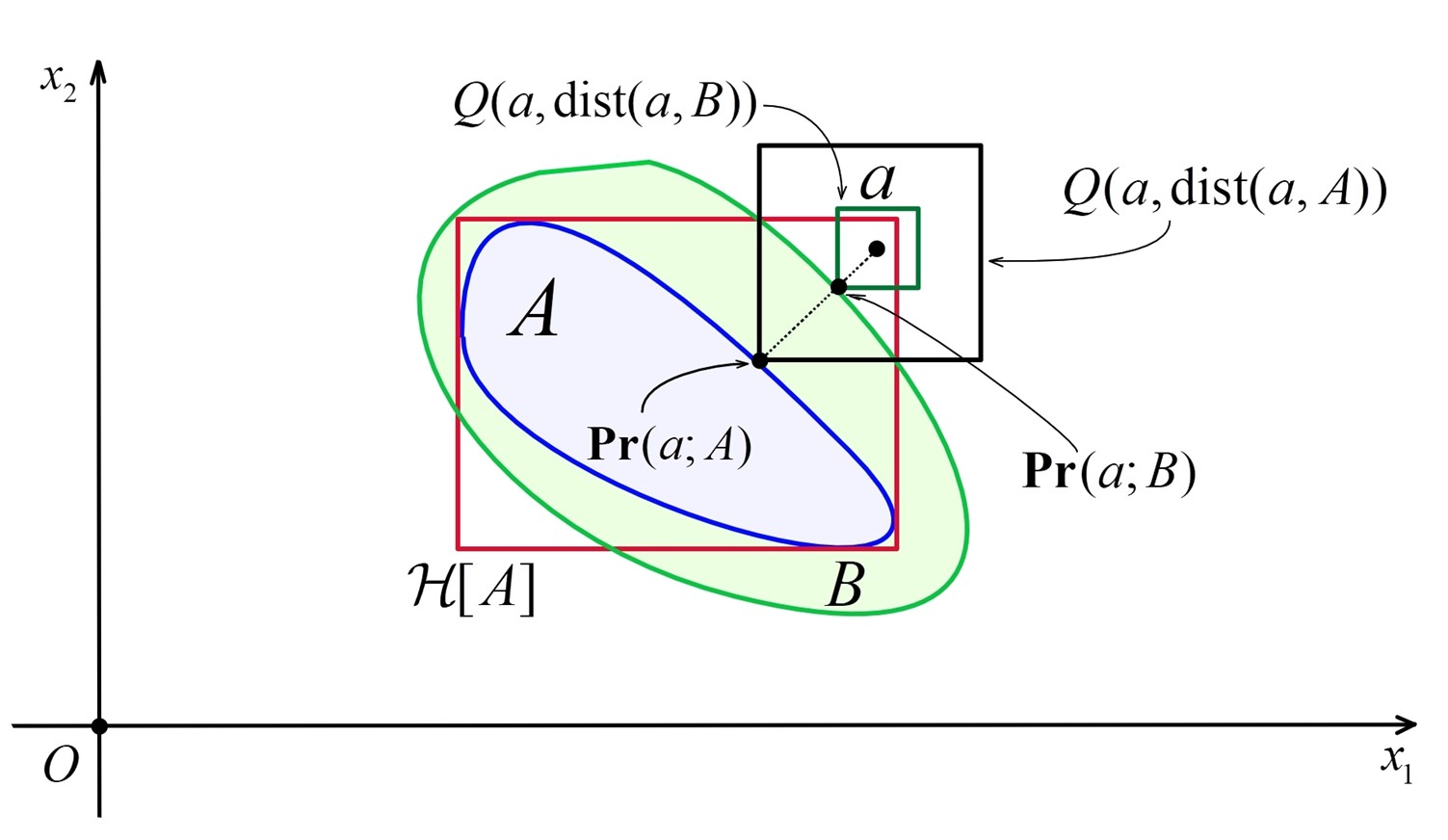}
\hspace*{10mm}
\caption{Lemma \reff{AB-PR}.}
\end{figure}
\end{lemma}
\par {\it Proof.} The proof of the lemma is given in \cite[p. 302]{S-2002}. See also \cite[p. 69]{S-2021-L}. For the convenience of the reader, we present it here.
\smsk
\par First, we note that if $a\in B$, the statement of the lemma is immediate from Lemma \reff{WMP-S}.
\par Suppose that $a\notin B$. In this case, Lemma \reff{WMP-S} tells us that $\Prm(a;A)$ is one of the vertices of the square $Q(a,r)$ with $r=\dist(a,A)>0$. Because $A\subset B$, the point $a\in \HR[B]$ so that, thanks to Lemma \reff{WMP-S}, $\Prm(a;B)$ is a vertex of the square $Q(a,\alpha)$ where $\alpha=\dist(a,B)>0$.
\par Using a suitable shift and dilation, without loss of generality, we can assume that
$$
a=(0,0),~~r=\dist(a,A)=1,~~~\text{and}~~~ \Prm(a;A)=(1,1).
$$
Clearly, in this case $0<\alpha\le 1$. Furthermore, under these conditions, the statement of the lemma is equivalent to the property
\bel{P-AL}
\Prm(a;B)=(\alpha,\alpha).
\ee
\par Suppose that this property does not hold, i.e.,
$\Prm(a;B)\in\{(\alpha,-\alpha),(-\alpha,\alpha),
(-\alpha,-\alpha)\}$.
\par In order to get a contradiction, we construct a straight line $\ell_{A}$ which passes through $(1,1)$ and separates (not strictly) the square $Q(a,r)=[-1,1]^{2}$ and $A$. This line determines two closed half-planes, $S^{+}_{A}$ and $S^{-}_{A}$, with the common boundary (i.e., the line $\ell_A$) such that $\RT=S^{+}_{A}\cup S^{-}_{A}$.
\par See Fig. 14.

\begin{figure}[H]
\hspace{20mm}
\includegraphics[scale=0.80]{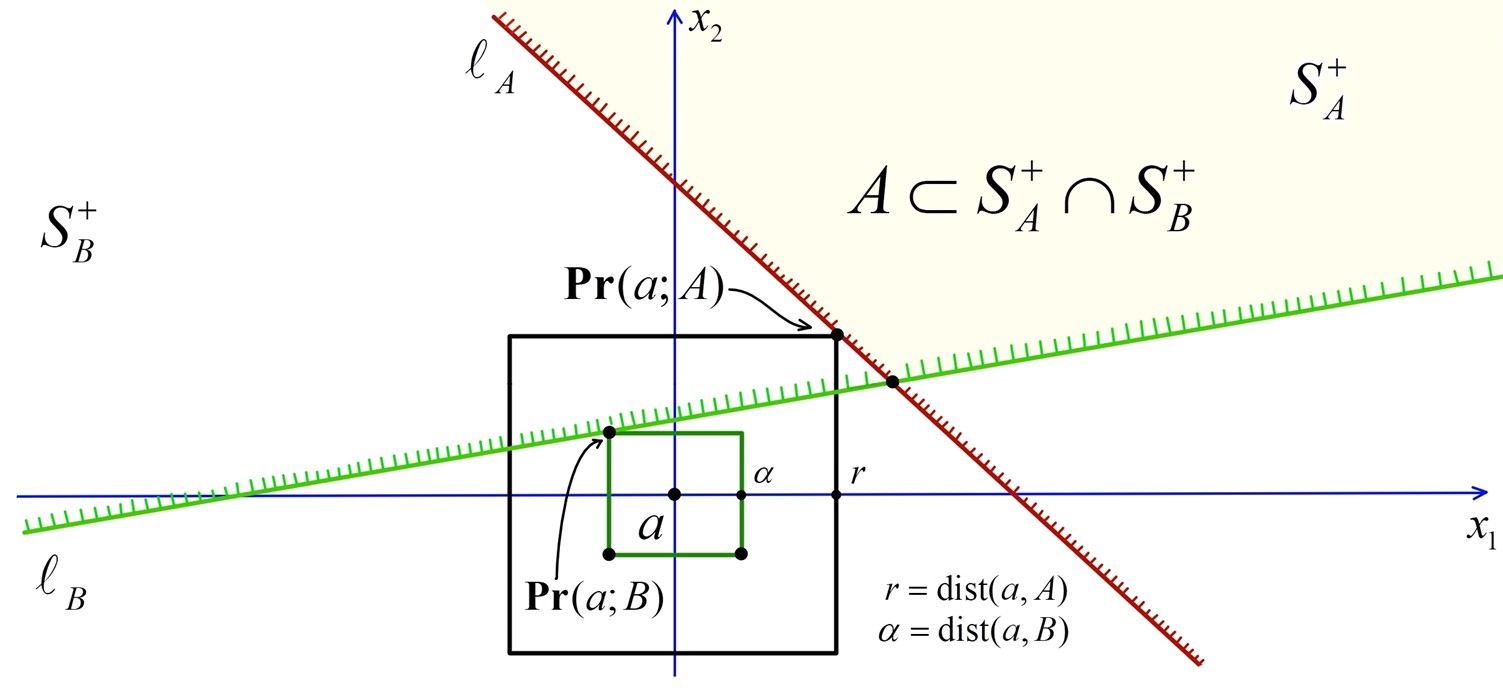}
\caption{Lemma \reff{AB-PR}.}
\end{figure}

\par One of them, say $S^{+}_{A}$, contains $A$, so that $S^{-}_{A}\supset Q(a,r)$. We know that $S^{+}_{A}$ contains $(1,1)$ and does not contain intrinsic points of the square $[-1,1]^{2}$, so that $Q(a,r)\cap\ell_A=(1,1)$. Therefore, the half-plane $S^{+}_{A}$ admits the following representation:
\bel{SA-R2}
S^{+}_{A}=\{x=(x_1,x_2)\in \R^{2}:
(x_1-1)\,h_1+(x_2-1)\,h_2\ge 0\}
\ee
where $h_{1},h_{2}>0$ are certain numbers.
\smsk
\par Let us assume that $\Prm(a;B)=(-\alpha,\alpha)$ and show that this assumption leads to a contradiction. We let $\ell_B$ denote a straight line which passes through the point $(-\alpha,\alpha)$ and separates the square $Q(a,\dist(a,B))=[-\alpha,\alpha]^2$ and the set $B$. Let $S^{+}_{B}$ be the one of the two half-planes determined by $\ell_B$ which contains $B$. Then the other half-plane, $S^{-}_{B}$, contains $Q(a,\dist(a,B))$, and $S^{+}_{B}\cap S^{-}_{B}=\ell_B$.
\smsk
\par We know that $S^{+}_{B}$ contains the point $(-\alpha,\alpha)$ on its boundary and does not contain intrinsic points of the square $[-\alpha,\alpha]^2$. Therefore, this half-plane can be represented in the form
\bel{SB-N}
S^{+}_{B}=\{(x_1,x_2)\in \R^{2}:
-(x_1+\alpha)\,s_{1}+(x_2-\alpha)\,s_2\ge 0\}
\ee
with certain $s_1,s_2>0$.
\par Thus, $A\subset S^{+}_{A}$ and $A\subset B\subset S^{+}_{B}$, so that $A\subset S^{+}_{A}\cap S^{+}_{B}$ proving that for every $x=(x_1,x_2)\in A$ we have
\bel{S-ABA}
(x_1-1)\,h_1+(x_2-1)\,h_2\ge 0~~~\text{and}~~~
-(x_1+\alpha)\,s_{1}+(x_2-\alpha)\,s_2\ge 0.
\ee
See \rf{SA-R2} and \rf{SB-N}. Note also that since  $S^{+}_{A}\cap S^{+}_{B}\supset A\ne\emp$, we have $h_2+s_2>0$.
\par Let us prove that inequalities \rf{S-ABA} imply the following inclusion:
\bel{X2-A}
A\subset \Hc_\alpha=\{x=(x_1,x_2)\in\RT:x_2\ge\alpha\}.
\ee
\par Indeed, it is easy to see that from \rf{S-ABA} we have
$$
x_2-\alpha\ge \frac{s_1((1+\alpha)h_1+(1-\alpha)h_2))}
{s_1h_2+s_2h_1}\ge 0,~~~~~~x=(x_1,x_2)\in A,
$$
proving \rf{X2-A}.
\par Let us note that $\Hc_\alpha\in\RCT$ so that  $\Hc[A]\subset\Hc_\alpha$. Therefore, thanks to the lemma's assumption, $a=(0,0)\in \Hc_\alpha$. But $\alpha>0$ so that $a=(0,0)\notin \Hc_\alpha$, a contradiction.
\par In a similar way, we get a contradiction provided
$\Prm(a;B)=(\alpha,-\alpha)$ or $\Prm(a;B)=(-\alpha,-\alpha)$ proving the required property \rf{P-AL} and the lemma.\bx
\begin{lemma}\lbl{D-23} (i) Let $u\in\Mc$, and let $a\in\HR[F^{[1]}[u:\tlm]]$. Then
\bel{DS-F1}
\dist(a,F^{[1]}[u:\tlm])=\max_{z\in\Mc}\,\dist(a,F(z)+
\tlm\,\rho(u,z)Q_0)
=\max_{z\in\Mc}\,[\dist(a,F(z))-\tlm\,\rho(u,z)]_+\,;
\ee
\par (ii) Let $u,v\in\Mc$, and let $a\in\HR[F^{[1]}[u:\tlm]]$ and $b\in\HR[F^{[1]}[v:\tlm]]$.
Then
\bel{DR-X}
|\dist(a,F^{[1]}[u:\tlm])-\dist(b,F^{[1]}[v:\tlm])|
\le \|a-b\|+\tlm\,\rho(u,v).
\ee
\end{lemma}
\par {\it Proof.} {\it (i)} Let
$$
A=F^{[1]}[u:\tlm]
$$
and, given $z\in\Mc$, let $A_z=F(z)+\tlm\,\rho(u,z)Q_0$. Then, thanks to \rf{WF1-3L}, $A=\cap\{A_z:z\in\Mc\}$. Our goal is to prove that
\bel{O-G}
\dist(a,A)=\max\{\dist(a,A_z):z\in\Mc\}~~~\text{provided}~~~
a\in\HR[A].
\ee
\par Lemma \reff{WF1-NE} tells us that $A$ is a nonempty convex and closed subset of $\RT$. Because $A\subset A_z$ for each $z\in\Mc$, the left hand side of the above equality majorizes its right hand side.
\smsk
\par Prove the converse inequality. If $a\in A$, nothing to prove.
\par Let $a\notin A$, and let $\ve\in(0,\dist(a,A))$ be arbitrary. We know that $a\in\HR[A]$ so that, thanks to Lemma \reff{AB-PR}, $\Prm(a,A)$ is a singleton. See Fig. 15.

\bsk
\begin{figure}[H]
\hspace{5mm}
\includegraphics[scale=0.7]{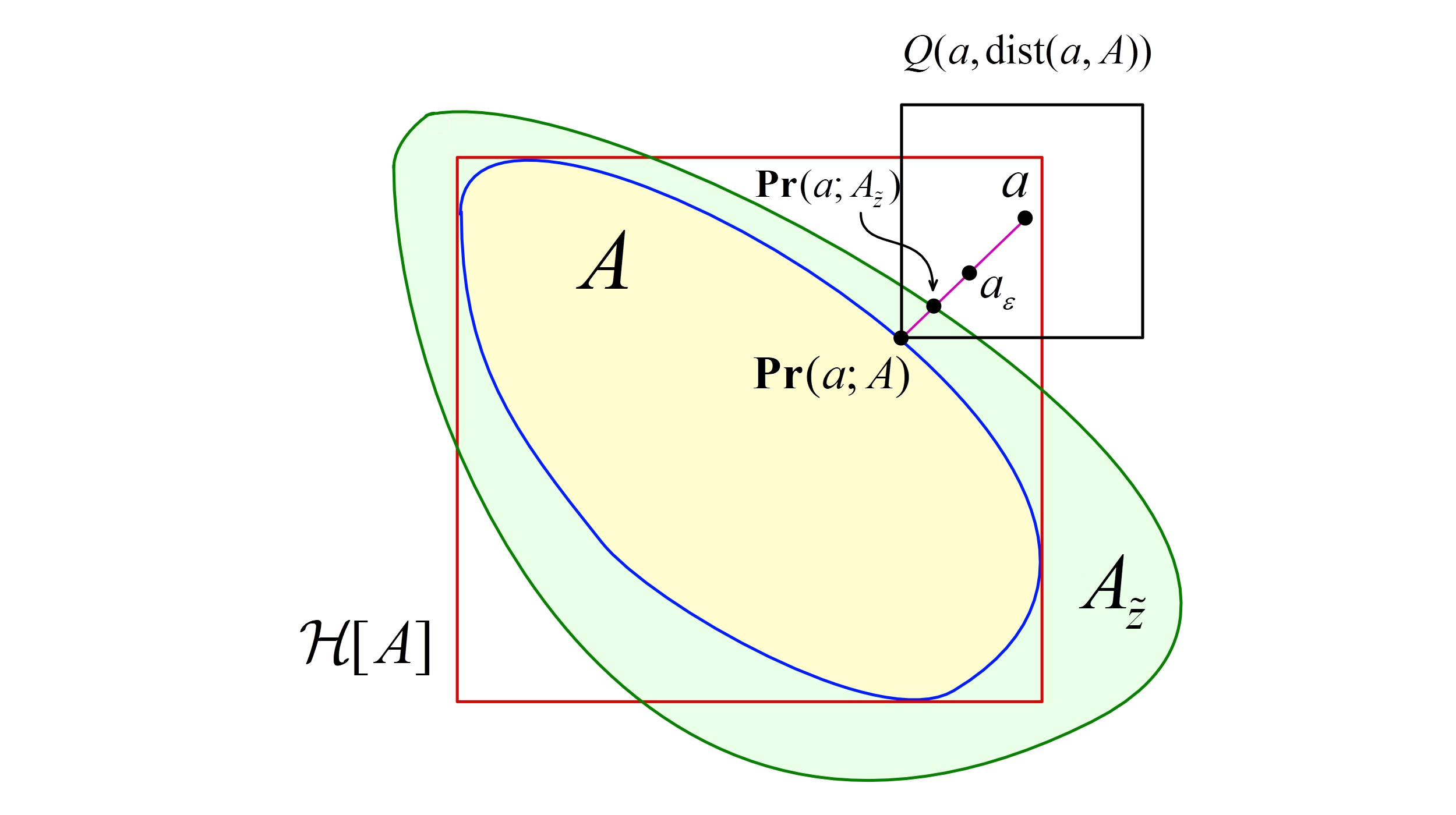}
\caption{Lemma \reff{D-23}.}
\end{figure}

\par We let $a_{\ve}$ denote a point on the interval $(\Prm(a,A),a]$ such that $\|a_{\ve}-\Prm(a,A)\|<\ve$. Because $a_\ve\notin A$ and $A=\cap\{A_z:z\in\Mc\}$, there exists an element $\tz\in\Mc$ such that $a_\ve\notin A_{\tz}$. Note that $A\subset A_{\tz}$. Lemma \reff{AB-PR} tells us that in this case $\Prm(a,A_{\tz})$ is a  singleton such that $\Prm(a,A_{\tz})\in[\Prm(a,A),a]$.
\par Then $\Prm(a,A_{\tz})\in [\Prm(a,A),a_\ve]$; otherwise
$a_\ve\in[\Prm(a,A),\Prm(a,A_{\tz})]\subset A_{\tz}$, a contradiction.
\par This proves that $\|\Prm(a,A)-\Prm(a,A_{\tz})\|<\ve$.
Hence,
\be
\dist(a,A)&=&\|a-\Prm(a,A)\|\le \|a-\Prm(a,A_{\tz})\|+
\|\Prm(a,A_{\tz})-\Prm(a,A)\|\nn\\
&\le&
\dist(a,A_{\tz})+\|a_\ve-\Prm(a,A)\|\le \dist(a,A_{\tz})+\ve.
\nn
\ee
\par Since $\ve>0$ can be chosen as small as desired, this implies the required inequality \rf{O-G} proving part {\it (i)} of the lemma.
\msk
\par {\it (ii)} Let $A=F^{[1]}[u:\tlm]$ and $B=F^{[1]}[v:\tlm]$.
Then, thanks to \rf{DS-F1},
\be
|\dist(a,A)-\dist(a,B)|&=&
|\sup_{z\in\Mc}\,[\dist(a,F(z))-\tlm\rho(u,z)]_+
-\sup_{z\in\Mc}\,[\dist(a,F(z))-\tlm\rho(v,z)]_+|
\nn\\
&\le&
\sup_{z\in\Mc}|\,[\dist(a,F(z))-\tlm\rho(u,z)]_+
-[\dist(a,F(z))-\tlm\rho(v,z)]_+|
\nn\\
&\le&
\tlm\,\sup_{z\in\Mc}|\rho(u,z)-\rho(v,z)|
\nn
\ee
so that, thanks to the triangle inequality,
\bel{D-AB1}
|\dist(a,A)-\dist(a,B)|\le\tlm\,\rho(u,v).
\ee
\par Next,
$$
|\dist(a,A)-\dist(b,B)|\le
|\dist(a,A)-\dist(a,B)|+|\dist(a,B)-\dist(b,B)|.
$$
Because $\dist(\cdot,B)$ is a Lipschitz function, from this and \rf{D-AB1}, we have \rf{DR-X} completing the proof of the lemma.\bx
\msk
\par Let $\delta\ge 0$, and let
\bel{H12-DF}
H_1~~\text{and}~~~H_2~~~\text{be two half-planes with}~~~ \dist(H_1,H_2)\le\delta.
\ee
\par Let
$$
\ell_i=\partial H_i~~~\text{be the boundary of the half-plane}~~~H_i,~~~i=1,2.
$$
\par Let us represent the half-planes $H_i$, $i=1,2$, in the form
$$
H_i=\{u\in\RT:\ip{\lh_i,u}+\alpha_i\le 0\}
$$
where $\lh_i$ is a unit vector and $\alpha_i\in\R$. (Recall that, given points $a=(a_1,a_2),b=(b_1,b_2)\in\RT$, by
$\ip{a,b}=a_1 b_1+a_2 b_2$ we denote the standard inner product in $\RT$.)
\smsk
\par Thus the vector
\bel{H-OL}
\lh_i~~\text{is directed outside of} ~~H_i~~\text{and}~~~\lh_i\perp \ell_i,~~~i=1,2.
\ee

\begin{proposition}\lbl{TWO-HP} Let
$$
S_1=H_1\cap(H_2+\delta Q_0)~~~~\text{and}~~~ S_2=H_2\cap(H_1+\delta Q_0).
$$
\par Given $a_1,a_2\in\RT$ suppose that for every $i=1,2$,
\bel{A-I}
a_i\in\HR[S_i]
\ee
and
\bel{A-II}
\Prm(a_i,H_i)\in S_i.
\ee
\par Then the following inequality
\bel{PRA-D}
\|\Prm(a_1,H_1)-\Prm(a_2,H_2)\|\le \delta+2\|a_1-a_2\|
\ee
holds.

\par See Fig. 16.
\end{proposition}

\begin{figure}[H]
\hspace{10mm}
\includegraphics[scale=0.7]{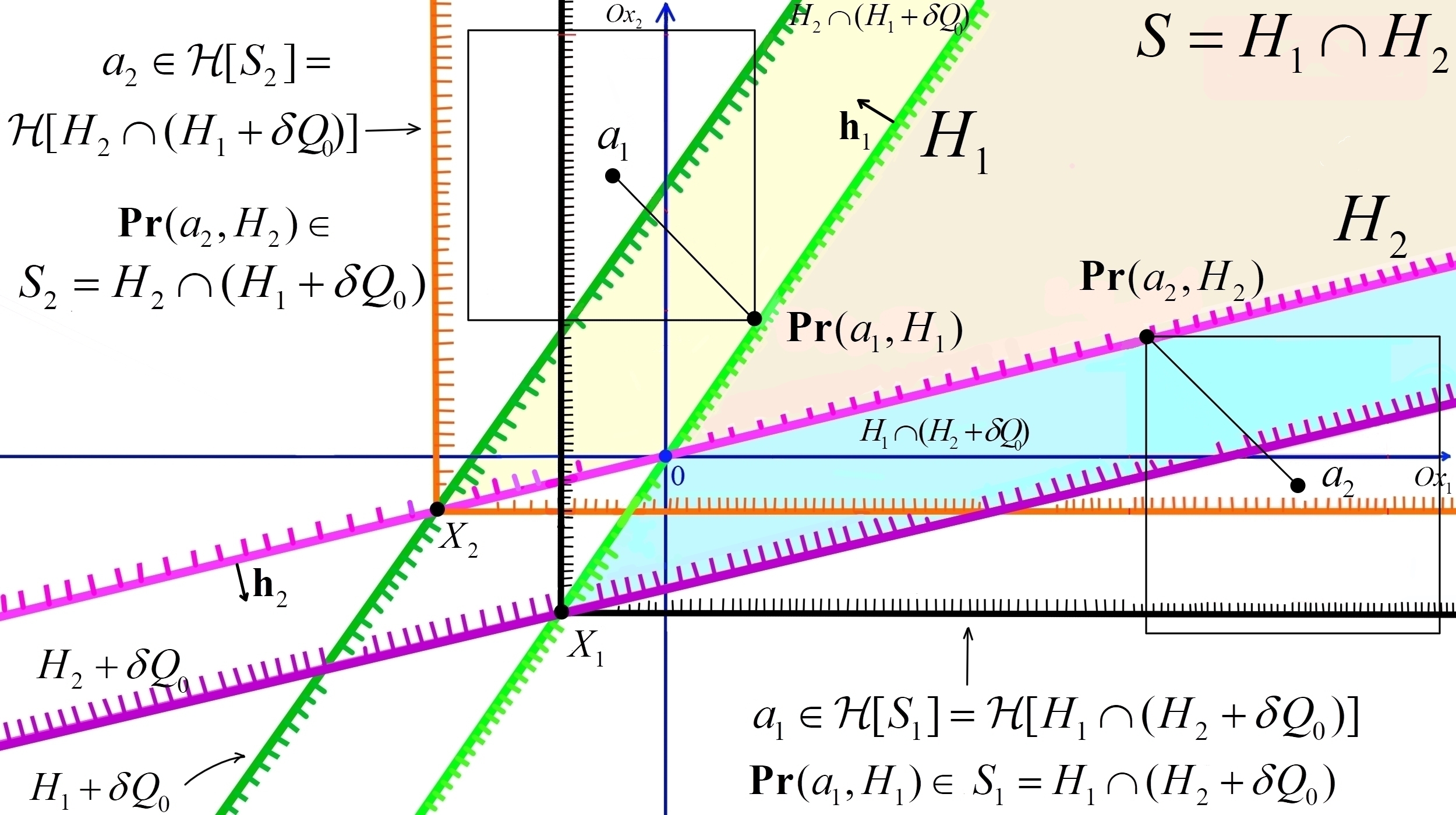}
\caption{Proposition \reff{TWO-HP}.}
\end{figure}
\msk
\msk
\par {\it Proof.} The proof of this proposition relies on several auxiliary statements.
\par Let us formulate the first of them. Let
$$
S=H_1\cap H_2.
$$
Recall that
\bel{S12-D}
S_1=H_1\cap(H_2+\delta Q_0)~~~~\text{and}~~~~
S_2=H_2\cap(H_1+\delta Q_0).
\ee
\par We know that $\dist(H_1,H_2)\le\delta$ so that $S_1\ne\emp$ and $S_2\ne\emp$. Furthermore, thanks to \rf{A-I}, \rf{A-II} and \rf{S12-D},
\bel{A12-S}
a_i\in\HR[S_i]~~~~\text{and}~~~\Prm(a_i,H_i)\in S_i,
~~~i=1,2.
\ee
\begin{lemma}\lbl{L1-P} Both $\Prm(a_1,H_1)$ and $\Prm(a_2,H_2)$ are singletons. Furthermore,
$\Prm(a_i,H_i)=\Prm(a_i,S_i)$ for every $i=1,2$, and the following inequality holds:
$$
|\dist(a_1,H_1)-\dist(a_2,H_2)|\le
\delta+\|a_1-a_2\|
$$
\end{lemma}
\par {\it Proof.} Thanks to \rf{A12-S}, the point $a_i\in \HR[S_i]$, so that  $a_i\in\HR[H_i]$ because $S_i\subset H_i$, $i=1,2$, see \rf{S12-D}. Therefore, thanks to Lemma \reff{WMP-S},
$$
\Prm(a_i,H_i)~~~\text{is a singleton for every}~~~i=1,2.
$$
\par Let us see that
$$
\dist(a_i,H_i)=\dist(a_i,S_i)~~~\text{and}~~~
\Prm(a_i,H_i)=\Prm(a_i,S_i),~~~~~i=1,2.
$$
\par Indeed, thanks to \rf{A-II}, $\Prm(a_i,H_i)\in S_i$. But $S_i\subset H_i$ so that
$$
\dist(a_i,H_i)\le\dist(a_i,S_i)\le \|a_i-\Prm(a_i,H_i)\|=\dist(a_i,H_i),~~~i=1,2.
$$
Furthermore, let $R=\dist(a_i,H_i)=\dist(a_i,S_i)$. Then
$$
\Prm(a_i,S_i)=Q(a_i,R)\cap S_i\subset Q(a_i,R)\cap H_i=\Prm(a_i,H_i).
$$
But $\Prm(a_i,H_i)$ is a singleton proving that $\Prm(a_i,H_i)=\Prm(a_i,S_i)$.
\bsk

\par Clearly,
\bel{HD-DL}
\dhf(S_1,S_2)=\dhf(H_1\cap[H_2+\delta Q_0],
H_2\cap[H_1+\delta Q_0])\le \delta.
\ee
\par Indeed, let $x\in H_1\cap[H_2+\delta Q_0]$. Then $x\in H_1$ and $\|x-y\|\le\delta$ for some $y\in H_2$.
Hence, $y\in H_2\cap[H_1+\delta Q_0]$ so that $\dist(x,H_2\cap[H_1+\delta Q_0])\le\delta$ proving \rf{HD-DL}.
\smsk
See \rf{HD-DF}. Therefore,
$$
|\dist(a_1,S_1)-\dist(a_1,S_2)|
\le\dhf(S_1,S_2)\le \delta.
$$
\par Note also that the function $\dist(\cdot,S_2)$ is Lipschitz. Therefore,
\be
|\dist(a_1,H_1)-\dist(a_2,H_2)|&=&
|\dist(a_1,S_1)-\dist(a_2,S_2)|
\nn\\
&\le& |\dist(a_1,S_1)-\dist(a_1,S_2)|+
|\dist(a_1,S_2)-\dist(a_2,S_2)|
\nn\\
&\le&
\delta+\|a_1-a_2\|
\nn
\ee
proving the lemma.\bx
\msk
\begin{remark}\lbl{PR-EL} {\em Recall that $\ell_i$ is the boundary of the half-plane $H_i$, $i=1,2$. Note that the assumption \rf{A-I} implies the following property: if
$$
\ell_i\parallel Ox_j
$$
for some $i,j\in\{1,2\}$, then $a_i\in H_i$. Indeed, suppose that this is true, say for $i=1$, i.e., either $\ell_1\parallel Ox_1$ or $\ell_1\parallel Ox_2$. Then $H_1\in\RCT$, i.e., $H_1$ is a rectangular. Then $\HR[H_1]=H_1$ so that, thanks to \rf{A-I}, $a_1\in \HR[S_1]\subset \HR[H_1]=H_1$.\rbx}
\end{remark}

\par Given a vector $\lh=(h_1,h_2)\in\RT$, we set
\bel{SN-H}
\SN(\lh)=(\sign(h_1),\sign(h_2)).
\ee
\par Let us note several useful properties of metric projections in the space $\LTI=(\RT,\|\cdot\|)$.
\par Let $\alpha\in \R$ and let $\lh=(h_1,h_2)$ be a unit vector. Let
$$
H=\{u\in\RT:\ip{\lh,u}+\alpha\le 0\}
$$
be a half-plane determined by $\alpha$ and $\lh$.
\begin{lemma}\lbl{SN-D} (i) For every $g=(g_1,g_2)\in\RT$, we have
\bel{DA-HN}
\dist(g,H)=[\ip{\lh,g}+\alpha]_+\,/\,\|\lh\|_1.
\ee
(Recall that $\|\lh\|_1=|h_1|+|h_2|$.)
\msk
\par (ii) Suppose that $h_1,h_2\ne 0$ and $g\notin H$. Let $f=(f_1,f_2)=\Prm(g,H)$. Then $f$ is a singleton having the following properties:
\bel{NP-HPN}
f_1=\frac{h_2(g_1-g_2)-\alpha}{h_1+h_2}~~~~\text{and}~~~~
f_2=\frac{-h_1(g_1-g_2)-\alpha}{h_1+h_2}~~~~~\text{if}~~~~
h_1\cdot h_2>0,
\ee
and
\bel{NP-HN}
f_1=\frac{-h_2(g_1+g_2)-\alpha}{h_1-h_2}~~~~\text{and}~~~~
f_2=\frac{h_1(g_1+g_2)+\alpha}{h_1-h_2}~~~~~\text{if}~~~~
h_1\cdot h_2<0.
\ee
\par Furthermore, the following equality holds:
\bel{FG-DS}
f=g-\dist(g,H)\,\SN(\lh).
\ee
See \rf{SN-H}.
\end{lemma}
\par {\it Proof.} Formula \rf{DA-HN} is a well known result from the functional analysis.
\smsk
\par Prove (ii). Because $h_1\ne 0$ and $h_2\ne 0$, the rectangular hull of $H$ coincides with $\RT$. Therefore, thanks to Lemma \reff{WMP-S}, $\Prm(g,H)$ is a singleton.
\par Because $g\notin H$, we have
\bel{G-NH}
\ip{\lh,g}+\alpha>0.
\ee
\par If $h_1\cdot h_2>0$ then
$$
\ip{\lh,f}+\alpha=h_1 f_1+h_2f_2+\alpha=
h_1\cdot \frac{h_2(g_1-g_2)-\alpha}{h_1+h_2}+
h_2 \cdot\frac{-h_1(g_1-g_2)-\alpha}{h_1+h_2}+\alpha=0,
$$
and if $h_1\cdot h_2<0$ then
$$
\ip{\lh,f}+\alpha=
h_1\cdot \frac{-h_2(g_1+g_2)-\alpha}{h_1-h_2}+
h_2 \cdot\frac{h_1(g_1+g_2)+\alpha}{h_1-h_2}+\alpha=0
$$
as well, proving that $f\in\partial H$. Furthermore,
if $h_1\cdot h_2>0$ then
$$
f_1-g_1=\frac{h_2(g_1-g_2)-\alpha}{h_1+h_2}-g_1=
-\frac{h_1g_1+h_2g_2}{h_1+h_2}
=-\frac{\ip{\lh,g}}{\|\lh\|_1}\sign(h_1)
$$
and
$$
f_2-g_2=\frac{h_1(g_1+g_2)+\alpha}{h_1-h_2}-g_2=
\frac{h_1g_1+h_2g_2}{h_1-h_2}
=-\frac{\ip{\lh,g}}{\|\lh\|_1}\sign(h_2).
$$
\par The reader can easily see that the same formulae hold provided $h_1\cdot h_2<0$. From this, \rf{DA-HN} and \rf{G-NH}, we have
$$
f-g=-\frac{\ip{\lh,g}+\alpha}{\|\lh\|_1}
(\sign h_1,\sign h_2)=
-\dist(g,H)\,(\sign h_1,\sign h_2)=-\dist(g,H)\SN(\lh).
$$
\par In particular, $\|f-g\|=\dist(g,H)$ proving that $f=\Prm(g,H)$ and the formula \rf{FG-DS}.
\par The proof of the lemma is complete.\bx

\begin{lemma}\lbl{L3-P} Inequality \rf{PRA-D} holds in the following two cases:
\par (i) either $a_1\in H_1$ or $a_2\in H_2$ and (ii) $\SN(\lh_1)=\SN(\lh_2)$.
\end{lemma}
\par {\it Proof.} Clearly, \rf{PRA-D} holds provided both $a_1\in H_1$ and $a_2\in H_2$. Therefore, without loss of generality, we may assume that $a_1\notin H_1$. Remark \reff{PR-EL} tells us that in this case $\ell_1\nparallel Ox_j$, $j=1,2$. But $\lh_1\perp \ell_1$ so that $\lh_1$ has non-zero coordinates.
\par Lemma \reff{SN-D} tells us that in this case
$$
\Prm(a_1,H_1)=a_1-\dist(a_1,H_1)\,\SN(\lh_1).
$$
\par Note that for the same reason, the formula
$$
\Prm(a_2,H_2)=a_2-\dist(a_2,H_2)\,\SN(\lh_2)
$$
holds provided $a_2\notin H_2$. It is also clear that this formula holds whenever $a_2\in H_2$.
\par Therefore, both in case (i) (i.e., if $a_2\in H_2$) and in case (ii) ($\SN(\lh_1)=\SN(\lh_2)$), we have
$$
\Prm(a_2,H_2)=a_2-\dist(a_2,H_2)\,\SN(\lh_1).
$$
\par Hence,
\be
\|\Prm(a_1,H_1)-\Prm(a_2,H_2)\|
&=&
\|a_1-\dist(a_1,H_1)\,\SN(\lh_1)-
(a_2-\dist(a_2,H_2)\,\SN(\lh_1))\|
\nn\\
&\le&
\|a_1-a_2\|+|\dist(a_1,H_1)-\dist(a_2,H_2)|.
\nn
\ee
\par Combining this inequality with Lemma \reff{L1-P}, we obtain the required inequality \rf{PRA-D}.\bx
\msk
\par Lemma \reff{L3-P} implies the following important statement:
$$
\text{\it Inequality \rf{PRA-D} holds provided either}~~a_1\in H_1~~\text{\it or}~~a_2\in H_2.
$$
\par Therefore, everywhere below, in the proof of inequality \rf{PRA-D}, we will assume that
\bel{A12-HN}
a_1\notin H_1~~~\text{and}~~~a_2\notin H_2.
\ee
\par Furthermore, thanks to Remark \reff{PR-EL}, in what follows, we may assume that
\bel{L-NP}
\ell_i\nparallel Ox_j~~~\text{for every}~~i,j=1,2.
\ee

\begin{remark} {\em We recall the representation of $H_1$ and $H_2$ in the form
$$
H_i=\{u\in\RT:\ip{\lh_i,u}+\alpha_i\le 0\},~~~i=1,2,
$$
where each $\lh_i$ is a unit vector and $\alpha_i\in\R$.
\par Thanks to \rf{H-OL}, $\lh_i\perp \ell_i\,(=\partial H_i)$, so that from \rf{L-NP} we have $\lh_i\nparallel Ox_j$ for every $i,j=1,2$. In particular, each $\lh_i$, $i=1,2$, has non-zero coordinates.\rbx}
\end{remark}

\begin{lemma}\lbl{L-PH2} Inequality \rf{PRA-D} holds provided
$$
\dist(a_i,H_1)+\dist(a_i,H_2)\le\delta~~~
\text{for some}~~~i\in\{1,2\}.
$$
\end{lemma}
\par {\it Proof.} For instance, suppose that
\bel{SD-A1}
\dist(a_1,H_1)+\dist(a_1,H_2)\le\delta,
\ee
and prove \rf{PRA-D}.
\par First, let us show that
\bel{PR-H2}
\|\Prm(a_1,H_2)-\Prm(a_2,H_2)\|\le 2\|a_1-a_2\|.
\ee
\par Indeed, thanks to \rf{L-NP},
$\ell_2\nparallel Ox_1$ and $\ell_2\nparallel Ox_2$
so that $\HR[H_2]=\RT$.  Hence, $a_1,a_2\in\HR[H_2]$.
\smsk
\par Therefore, thanks to \rf{FG-DS},
$$
\Prm(a_i,H_2)=a_i-\dist(a_i,H_2)\,\SN(\lh_2), ~~~~i=1,2.
$$
Hence,
\be
\|\Prm(a_1,H_2)-\Prm(a_2,H_2)\|
&=&
\|a_1-\dist(a_1,H_2)\,\SN(\lh_2)-
(a_2-\dist(a_2,H_2)\,\SN(\lh_2))\|
\nn\\
&\le&
\|a_1-a_2\|+|\dist(a_1,H_2)-\dist(a_2,H_2)|.
\nn
\ee
But $\dist(\cdot,H_2)$ is a Lipschitz function, and the proof of \rf{PR-H2} is complete.
\smsk
\par Now, from \rf{SD-A1}, \rf{PR-H2} and the triangle inequality, we have
\be
\|\Prm(a_1,H_1)-\Prm(a_2,H_2)\|
&\le&
\|(\Prm(a_1,H_1)-a_1)-(\Prm(a_1,H_2)-a_1)\|+
\|\Prm(a_1,H_2)-\Prm(a_2,H_2)\|
\nn\\
&\le&
\|\Prm(a_1,H_1)-a_1\|+\|\Prm(a_1,H_2)-a_1\|+2\|a_1-a_2\|
\nn\\
&=&
\dist(a_1,H_1)+\dist(a_1,H_2)+2\|a_1-a_2\|\le \delta+2\|a_1-a_2\|.
\nn
\ee
\par This proves the required inequality \rf{PRA-D} and the lemma.\bx
\begin{lemma} Inequality \rf{PRA-D} holds provided
$$
\text{either}~~~a_1\in S_2=H_2\cap (H_1+\delta Q_0)
~~~\text{or}~~~
a_2\in S_1=H_1\cap (H_2+\delta Q_0).
$$
\end{lemma}
\par {\it Proof.} For instance, suppose that
\bel{AA-P}
a_1\in H_2\cap (H_1+\delta Q_0).
\ee
We have
\be
\|\Prm(a_1,H_1)-\Prm(a_2,H_2)\|&\le& \|\Prm(a_1,H_1)-a_1\|+\|a_1-a_2\|+\|\Prm(a_2,H_2)-a_2\|
\nn\\
&=&
\dist(a_1,H_1)+\|a_1-a_2\|+\dist(a_2,H_2).
\nn
\ee
\par But $a_1\in H_1+\delta Q_0$ so that $\dist(a_1,H_1)\le\delta$. Furthermore, thanks to \rf{AA-P},
$$
\dist(a_2,H_2)\le \dist(a_2,H_2\cap (H_1+\delta Q_0))\le \|a_2-a_1\|
$$
because $a_1\in H_2\cap (H_1+\delta Q_0)$. Hence,
$$
\|\Prm(a_1,H_1)-\Prm(a_2,H_2)\|
\le\dist(a_1,H_1)+\|a_1-a_2\|+\dist(a_2,H_2)
\le \delta+2\|a_1-a_2\|
$$
proving \rf{PRA-D} and the lemma.\bx
\smsk
\begin{lemma} Inequality \rf{PRA-D} holds provided $\ell_1\parallel\ell_2$.
\end{lemma}
\par {\it Proof.} Because
$$
\ell_1\parallel\ell_2~~~~\text{and}~~~~\ell_i\nparallel Ox_j,~~~ i,j=1,2,
$$
(see \rf{L-NP}), we have
$$
\HR[H_1\cap(H_2+\delta Q_0)]\supset\HR[\ell_1]=\RT.
$$
Hence,
$$
\HR[H_1\cap(H_2+\delta Q_0)]=\RT.
$$
\par In the same way, we show that
$$
\HR[H_2\cap(H_1+\delta Q_0)]=\RT.
$$

\par Let us note that, since $\ell_1\parallel \ell_2$ and $\lh_i\perp \ell_i$ (see \rf{H-OL}), the vectors $\lh_1$ and $\lh_2$ are collinear unit vectors. Therefore,
$$
\text{either}~~~\lh_1=\lh_2~~~\text{or}~~~\lh_1=-\lh_2.
$$

\par If $\lh_1=\lh_2$, then
$$
\SN(\lh_1)=\SN(\lh_2)
$$
so that, thanks to Lemma \reff{L3-P}, inequality \rf{PRA-D} holds.
\par See Fig. 17.
\bsk\msk
\begin{figure}[H]
\hspace{19mm}
\includegraphics[scale=0.6]{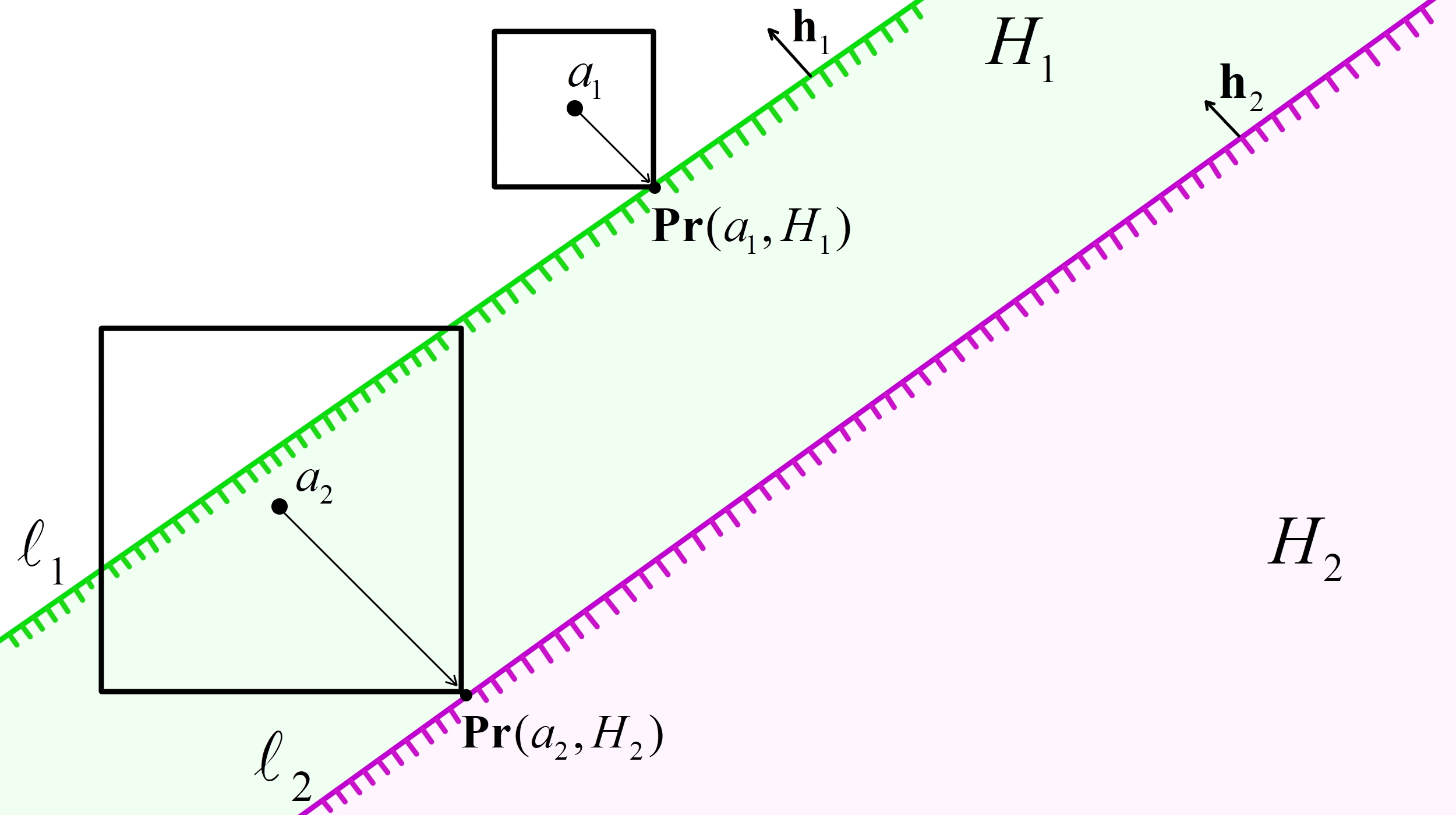}
\hspace*{-7mm}
\vspace*{3mm}
\caption{Metric projections onto $H_1$ and $H_2$ with the parallel boundaries: the first case.}
\end{figure}
\msk
\par Let us prove \rf{PRA-D} for $\lh_2=-\lh_1$. Consider three cases.
\msk
\par (a). Suppose that
$$
\text{either}~~~a_1\notin H_2~~~\text{or}~~~a_2\notin H_1.
$$
Prove that in this case \rf{PRA-D} holds.
\msk
\par For instance, let us assume that $a_1\notin H_2$.
See Fig. 18.

\begin{figure}[H]
\hspace{19mm}
\includegraphics[scale=0.6]{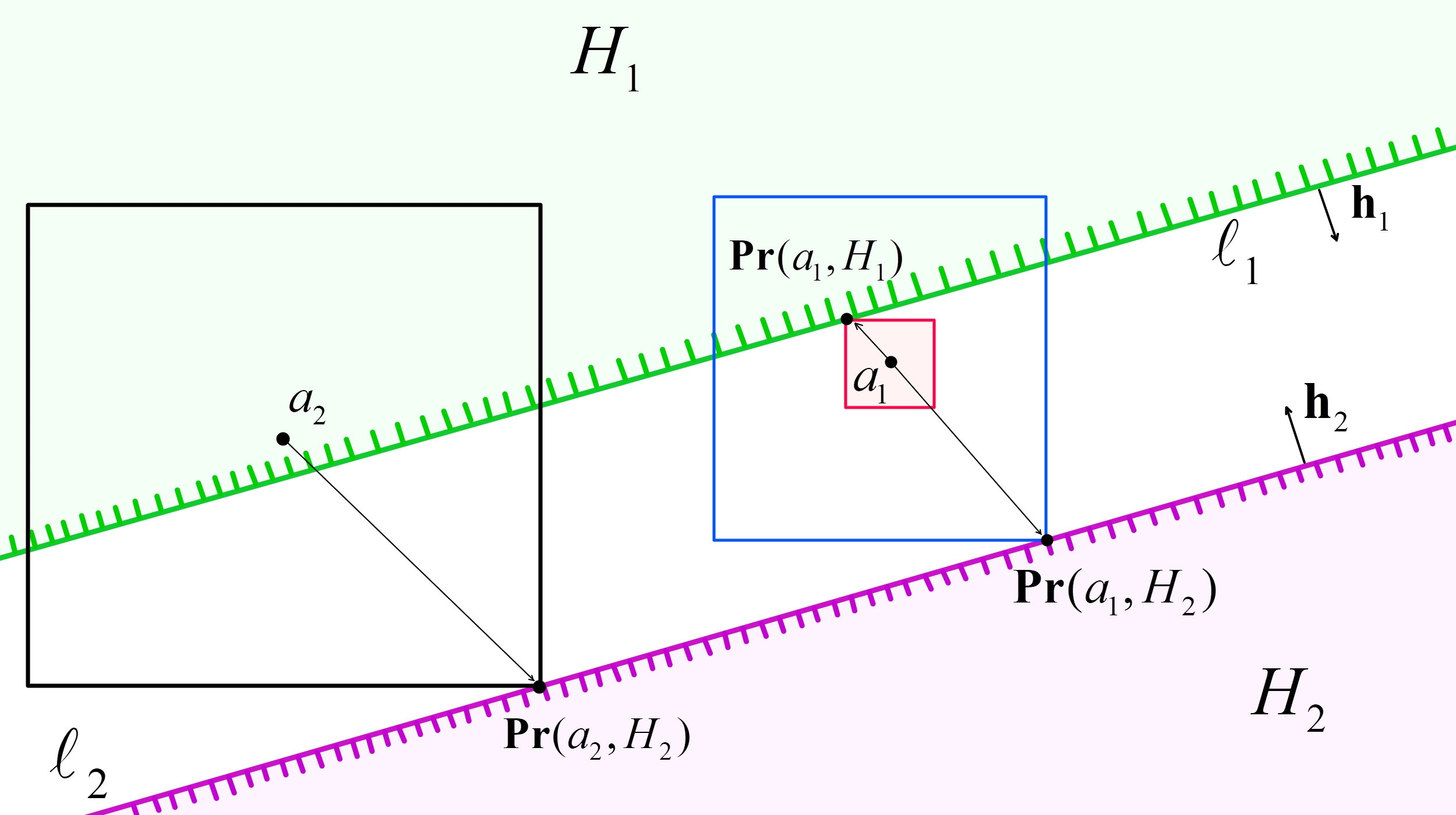}
\hspace*{-7mm}
\vspace*{3mm}
\caption{Metric projections onto $H_1$ and $H_2$ with the parallel boundaries: the second case.}
\end{figure}
\msk

\par In this case, $a_1\notin H_1\cup H_2$ (because  $a_1\notin H_1$, see \rf{A12-HN}). Clearly,
$H_1\cap H_2=\emp$.
\par Recall that, thanks to assumption \rf{H12-DF},
$$
H_1~~\text{and}~~H_2~~\text{are two half-planes with}~~ \dist(H_1,H_2)\le\delta.
$$
\par Because $H_1\cap H_2=\emp$ and
$\ell_1\parallel \ell_2$, we have
$$
\dist(\ell_1,\ell_2)\le\delta.
$$

\par Let us prove that in the case under consideration, we have
\bel{RW-4}
\dist(a_1,H_1)+\dist(a_1,H_2)\le\delta.
\ee
\par Indeed, let
$$
T=\clos(\RT\setminus (H_1\cup H_2)),
$$
i.e., $T$ is the strip between the half-planes $H_1$ and $H_2$. Thus, $\partial T=\ell\cup\ell_2$.
\par Because $\ell_1\parallel \ell_2$ and $\dist(\ell_1,\ell_2)\le \delta$, we have
$$
\dist(x,H_2)=\dist(x,\ell_2)\le\delta,~~~x\in\ell_1,
~~~~\text{and}~~~ \dist(x,H_1)=\dist(x,\ell_1)\le\delta,~~~x\in\ell_2.
$$

\smsk
\par We define a function $f$ on $T$ by letting
$$
f(x)=\dist(x,H_1)+\dist(x,H_2).
$$
\par Clearly, $f$ is a convex continuous function on $T$. Therefore, $\sup_T f=\sup_{\partial T} f$. But
$$
f(x)=\dist(x,\ell_2)\le\delta~~\text{on}~~\ell_1
~~~\text{and}~~~
f(x)=\dist(x,\ell_1)\le\delta~~\text{on}~~\ell_2
$$
so that $\sup_{\partial T} f\le\delta$. Hence, $\sup_{T} f\le\delta$ proving \rf{RW-4}.
\par Therefore, thanks to Lemma \reff{L-PH2}, inequality \rf{PRA-D} holds.

\msk
\par (b). Suppose that
$$
a_1\in H_2,~~a_2\in H_1~~~~\text{and}~~~~H_1\cap H_2=\emp.
$$
See Fig. 19.

\bsk\bsk
\begin{figure}[H]
\hspace{12mm}
\includegraphics[scale=0.6]{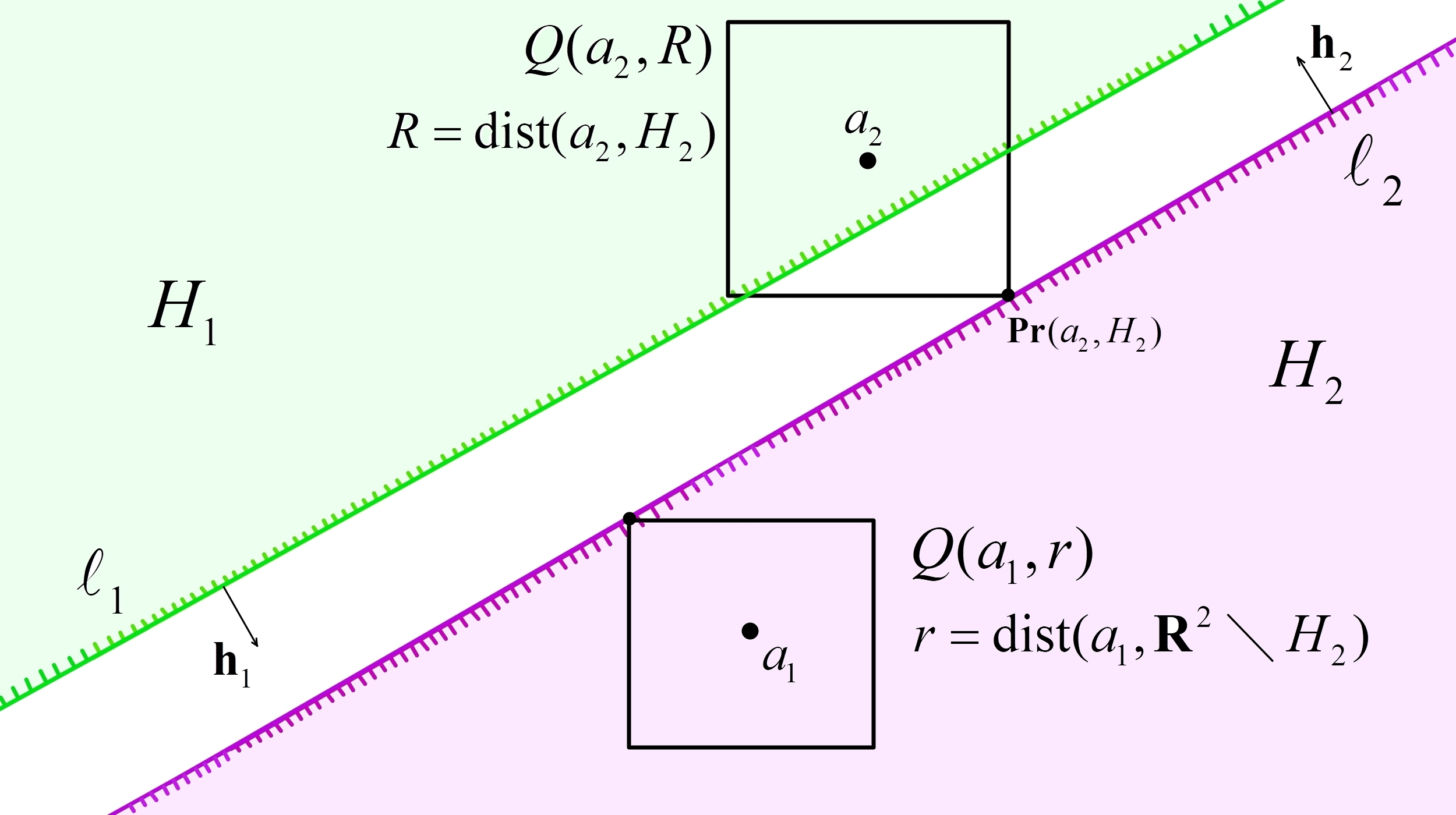}
\hspace*{-7mm}
\vspace*{3mm}
\caption{Half-planes $H_1$ and $H_2$ with the parallel boundaries: the third case.}
\end{figure}
\msk
\par Prove that in this case \rf{PRA-D} holds.

\par Let
$$
R=\dist(a_2,H_2)~~~\text{and}~~~
r=\dist(a_1,\RT\setminus H_2).
$$
Because $a_1\in H_2$ and $a_2\notin H_2$, the squares $Q(a,r)$ and $Q(a,R)$ are non-overlapping (i.e., they do not have common interior points). Hence,
$$
r+R\le\|a_1-a_2\|.
$$
\par Clearly, in the case under consideration,
$$
\RT\setminus H_2=H_1+p~~~~\text{with some}~~~p\in\RT~~~\text{with}~~~ \|p\|\le\delta.
$$
Therefore,
$$
\dist(a_1,H_1)\le \dist(a_1,\RT\setminus H_2)+\delta=r+\delta.
$$
Hence,
$$
\dist(a_1,H_1)+\dist(a_2,H_2)\le r+\delta+R
\le \delta+\|a_1-a_2\|.
$$
\par Finally,
\be
\|\Prm(a_1,H_1)-\Prm(a_2,H_2)\|&\le&
\|\Prm(a_1,H_1)-a_1\|+\|a_1-a_2\|
+\|\Prm(a_2,H_2)-a_2\|\nn\\
&=&
\dist(a_1,H_1)+\dist(a_2,H_2)+\|a_1-a_2\|\nn\\
&\le&
\delta+2\|a_1-a_2\|
\nn
\ee
proving the required inequality \rf{PRA-D}.

\msk
\par (c). Suppose that
$$
a_1\in H_2,~~a_2\in H_1~~~\text{and}~~~
H_1\cap H_2\ne\emp.
$$
See Fig. 20.
\bsk\bsk

\begin{figure}[H]
\hspace{19mm}
\includegraphics[scale=0.6]{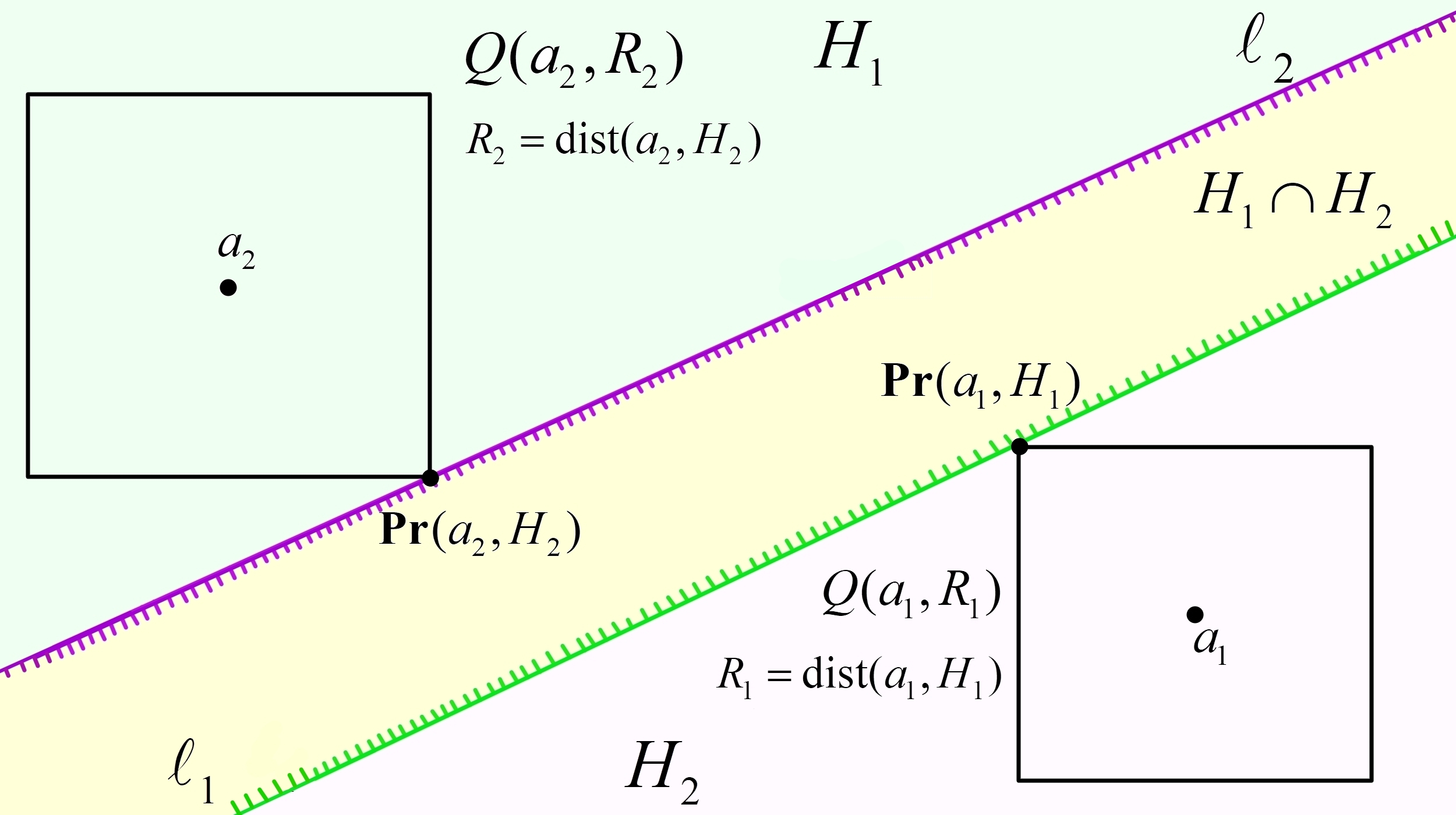}
\hspace*{-7mm}
\vspace*{3mm}
\caption{Half-planes $H_1$ and $H_2$ with the parallel boundaries: the fourth case.}
\end{figure}
\msk\msk
\par Prove that in this case \rf{PRA-D} holds. Indeed, let
$$
R_1=\dist(a_1,H_1)~~~~\text{and}~~~~
R_2=\dist(a_2,H_2).
$$
Then the squares $Q(a_1,R_1)$ and $Q(a_2,R_2)$ are non-overlapping so that
$$
R_1+R_2=\dist(a_1,H_1)+\dist(a_2,H_2)\le \|a_1-a_2\|.
$$
Hence,
\be
\|\Prm(a_1,H_1)-\Prm(a_2,H_2)\|&\le&
\|\Prm(a_1,H_1)-a_1\|+\|a_1-a_2\|
+\|\Prm(a_2,H_2)-a_2\|\nn\\
&=&
\dist(a_1,H_1)+\dist(a_2,H_2)+\|a_1-a_2\|\nn\\
&\le&
2\|a_1-a_2\|
\nn
\ee
proving \rf{PRA-D}.
\par The proof of the lemma is complete.\bx

\bsk
\par Summarizing the conditions imposed on the points $a_1$ and $a_2$ and the results obtained in this section, we conclude that to complete the proof of Proposition \reff{TWO-HP} it is sufficient to prove \rf{PRA-D} provided 
$$
\ell_1\nparallel \ell_2,~~~~\ell_i\nparallel Ox_j,~~~i,j=1,2,
$$
and $a_1$ and $a_2$ satisfy the following conditions:
\bel{A-HR1}
a_1\in \Hc[S_1]\setminus (H_1\cup S_2)
~~~\text{and}~~~
\Prm(a_1,H_1)\in S_1,
\ee
and
$$
a_2\in \Hc[S_2]\setminus (H_2\cup S_1)
~~~\text{and}~~~
\Prm(a_2,H_2)\in S_2.
$$
Recall that
$$
S_1=H_1\cap(H_2+\delta Q_0)~~~~\text{and}~~~~
S_2=H_2\cap(H_1+\delta Q_0).
$$
See \rf{S12-D}.
\smsk
\par Clearly, without loss of generality, we may assume that $\ell_1\cap\ell_2=\{0\}$. Then 
\bel{H12-DFN}
H_i=\{u\in\RT:\ip{\lh_i,u}\le 0\}, ~~~~i=1,2,
\ee
where $\lh_i$, $i=1,2$, are unit vectors. In these settings,
\bel{L12-DF}
\ell_i=\{u\in\RT: \ip{\lh_i,u}=0\}, ~~~i=1,2.
\ee
We know that
$$
\ell_1\nparallel \ell_2~~~~\text{and}
~~~~\ell_i\nparallel Ox_j,~~~i,j=1,2.
$$
We also know that
$$
\lh_i\perp \ell_i,~~ i=1,2.
$$
\par Therefore,
$$
\lh_1\nparallel \lh_2,~~~~\text{and each of the vectors}~~~~
\lh_i,~i=1,2,~~~~\text{has non-zero coordinates.}
$$

\par  Note that, in the case under consideration, the set
$$
S=H_1\cap H_2
$$
is a convex cone with the vertex at $0$. Moreover, the sets
$$
S_1=H_1\cap(H_2+\delta Q_0)~~~~\text{and}~~~~
S_2=H_2\cap(H_1+\delta Q_0)
$$
are convex cones in $\RT$. Let
\bel{VS-12}
X_1=(s_1,s_2)~~~~\text{and}~~~~X_2=(t_1,t_2)
\ee
be the vertices of the cones $S_1$ and $S_2$ respectively.
Thus,
\bel{X1-L}
X_1=\ell_1\cap \tell_2~~~\text{where}~~~\tell_2=\partial(H_2+\delta Q_0),
\ee
and
\bel{X2-L}
X_2=\ell_2\cap \tell_1~~~\text{where}~~~\tell_1=\partial(H_1+\delta Q_0).
\ee
\par Moreover, we have the following representations of the cones $S_1$ and $S_2$:
\bel{S-X12}
S_1=S+X_1=H_1\cap H_2+X_1,~~~ S_2=S+X_2=H_1\cap H_2+X_2.
\ee
\msk
\par See Fig. 21.

\begin{figure}[H]
\hspace*{5mm}
\includegraphics[scale=0.7]{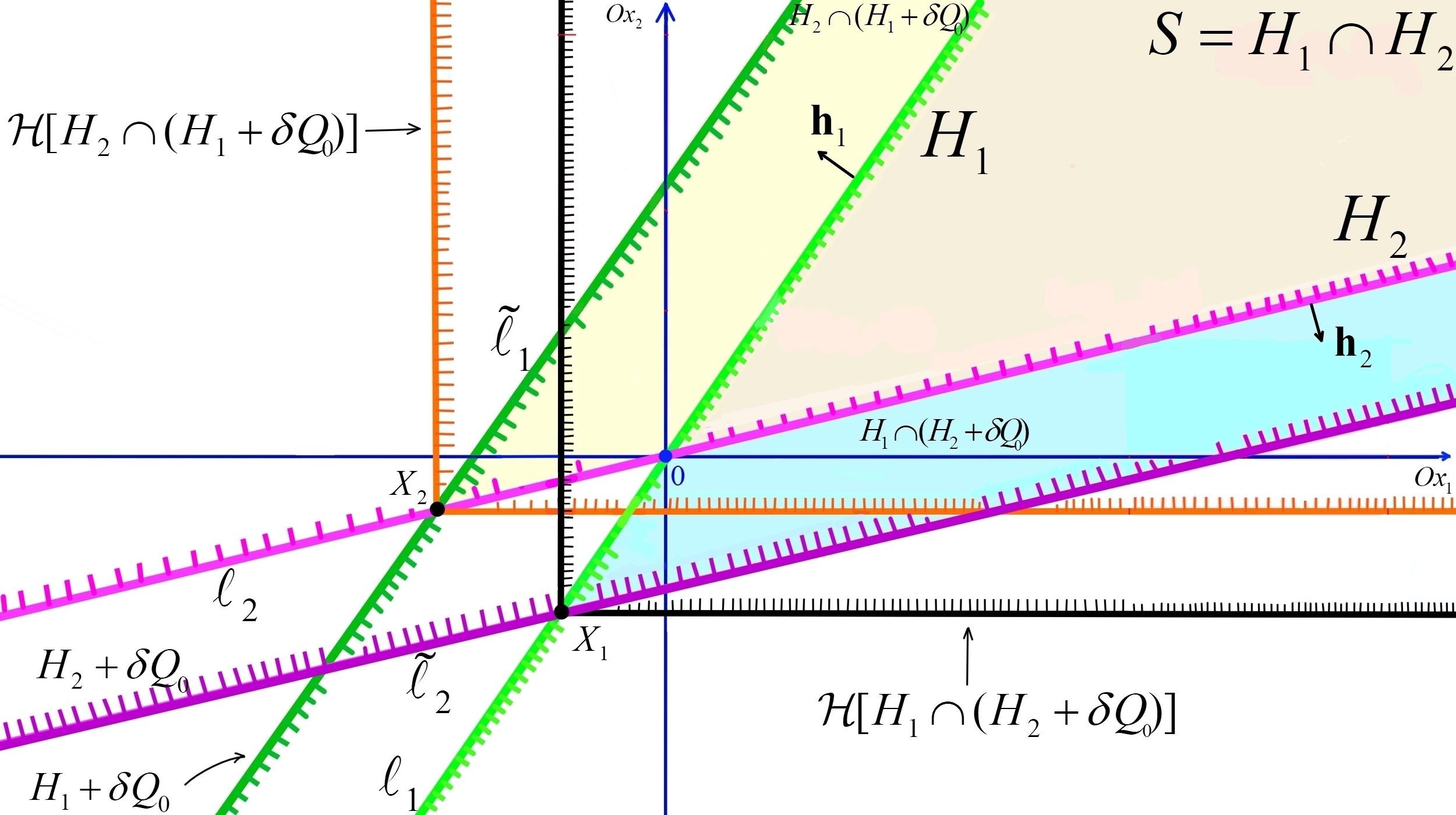}
\hspace*{12mm}
\caption{The half-planes $H_1$ and $H_2$ and points $X_1$ and $X_2$.}
\end{figure}
\bsk

\par Let us represent the vectors $\lh_1$ and $\lh_2$ in the polar form:
\bel{H-12}
\lh_i=(\cos\vf_i,\sin\vf_i)~~~\text{where the angle}~~~ \vf_i\in[0,2\pi),~~~i=1,2.
\ee

\par Let us see that the following
explicit formulae for the points $X_1$ and $X_2$ hold:
\bel{X-1}
X_1=\frac{\delta}{\Delta}\,\|\lh_2\|_1\,(-\sin \vf_1,\cos \vf_1).
\ee
and
\bel{X-2}
X_2=\frac{\delta}{\Delta}\,\|\lh_1\|_1\,
(\sin \vf_2,-\cos \vf_2).
\ee
where
$$
\Delta=\sin(\vf_2-\vf_1).
$$

\par Indeed, the half plane $H_i+\delta Q_0$, $i=1,2$, has the following representation:
$$
H_i+\delta Q_0=\{u\in\RT:\ip{\lh_i,u}\le \delta\,\|\lh_i\|_1\}, ~~~i=1,2.
$$
(Recall that, given $u=(u_1,u_2)\in\RT$ we let  $\|u\|_1=|u_1|+|u_2|$ denote the $\ell^1_2$-norm in $\RT$.)
Hence,
\bel{E-LW}
\tell_i=\{u\in\RT: \ip{\lh_i,u}=\delta\, \|\lh_i\|_1\}, ~~~i=1,2.
\ee
\par Let
\bel{A-DL}
A=\left(
\begin{array}{ll}
\cos \vf_1& \sin \vf_1\\
\cos \vf_2& \sin \vf_2
\end{array}
\right)
~~~\text{and let}~~~\Delta=\det A=\sin(\vf_2-\vf_1).
\ee
(Clearly, $\Delta\ne 0$ because $\lh_1\nparallel\lh_2$.) We know that $$X_1=(s_1,s_2)=\ell_1\cap \tell_2$$ so that, thanks to \rf{L12-DF} and \rf{E-LW}, the vector $(s_1,s_2)$ is the solution of the system of linear equations
$$
A\left(
\begin{array}{l}
s_1\\s_2
\end{array}
\right)=\left(
\begin{array}{l}
\hspace*{5mm} 0\\
\delta\,\|\lh_2\|_1
\end{array}
\right).
$$
Therefore,
$$
s_1=\frac{1}{\Delta}\,\left|
\begin{array}{ll}
0& \sin \vf_1\\
\delta\,\|\lh_2\|_1& \sin \vf_2
\end{array}
\right|
=-\frac{\delta}{\Delta}\,\|\lh_2\|_1\,\sin \vf_1
~~~\text{and}~~~
s_2=\frac{\delta}{\Delta}\,\|\lh_2\|_1\,\cos \vf_1
$$
proving \rf{X-1}.
\par In the same way we prove \rf{X-2}.

\bsk

\par Let $Ox^+_i$ and $Ox^-_i$ be the non-negative and non-positive semi-axes of the axis $Ox_i$ respectively, $i=1,2$. Thus,
$$
Ox^+_i=\{t e_i:t\ge 0\}~~~~\text{and}~~~~
Ox^-_i=\{t e_i:t\le 0\}, ~~~~i=1,2,
$$
where $e_1=(1,0)$ and $e_2=(0,1)$.
\msk


\par Let us prove that inequality \rf{PRA-D} holds in each of the following cases:
\msk

\par {\bf Case 1.} The cone $S=H_1\cap H_2$ contains two distinct coordinate semi-axes.
\smsk
\par {\bf Case 2.} The cone $S$ contains exactly one coordinate semi-axis.
\smsk
\par {\bf Case 3.} The cone $S$ does not contain any coordinate semi-axis.
\msk
\par We begin the proof with the first case.

\begin{lemma} Inequality \rf{PRA-D} holds provided the cone $S=H_1\cap H_2$ contains two distinct semi-axes.
\end{lemma}
\par {\it Proof.} Without loss of generality we may assume that $S$ contains the semi-axes $Ox^+_1=\{t e_1:t\ge 0\}$ and $Ox^-_2=\{t e_2:t\ge 0\}$. See Fig. 22.

\bsk
\begin{figure}[H]
\hspace{14mm}
\includegraphics[scale=0.22]{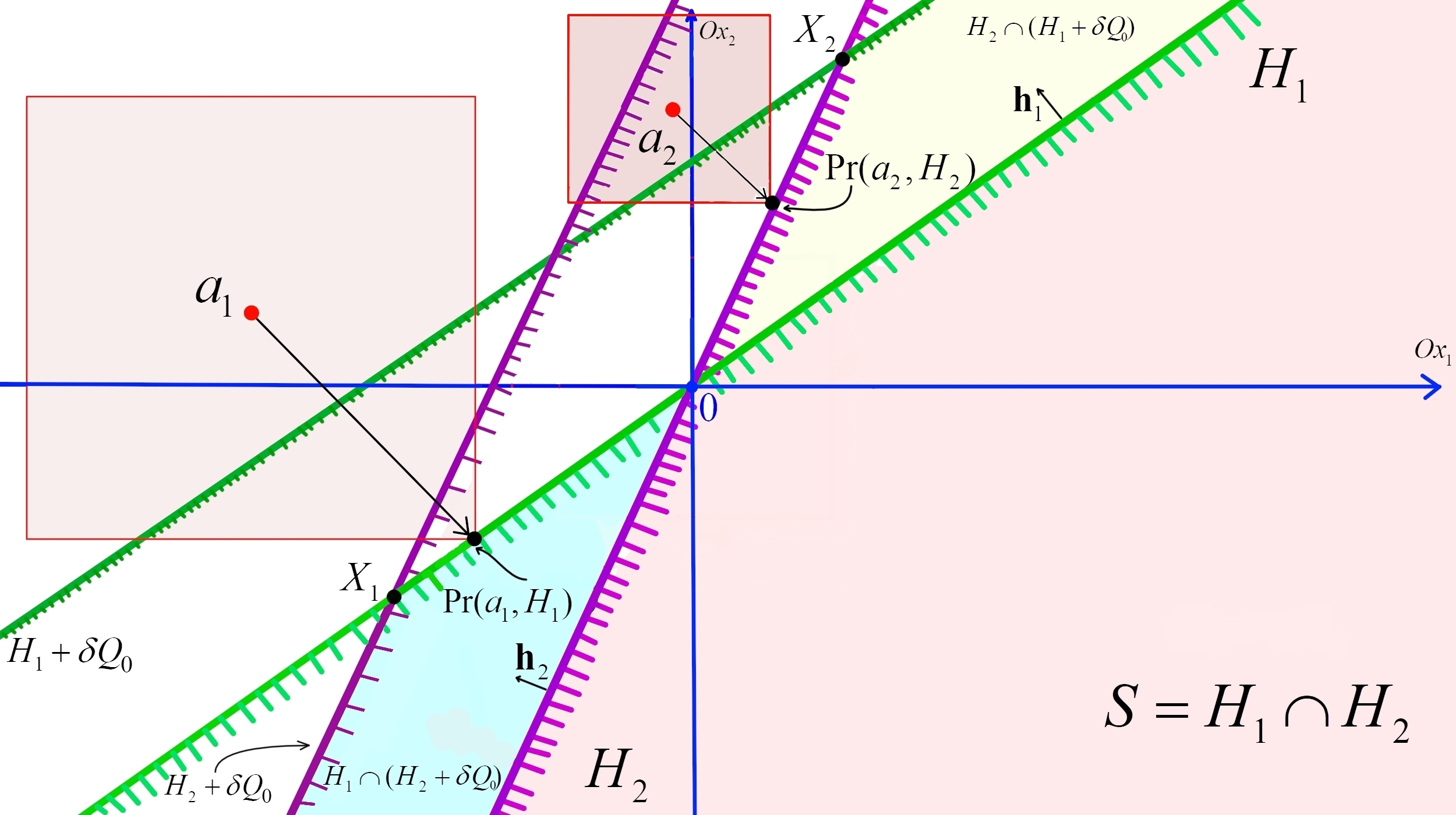}
\hspace*{-10mm}
\caption{The half-planes $H_1$ and $H_2$ with the non-parallel boundaries: the first case.}
\end{figure}

\par Thus, $S\supset Ox^+_1\cup Ox^-_2$ so that $S$ contains the forth coordinate quadrant
$$
Q_4=\{(s,t)\in\RT: s\ge 0, t\le 0\}.
$$
\par Because $H_i\supset S=H_1\cap H_2$, $i=1,2$, the half-plane $H_i\supset Q_4$ as well.
\par Therefore, the vector $\lh_i$ belongs to the second coordinate quadrant $Q_2=\{(s,t)\in\RT: s\le 0, t\ge 0\}$ proving that $\SN(\lh_1)=\SN(\lh_2)=(-1,1)$. Lemma \reff{L3-P} tells us that in this case inequality \rf{PRA-D} holds, and the proof of the lemma is complete.\bx
\bsk
\par We turn to the proofs of inequality \rf{PRA-D} in {\bf Case 2} and {\bf Case 3}. Below we present several auxiliary results which we use in these proofs.
\begin{lemma}\lbl{TR-X12} Inequality \rf{PRA-D} holds provided either $a_1$ or $a_2$ belongs to the triangle
$$
T=\Delta(X_1,X_2,O)
$$
with the vertices at points $X_1$, $X_2$ and $O$.
\end{lemma}
\par {\it Proof.} We define a function $D:T\to\R$ by letting
$$
D(x)=\dist(x,H_1)+\dist(x,H_2),~~~x\in T.
$$
\par Prove that $D\le\delta$ on $T$. Indeed, $D$ is a convex continuous function so that its maximum on $T$ is attained on the set of vertices of the triangle $T$, i.e., at the points $O,X_1$ and $X_2$. But $D(O)=0$, $D(X_1)=\delta$ because
$$
X_1\in\tell_2=\partial(H_2+\delta\,Q_0),
$$
and $D(X_2)=\delta$ because
$$
X_2\in\tell_1=\partial(H_1+\delta\,Q_0).
$$
\par Therefore, if $a_i\in T$ for some $i=1,2$, then
$$
D(a_i)=\dist(a_i,H_1)+\dist(a_i,H_2)\le\delta.
$$
\par Lemma \reff{L-PH2} tells us that in this case \rf{PRA-D} holds, and the proof of the lemma is complete.\bx
\begin{lemma}\lbl{TG-A} Let $k_1,k_2\in\R$, $0<|k_2|\le k_1$, and let $G_1,G_2\subset\RT$ be two half-planes determined as follows:
$$
G_1=\{(x_1,x_2)\in\RT:-k_1x_1+x_2\le 0\}~~~~
\text{and}~~~~
G_2=\{(x_1,x_2)\in\RT: k_2x_1-x_2\le 0\}.
$$
\par Let $a=(a_1,a_2)\in\RT$, $a_1\ge 0$, $a\notin G_1$, and let a point $b\notin G_2$. Then
$$
\dist(a,G_1)+\dist(b,G_2)\le\|a-b\|.
$$
\end{lemma}
\par {\it Proof.} Let
$$
d_1=\dist(a,G_1), ~~~r=\dist(a,\RT\setminus G_2), ~~~~\text{and}~~~d_2=\dist(b,G_2).
$$
See Fig. 23

\begin{figure}[H]
\hspace{8mm}
\includegraphics[scale=0.7]{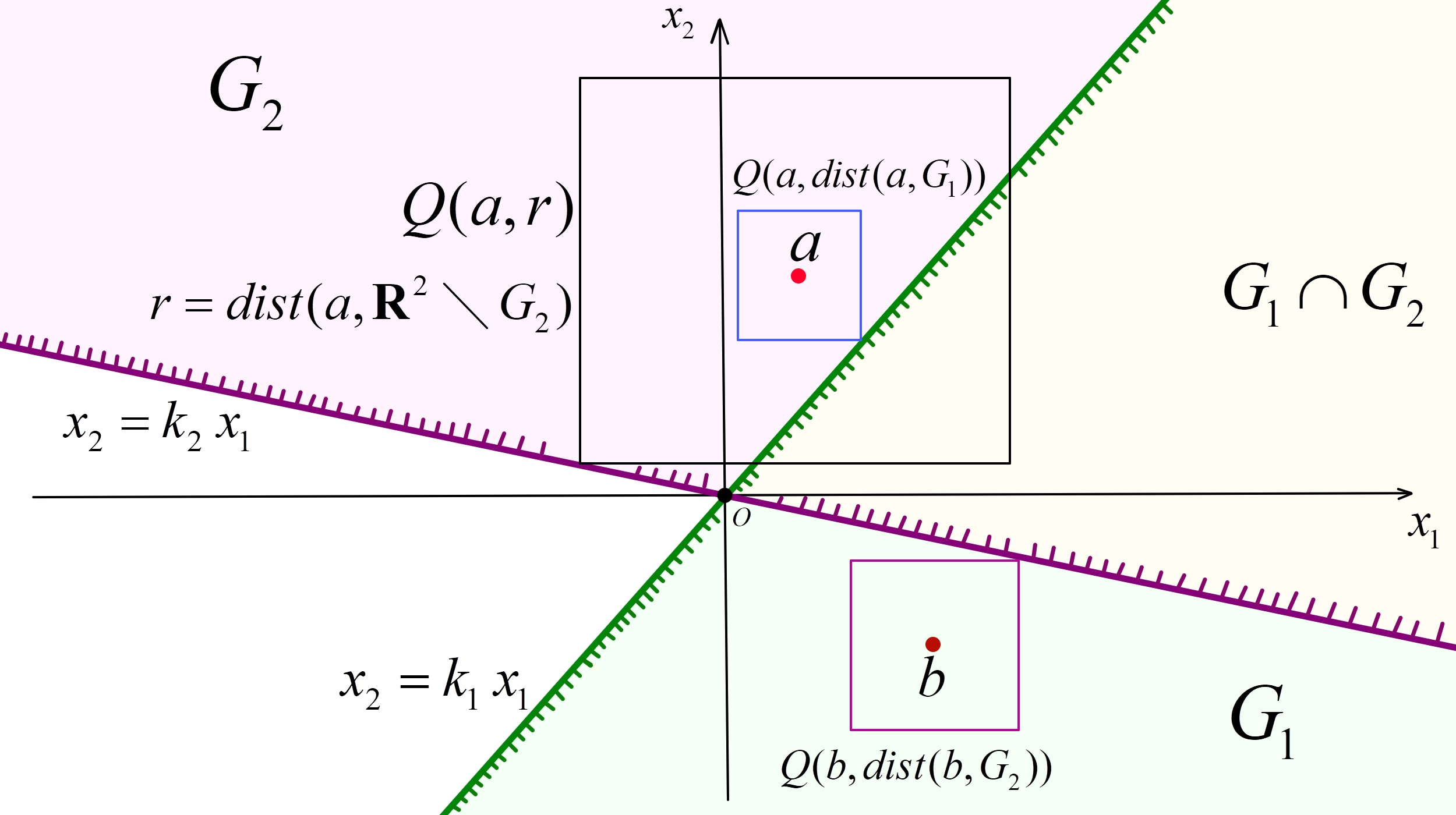}
\hspace*{10mm}
\caption{Lemma \reff{TG-A}.}
\end{figure}
\msk
\par We know that
$$
G_1=\{x\in\RT:\ip{g_1,x}\le 0\} ~~~\text{with}~~~g_1=(-k_1,1),
$$
and
$$
\RT\setminus G_2=\{x\in\RT:\ip{g_2,x}< 0\} ~~~\text{with}~~~g_2=(-k_2,1).
$$
Hence,
$$
d_1=\dist(a,G_1)=[\ip{g_1,a}]_+/\|g_1\|_1=
[a_2-k_1a_1]_+/(|k_1|+1)
$$
and
$$
r=\dist(a,\RT\setminus G_2)=[\ip{g_2,a}]_+/\|g_2\|_1=
[a_2-k_2a_1]_+/(|k_2|+1).
$$
\par Recall that
$$
a=(a_1,a_2)\notin G_1,~~a_1\ge 0,
~~~~\text{and}~~~~ |k_2|\le k_1.
$$
Hence, $a_2-k_1a_1\ge 0$ so that
$$
d_1=\frac{a_2-k_1a_1}{|k_1|+1}
\le \frac{a_2-k_1a_1}{|k_2|+1}
\le \frac{a_2-k_2a_1}{|k_2|+1}=r.
$$
\par From this we have
$$
Q(a,d_1)\subset Q(a,r)=Q(a,\dist(a,\RT\setminus G_2))
\subset G_2.
$$
\par Also, let us note that $b\notin G_2$ so that
$$
Q(b,d_2)=Q(b,\dist(b,G_2))\subset \clos(\R^2\setminus G_2).
$$
\par Thus, the straight line $x_2=k_2x_1$ separates (not strictly) the squares $Q(a,d_1)$ and $Q(b,d_2)$. Therefore, these squares are non-overlapping so that
$$
d_1+d_2\le \|a-b\|.
$$
\par The proof of the lemma is complete.\bx

\begin{lemma}\lbl{FL-P} Let $\alpha_i\in(0,\pi/2)$, $i=1,2$, be the angle between the positive direction of the axis $Ox_1$ and the straight line $\ell_i=\partial H_i$, and let $\alpha_2\le\alpha_1$.
\par Suppose that either the cone
$$
S=H_1\cap H_2\subset \RT_+=Ox_1^+\times Ox_2^+
$$
or
$$
S~~\text{contains a single coordinate semi-axis, say}~~ Ox_1^+.
$$

\par Under these conditions, if
\bel{AM-1}
a_1\in\Hc[S_2]=\Hc[H_2\cap(H_1+\delta Q_0)], ~~~~
a_1\notin H_1+\delta Q_0,
\ee
then inequality \rf{PRA-D} holds.
\end{lemma}
\par {\it Proof.} Let
$$
K=\Hc[S_2]\setminus (H_1+\delta Q_0)=\Hc[H_2\cap(H_1+\delta Q_0)]\setminus (H_1+\delta Q_0).
$$
Then, thanks to \rf{AM-1},
\bel{A1-BK}
a_1\in K.
\ee
\par On Fig. 24 we show the set $K$ and the relative positions of the half-spaces $H_1$ and $H_2$ provided $S$ contains the (single) semi-axis $Ox_1^+$. (This is
{\bf Case 2} of our approach.)

\bsk\bsk\bsk
\begin{figure}[H]
\hspace{8mm}
\includegraphics[scale=0.7]{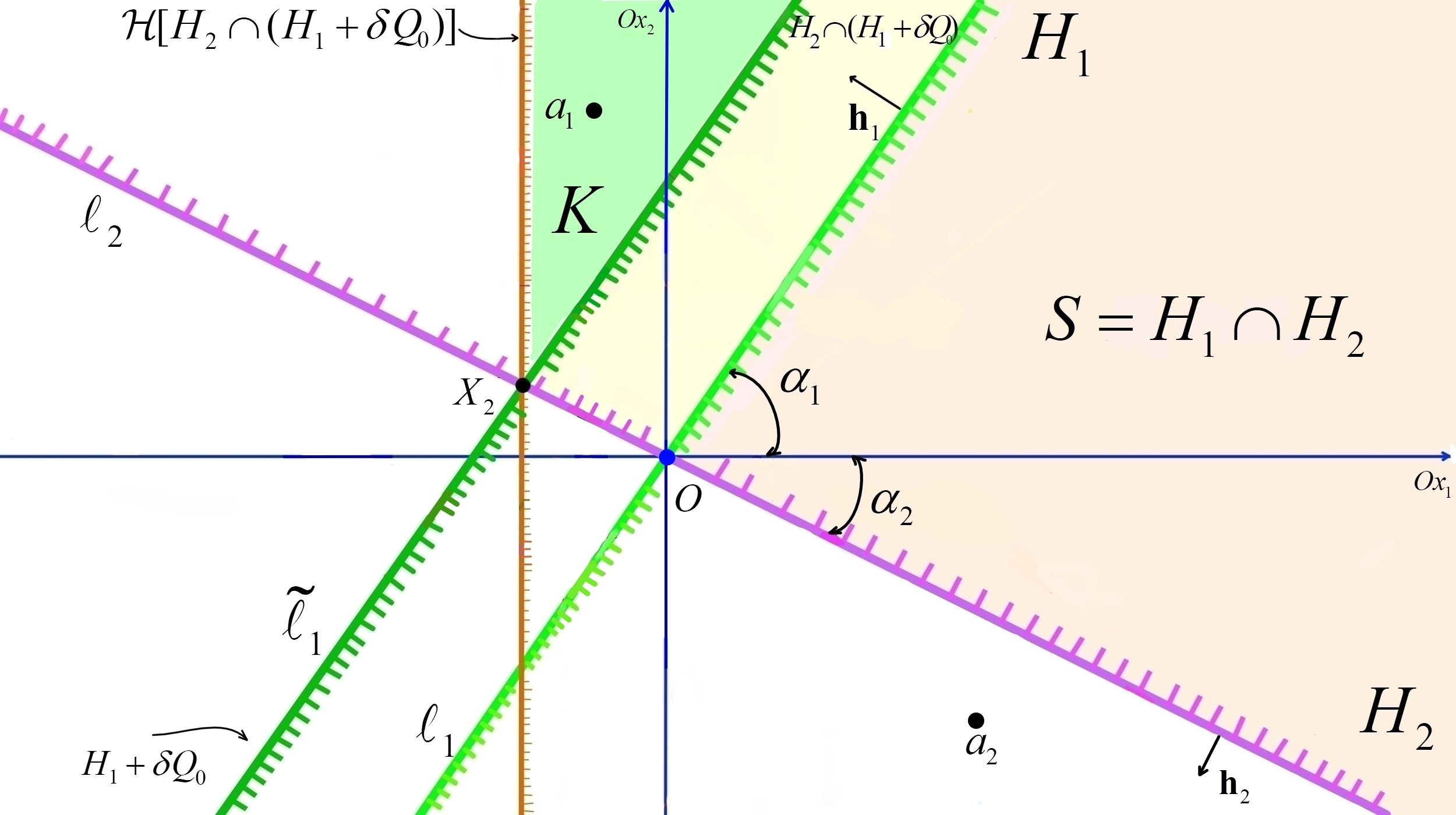}
\hspace*{10mm}
\caption{Lemma \reff{FL-P}.}
\end{figure}


\par On Fig. 25 we show the set $K$ and the relative positions of the half-spaces $H_1$ and $H_2$ provided
$S=H_1\cap H_2\subset \RT_+$, i.e., $S$ does not contain any coordinate axis ({\bf Case 3} of our approach.)
\bsk\msk

\begin{figure}[H]
\hspace{8mm}
\includegraphics[scale=0.7]{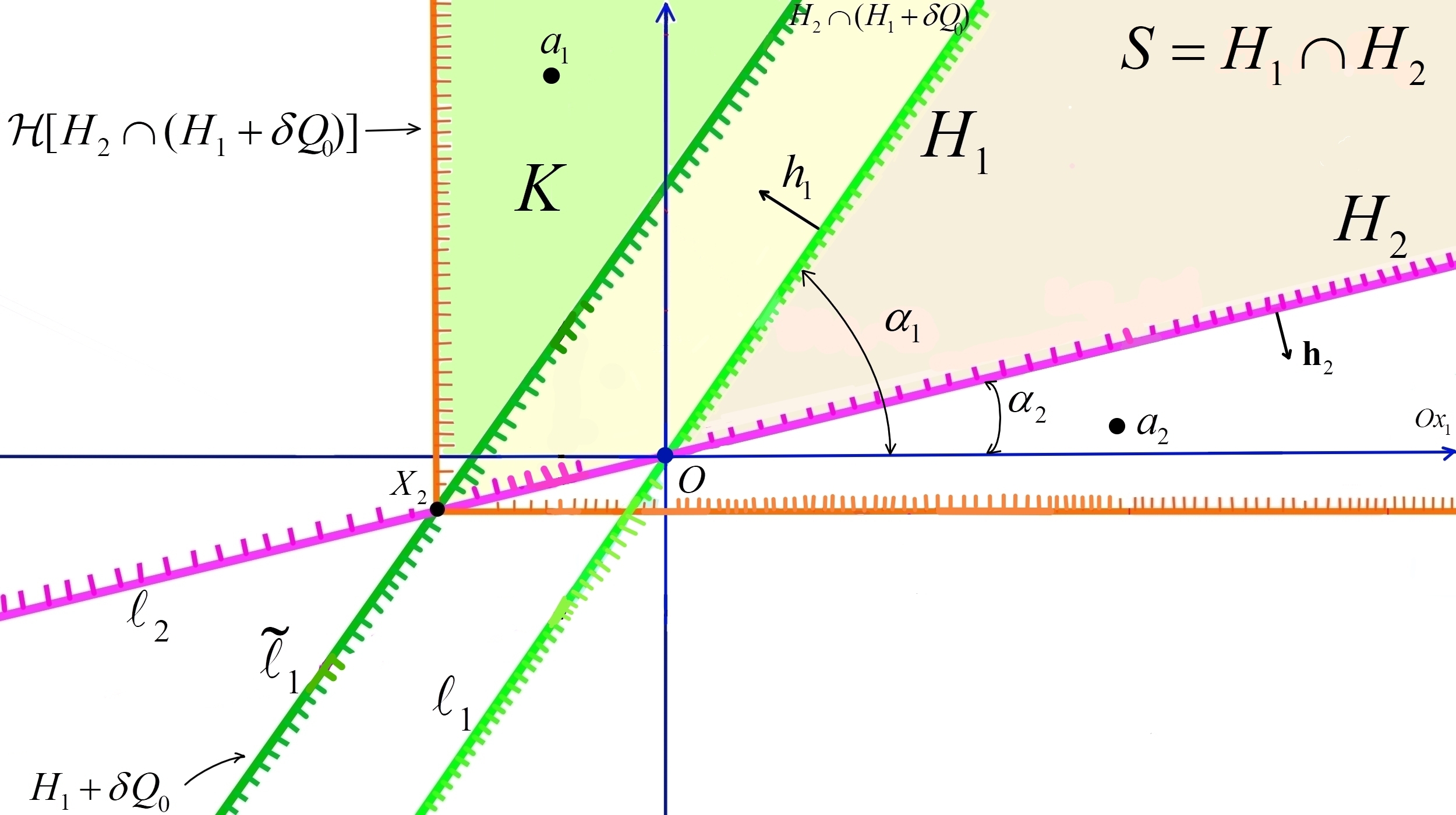}
\hspace*{10mm}
\caption{Lemma \reff{FL-P}.}
\end{figure}
\msk
\bsk

\par On figures Fig. 24 and 25, $X_2$ is the point defined by \rf{X2-L}, i.e.,
$$
X_2=\ell_2\cap \tell_1 ~~~\text{where}~~~\tell_1=\partial(H_1+\delta Q_0).
$$
Also recall that $\ell_2=\partial H_2$.
\par Let $k_1=\tan(\alpha_1)$ and $k_2=\tan(\alpha_2)$. Then $k_1>0$, $|k_2|\le k_1$. Furthermore,
$$
G_1=\{(x_1,x_2)\in\RT:-k_1x_1+x_2\le 0\}=H_1~~~~~
$$
and
$$
G_2=\{(x_1,x_2)\in\RT:k_2x_1-x_2\le 0\}=H_2.
$$
\par Let
$$
\tK=\{(x_1,x_2)\in\RT: x_1\ge 0,~(x_1,x_2)\notin H_1\}.
$$
\par Then $\tK=K-X_2$. It is also clear that
\bel{H1-SH}
H_1+\delta Q_0=H_1+X_2
~~~\text{(because}~~~
X_2\in\partial(H_1+\delta Q_0)).
\ee
\par Let
$$
a=(a^{(1)},a^{(2)})=a_1-X_2.
$$
\par Since $a_1\in K$, see \rf{A1-BK}, and $\tK=K-X_2$, the point $a\in\tK$.
\smsk
\par Let $b=a_2-X_2$. Because $X_2\in\partial H_2$
and $a_2\notin H_2$, the point $b_2\notin H_2=G_2$.
\par Lemma \reff{TG-A} tells us that in this case
$$
\dist(a,G_1)+\dist(b,G_2)\le\|a-b\|=\|a_1-a_2\|.
$$
Furthermore, thanks to \rf{H1-SH},
$$
\dist(a,G_1)=\dist(a_1-X_2,H_1)=\dist(a_1,H_1+X_2)
=\dist(a_1,H_1+\delta Q_0)
$$
and
$$
\dist(b,G_2)=\dist(a_2-X_2,H_2)=\dist(a_2,H_2+X_2)
=\dist(a_2,H_2).
$$
\par Then, thanks to Lemma \reff{TG-A},
\bel{G-4}
\dist(a_1,H_1+\delta Q_0)+\dist(a_2,H_2)\le\|a_1-a_2\|.
\ee

\par A geometrical proof of this inequality can be given with the help of Fig. 26 below.
\bsk

\begin{figure}[H]
\hspace{8mm}
\includegraphics[scale=0.7]{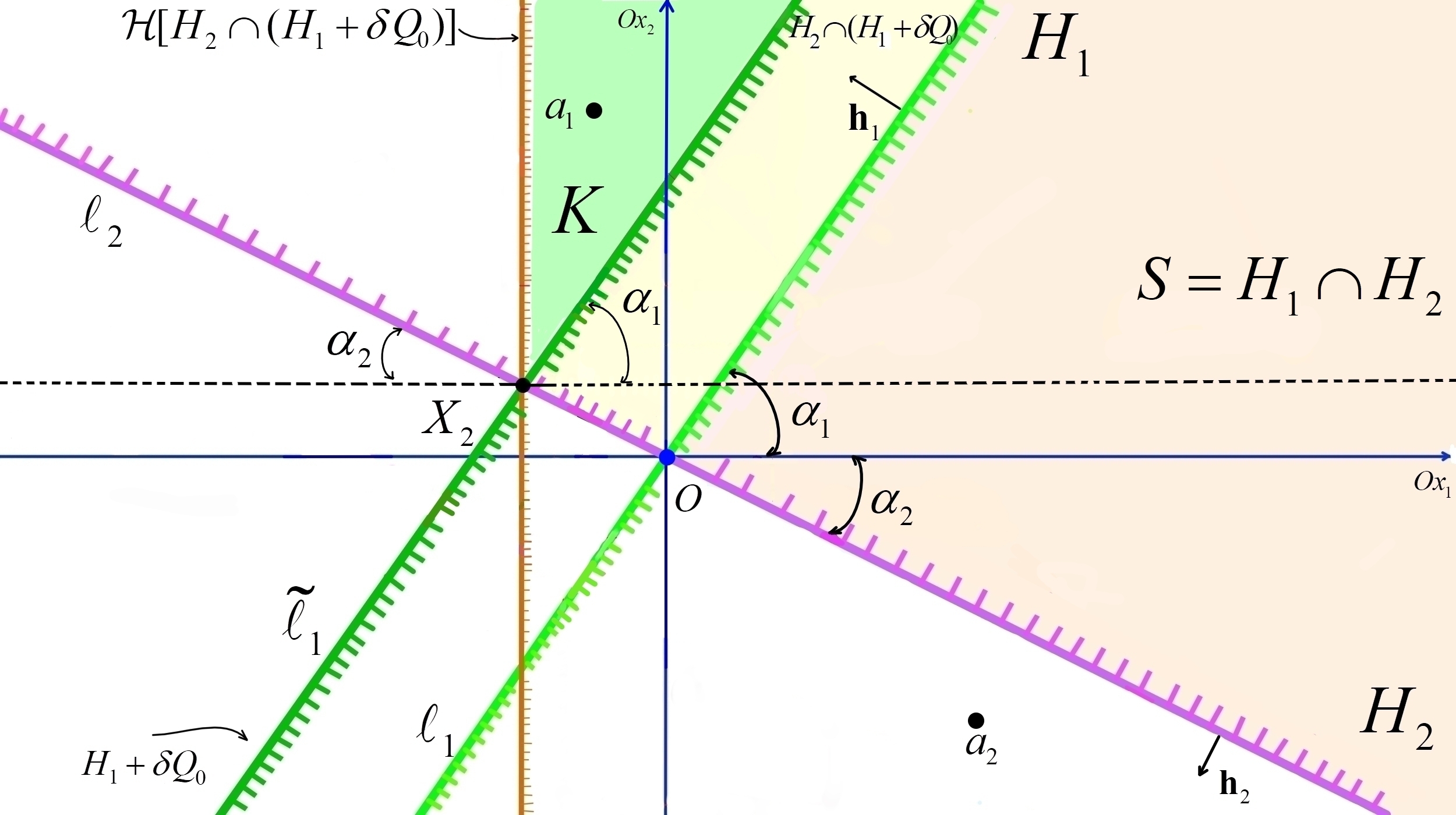}
\hspace*{10mm}
\caption{Lemma \reff{FL-P}.}
\end{figure}
\msk
\par We introduce a new coordinate system in $\RT$ by  shifting of the axis $Ox_1$ to the dotted line shown on Fig. 26. In other words, we shift the origin to the point $X_2$.  In this new coordinate system we set $G_1=H_1+\delta Q_0$ and $G_2=H_2$. Then we apply to these settings Lemma \reff{TG-A}, and get the required inequality \rf{G-4}.
\msk
\par Now, we complete the proof of Lemma \reff{FL-P} as follows. We note that
$$
\dist(a_1,H_1)\le\dist(a_1,H_1+\delta Q_0)+\delta.
$$
\par From this and \rf{G-4}, we have
$$
\dist(a_1,H_1)+\dist(a_2,H_2)\le
\dist(a_1,H_1+\delta Q_0)+\dist(a_2,H_2)+\delta
\le \delta+\|a_1-a_2\|.
$$
Hence,
\be
\|\Prm(a_1,H_1)-\Prm(a_2,H_2)\|&\le&
\|\Prm(a_1,H_1)-a_1\|+\|a_1-a_2\|
+\|\Prm(a_2,H_2)-a_2\|\nn\\
&=&
\dist(a_1,H_1)+\dist(a_2,H_2)+\|a_1-a_2\|\nn\\
&\le&
\delta+2\|a_1-a_2\|
\nn
\ee
proving the lemma.\bx
\msk
\par The next lemma states that \rf{PRA-D} holds in {\bf Case 2}.
\begin{lemma}\lbl{L6-P} Inequality \rf{PRA-D} holds provided the cone $S=H_1\cap H_2$ contains exactly one semi-axis.
\end{lemma}
\par {\it Proof.} Without loss of generality, we may assume that $S$ contains the semi-axis $Ox_1^+$.
\par Let $\alpha_i\in(0,\pi/2)$, $i=1,2$, be the angle between the straight line $\ell_i$ and the axis $Ox_1$.
\par Because the uniform norm on the plane is invariant under reflections with respect to the coordinate axes, without loss of generality we can also assume that $\alpha_2\le\alpha_1$.

\par See Fig. 27.

\bsk\msk
\begin{figure}[H]
\hspace{8mm}
\includegraphics[scale=0.7]{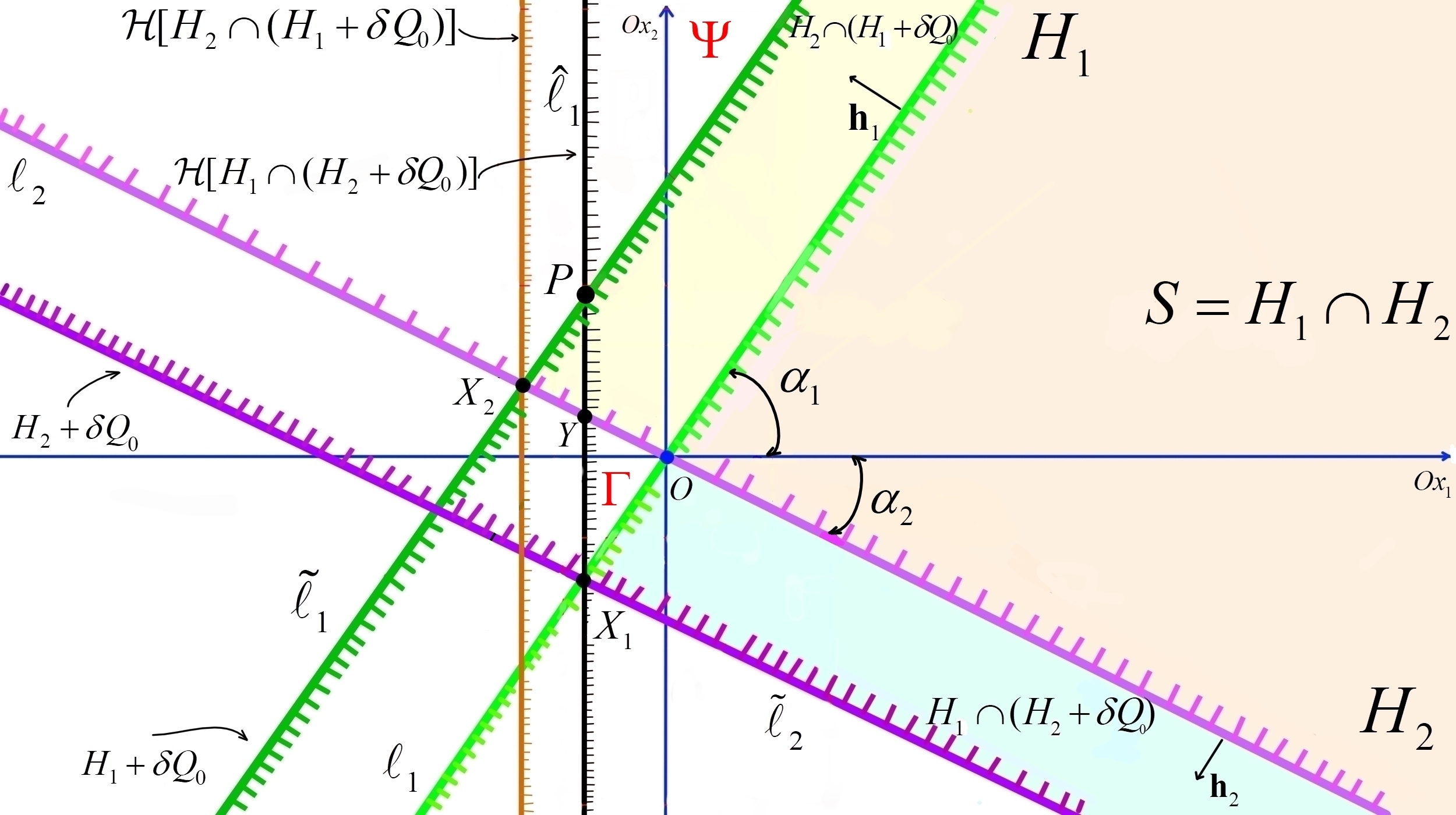}
\hspace*{10mm}
\caption{The half-planes $H_1$ and $H_2$ with the non-parallel boundaries: {\bf Case 2}.}
\end{figure}
\msk

\par We know that
\bel{A-HR1N}
a_1\in \Omega=\Hc[S_1]\setminus (H_1\cup S_2)
\ee
where, $S_1=H_1\cap(H_2+\delta Q_0)$ and
$S_2=H_2\cap(H_1+\delta Q_0)$.
See \rf{A-HR1}.
\par We also recall that
$$
\lh_1=(\cos\vf_1,\sin\vf_1)~~~\text{and}~~~
\lh_2=(\cos\vf_2,\sin\vf_2)
$$
are the vectors determining the half-spaces $H_1$ and $H_2$
respectively. See \rf{H12-DFN} and \rf{H-12}. (I.e., $\lh_i$ is a vector orthogonal to the straight line $\ell_i$, the boundary of $H_i$, see \rf{L12-DF}, and directed outside of $H_i$.) In the case under consideration, we have
\bel{FI-C1}
\vf_1\in (\pi/2,\pi)~~~~~\text{and}~~~~\vf_2\in (\pi,3/2\pi).
\ee
\par Thus,
\bel{AL-12}
\vf_1=\alpha_1+\pi/2~~~~~\text{and}~~~~~ \vf_2=\tfrac32\pi-\alpha_2.
\ee
\par Furthermore, because the convex cone $S=H_1\cap H_2$ contains the positive semi-axis $Ox_1^+$, the rectangular hulls of the sets $S_1=H_1\cap (H_2+\delta Q_0)$ and  $S_2=H_2\cap (H_1+\delta Q_0)$ have the following representations:
$$
\HR[S_1]=\{(u_1,u_2)\in\RT:u_1\ge s_1\}~~~\text{and}~~~
\HR[S_2]=\{(u_1,u_2)\in\RT:u_1\ge t_1\}.
$$
\par We recall that
$$
X_1=\ell_1\cap \tell_2=(s_1,s_2)~~~~\text{and}~~~~
X_2=\ell_2\cap \tell_1=(t_1,t_2)
$$
are the vertices of the cones $S_1$ and $S_2$ respectively. (Here $\ell_i=\partial H_i$ and
$\tell_i=\partial(H_i+\delta Q_0)$, $i=1,2$.) See \rf{VS-12}, \rf{X1-L} and \rf{X2-L}.
\smsk
\par Let
\bel{LH-1}
\hat{\ell}_1=\{(s_1,u)\in\RT:u\in\R\}
\ee
be the boundary of the rectangle $\HR[S_1]$. Also, we let
$$
P=\hat{\ell}_1\cap\tell_1
$$
denote the point of intersection of the straight lines $\hat{\ell}_1$ and $\tell_1=\partial(H_1+\delta Q_0)$. See \rf{X2-L} and Fig. 27.
\msk

\par Let us prove that
\bel{SGT}
t_1\le s_1\le 0.
\ee

\par Recall that
\bel{X12-N}
X_1=\frac{\delta}{\Delta}\,\|\lh_2\|_1\,(-\sin \vf_1,\cos \vf_1)~~~~\text{and}~~~~
X_2=\frac{\delta}{\Delta}\,\|\lh_1\|_1\,
(\sin \vf_2,-\cos \vf_2)
\ee
where $\Delta=\sin(\vf_2-\vf_1)$. See \rf{X-1} and \rf{X-2}.
\msk
\par Note that  $\Delta=\sin(\vf_2-\vf_1)>0$ (because $\vf_2-\vf_1\in(0,\pi)$) and  $\sin \vf_1>0$ (because $\vf_1\in(\pi/2,\pi))$. Therefore, thanks to the first formula in \rf{X12-N},
$$
s_1=-(\delta/\Delta)\,\|\lh_2\|_1\,\sin \vf_1\le 0.
$$
\par Let us see that
\bel{ST-F}
s_1-t_1=\frac{\delta}{\Delta}\sin(\vf_1+\vf_2)=
\frac{\delta}{\Delta}\sin(\alpha_1-\alpha_2).
\ee
\par Indeed, $\lh_i=(\cos \vf_i,\sin \vf_i)$,~$i=1,2$, and, thanks to \rf{FI-C1},
$$
\vf_1\in (\pi/2,\pi),~~~~\vf_2\in (\pi,3/2\pi)
$$
so that
$$
\|\lh_1\|_1= -\cos \vf_1+\sin \vf_1~~~\text{and}~~~
\|\lh_2\|_1= -\cos \vf_2-\sin \vf_2.
$$
Hence, thanks to formulae \rf{X12-N},
\be
s_1-t_1&=&(\delta/\Delta)
\{-\|\lh_2\|_1\sin\vf_1-\|\lh_1\|_1\sin\vf_2\}
\nn\\
&=&
(\delta/\Delta)\{-(-\cos \vf_2-\sin \vf_2)\sin\vf_1-
(-\cos \vf_1+\sin \vf_1)\sin\vf_2\}
\nn\\
&=&
(\delta/\Delta)
\{\cos \vf_2\sin\vf_1+\cos \vf_1\sin\vf_2\}
=(\delta/\Delta)\sin(\vf_1+\vf_2).
\nn
\ee
\par Recall that, thanks to \rf{AL-12},
$$
\alpha_1=\vf_1-\pi/2~~~~\text{and}~~~~ \alpha_2=\tfrac32\pi-\vf_2
$$
so that
$$
\vf_1+\vf_2=\alpha_1-\alpha_2+2\pi
$$
proving \rf{ST-F}.
\par It remains to note that $\sin(\alpha_1-\alpha_2)>0$ (because $0<\alpha_2\le\alpha_1<\pi/2$) and $\Delta>0$ so that, thanks \rf{ST-F}, $t_1\le s_1$ proving \rf{SGT}.
\smsk
\smsk
\par In particular, this inequality implies the inclusion
\bel{HR-V}
\HR[S_1]\subset \HR[S_2]
\ee
as it shown on Fig. 27.
\par Let us recall property \rf{A-HR1N}:
$$
a_1\in \Omega=\Hc[S_1]\setminus (H_1\cup S_2).
$$

\par Let $\Gamma=\Delta(Y,X_1,O)$ be the triangle with vertices in the points $Y$, $X_1$ and $O$ where
$$
Y=\ell_2\cap \hat{\ell}_1
$$
is the point of intersection of the line $\ell_2$ (the boundary of $H_2$) and the line $\hat{\ell}_1$ defined by \rf{LH-1}. Let
\bel{DF-PS}
\Psi=\HR[S_1]\setminus (H_1+\delta Q_0).
\ee
Clearly, $\Psi$ is a cone in $\RT$ with the vertices at the point $P$ determined by the straight lines  $\hat{\ell}_1$ and $\tell_1$.
\msk
\par Then the set $\Omega$ is a disjoint union of the triangle $\Gamma=\Delta(Y,X_1,O)$ and the cone $\Psi$, i.e.,
$$
\Omega=\Gamma\cup\Psi~~~~\text{and}~~~~\Gamma\cap\Psi=\emp.
$$
See Fig. 27.
\smsk
\par Let assume that $a_1\in \Gamma$. Since
$$
\Gamma=\Delta(Y,X_1,O)\subset T=\Delta(X_1,X_2,O),
$$
the point $a_1\in T$. Lemma \reff{TR-X12} tells us that in this case inequality \rf{PRA-D} holds.
\par Now, let $a_1\in \Psi$. Inclusion \rf{HR-V} and definition \rf{DF-PS} imply the following:
$$
a_1\in \Psi=\HR[S_1]\setminus (H_1+\delta Q_0)\subset
\HR[S_2]\setminus (H_1+\delta Q_0).
$$
\par Lemma \reff{FL-P} tells us that in this case \rf{PRA-D} holds completing the proof of the lemma.\bx
\msk
\par We turn to the proof of \rf{PRA-D} in {\bf Case 3}.

\begin{lemma}\lbl{L7-P} Inequality \rf{PRA-D} holds provided the cone $S=H_1\cap H_2$ does not contain any coordinate semi-axis.
\end{lemma}
\par {\it Proof.} Without loss of generality we may assume that
$$
S\subset \RT_+=Ox_1^+\times Ox_2^+,
$$
i.e., $S=H_1\cap H_2$ is located inside of the first coordinate quadrant. We can also assume that in $\RT_+$, the straight line $\ell_1=\partial H_1$ is located above the straight line $\ell_2=\partial H_2$, as it shown on Fig. 28.
\bsk\msk
\begin{figure}[H]
\hspace*{14mm}
\includegraphics[scale=0.65]{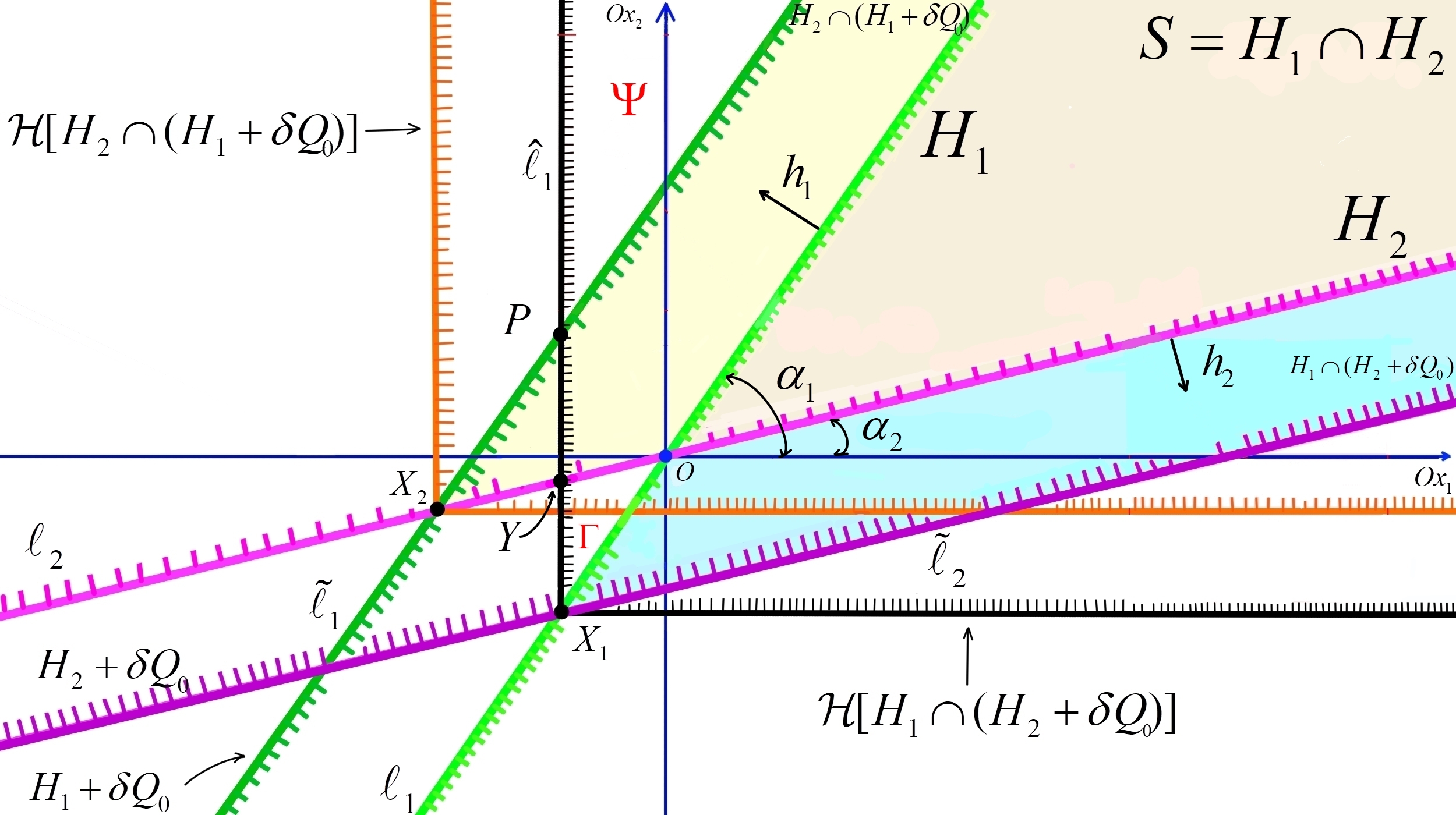}
\hspace*{12mm}
\caption{The half-planes $H_1$ and $H_2$ with the non-parallel boundaries: {\bf Case 3}.}
\end{figure}
\msk

\par Thus,
$$
0<\alpha_2<\alpha_1<\pi/2
$$
where, as in the previous lemma, $\alpha_i\in(0,\pi/2)$, $i=1,2$, is the angle between the straight line $\ell_i$ and the axis $Ox_1$.

\par Let
$$
\lh_1=(\cos\vf_1,\sin\vf_1)~~~\text{and}~~~
\lh_2=(\cos\vf_2,\sin\vf_2)
$$
be the vectors determining the half-spaces $H_1$ and $H_2$
respectively. See \rf{H12-DFN} and \rf{H-12}. In the case under consideration, we have
\bel{FI-12}
\vf_1\in (\pi/2,\pi)~~~~\text{and}~~~~
\vf_2\in (3/2\pi,2\pi).
\ee
Furthermore,
\bel{AL-C2}
\vf_1=\alpha_1+\pi/2~~~~\text{and}~~~~
\vf_2=\alpha_2+3/2\pi.
\ee
\par Moreover, the rectangular hulls of the sets $S_1=H_1\cap (H_2+\delta Q_0)$ and  $S_2=H_2\cap (H_1+\delta Q_0)$ have the following representations:
$$
\HR[S_1]=\{(u_1,u_2)\in\RT:u_1\ge s_1,\, u_2\ge s_2\}
~~~\text{and}~~~
\HR[S_2]=\{(u_1,u_2)\in\RT:u_1\ge t_1,\, u_2\ge t_2\}.
$$

\par We recall that $X_1=(s_1,s_2)$ and $X_2=(t_1,t_2)$ are the vertices of the cones $S_1$ and $S_2$ respectively. We also recall the formulae \rf{X-1} and \rf{X-2} for the points $X_1=(s_1,s_2)$ and $X_2=(t_1,t_2)$:
$$
X_1=\frac{\delta}{\Delta}\,\|\lh_2\|_1\,(-\sin \vf_1,\cos \vf_1)~~~~\text{and}~~~~
X_2=\frac{\delta}{\Delta}\,\|\lh_1\|_1\,
(\sin \vf_2,-\cos \vf_2)
$$
where $\Delta=\sin(\vf_2-\vf_1)$.
\par Using these formulae, let us show that
\bel{X1-X2}
X_1-X_2=\delta(1,-1).
\ee

\par Indeed, we know that $\vf_1\in (\pi/2,\pi)$ and $\vf_2\in (3/2\pi,2\pi)$. See \rf{FI-12}. Therefore,
$$
\|\lh_1\|_1=|\cos\vf_1|+|\sin\vf_1|=
-\cos\vf_1+\sin\vf_1,~~~
\|\lh_2\|_1=|\cos\vf_2|+|\sin\vf_2|=
\cos\vf_2-\sin\vf_2.
$$
Therefore,
$$
X_1=\frac{\delta}{\Delta}\,(\cos\vf_2-\sin\vf_2)\,(-\sin \vf_1,\cos \vf_1)
~~~\text{and}~~~
X_2=\frac{\delta}{\Delta}\,(-\cos\vf_1+\sin\vf_1)\,
(\sin \vf_2,-\cos \vf_2).
$$
Hence,
\be
X_1-X_2&=&
\frac{\delta}{\Delta}\,((\cos\vf_2-\sin\vf_2)\,(-\sin \vf_1)-(-\cos\vf_1+\sin\vf_1)\sin \vf_2,
\nn\\
&&
(\cos\vf_2-\sin\vf_2)
\cos \vf_1-(-\cos\vf_1+\sin\vf_1)(-\cos \vf_2))
\nn\\
&=&
\frac{\delta}{\Delta}\,(-\cos\vf_2\sin\vf_1
+\cos\vf_1\sin\vf_2,
-\sin\vf_2\cos\vf_1+\sin\vf_1\cos\vf_2)
\nn\\
&=&\frac{\delta}{\Delta}\sin(\vf_2-\vf_1)(1,-1).
\nn
\ee
\par Thanks to \rf{A-DL}, $\Delta=\sin(\vf_2-\vf_1)$, and the proof of \rf{X1-X2} is complete.
\msk
\par Thanks to this equality, $t_1=s_1-\delta$ and $t_2=s_2+\delta$. Furthermore, \rf{X1-X2} and \rf{S-X12} imply the following:
$$
S_2=S_1+\delta(-1,1)~~~\text{and}~~~
\HR[S_2]=\HR[S_1]+\delta(-1,1).
$$
\par Let us prove that
\bel{ST-CM}
t_1\le s_1\le 0~~~~\text{and}~~~~s_2\le t_2\le 0.
\ee
\par In fact, we know that $\vf_1\in(\pi/2,\pi)$ and $\vf_2\in(3/2\,\pi,2\pi)$, see \rf{FI-12}. Hence, $\pi/2<\vf_2-\vf_1$.
\par On the other hand, thanks to \rf{AL-C2},
$$
\alpha_2=\vf_2-3/2\pi<\alpha_1=\vf_1-\pi/2
$$
proving that $\vf_2-\vf_1\in(\pi/2,\pi)$. Hence,
$\Delta=\sin(\vf_2-\vf_1)>0$.
\par Let also note that $\sin \vf_1>0$ (because $\vf_1\in(\pi/2,\pi)$) and, thanks to \rf{X-1},
$$
s_1=-(\delta/\Delta)\,\|\lh_2\|_1\,\sin \vf_1.
$$
Hence, $s_1\le 0$. In addition,
$$
t_1=s_1-\delta\le s_1
$$
proving the first inequality in \rf{ST-CM}.
\smsk
\par Next, thanks to \rf{X-2},
$$
t_2=-(\delta/\Delta)\,\|\lh_1\|_1\,\cos\vf_2.
$$
But $\Delta>0$ and $\cos\vf_2>0$ (because $\vf_2\in(3/2\,\pi,2\pi)$) so that $t_2\le 0$. Moreover, $s_2=t_2-\delta$, so that $s_2\le t_2$, and the proof of \rf{ST-CM} is complete.
\smsk

\par Recall that
$$
a_1\in \Omega=\Hc[S_1]\setminus (H_1\cup S_2).
$$
(Recall also that $S_1=H_1\cap(H_2+\delta Q_0)$ and
$S_2=H_2\cap(H_1+\delta Q_0)$.)
\smsk
\par We complete the proof of the lemma in the same way as we have completed the proof of Lemma \reff{L6-P}.
\par Let
$$
\Gamma=\Delta(Y,X_1,O)
$$
be the triangle with vertices in the points $Y$, $X_1$ and $O$ where $Y$ is the point of intersection of the line $\ell_2$ (the boundary of $H_2$) and the line $\hat{\ell}_1$ defined by \rf{LH-1}. Let
$$
\Psi=\HR[S_1]\setminus (H_1+\delta Q_0).
$$
Clearly, $\Psi$ is a cone in $\RT$ with the vertices at the point $P$ determined by the straight lines  $\hat{\ell}_1$ and $\tell_1$.
\msk
\par Then the set $\Omega$ is a disjoint union of the triangle $\Gamma=\Delta(Y,X_1,O)$ and the cone $\Psi$, i.e.,
$$
\Omega=\Gamma\cup\Psi~~~~\text{and}~~~~\Gamma\cap\Psi=\emp.
$$
See Fig. 28.
\smsk
\par If $a_1\in \Gamma$, then $a_1\in T=\Delta(O,X_1,X_2)$. Lemma \reff{TR-X12} tells us that in this case inequality \rf{PRA-D} holds.
\par Now, let $a_1\in \Psi$. Then
$$
a_1\in \Psi=\HR[S_1]\setminus (H_1+\delta Q_0)\subset
\HR[S_2]\setminus (H_1+\delta Q_0).
$$
\par Lemma \reff{FL-P} tells us that in this case inequality \rf{PRA-D} holds as well, and the proof of Lemma \reff{L7-P} is complete.\bx
\msk
\par Finally, the results of Lemmas \reff{L1-P} -- \reff{L7-P} imply the required inequality \rf{PRA-D} completing the proof of Proposition \reff{TWO-HP}.\bx

\msk
\par We are in a position to complete the proof of Proposition \reff{LSEL-F}.
\smsk
\par \underline{\sc Proof of inequality \rf{LIP-FL}.} Recall that we have a mapping $g:\Mc\to\RT$ such that
\bel{WG-HFN}
g(x)\in \Tc_{F,\tlm}(x)=\HR[F^{[1]}[x:\tlm]]~~~~~\text{for every}~~~~~x\in\Mc,
\ee
and
\bel{G-LIPN}
\|g(x)-g(y)\|\le \lambda\,\rho(x,y)~~~~~\text{for all}~~~~~x,y\in\Mc.
\ee
\par We recall that, thanks to Lemma \reff{WF1-NE},
$$
F^{[1]}[x:\tlm]\ne\emp~~~~\text{for every}~~~~x\in\Mc.
$$
\par We define a mapping
\bel{FX-PRN}
f(x)=\Prm\left(g(x),F^{[1]}[x:\tlm]\right),~~~~~~x\in\Mc.
\ee
\par We have to prove that
$f$ is well defined and
$$
\|f\|_{\Lip(\Mc)}\le 2\lambda+\tlm.
$$
\par Let us fix elements $x,y\in\Mc$. We set
$$
A_1=F^{[1]}[x:\tlm],~~~A_2=F^{[1]}[y:\tlm]
$$
and
$$
a_1=g(x),~~~~~a_2=g(y).
$$
\par We also recall that the mapping $F^{[1]}[\cdot:\tlm]$ is defined by \rf{WF1-3L}. Thus,
\bel{A-R12}
A_i=\bigcap_{u\in\Mc}\, A_i^{[u]}\,,~~~~~i=1,2,
\ee
where given $u\in\Mc$, we set
\bel{AU-12}
A_1^{[u]}=F(u)+\tlm\,\rho(u,x)Q_0
~~~\text{and}~~~
A_2^{[u]}=F(u)+\tlm\,\rho(u,y)Q_0.
\ee
\par Lemma \reff{WF1-NE} tells us that
$$
A_1\ne\emp~~~~\text{and}~~~~A_2\ne\emp.
$$
\par Thanks to inequality \rf{G-LIPN}, we have
\bel{A-12}
\|a_1-a_2\|\le\lambda\,\rho(x,y),
\ee
and, thanks to \rf{WG-HFN},
\bel{HA-12}
a_1\in\HR[A_1]~~~~~\text{and}~~~~~a_2\in\HR[A_2].
\ee
Furthermore, formula \rf{FX-PRN} tells us that
$$
f(x)=\Prm(a_1,A_1)~~~~\text{and}~~~~f(y)=\Prm(a_2,A_2).
$$
(We also recall that, thanks to Lemma \reff{WMP-S},
the metric projection $\Prm(a_i,A_i)$ is well defined, i.e., $\Prm(a_i,A_i)$ is a singleton.)
\par In these settings, the required inequality \rf{LIP-FL}
reads as follows:
\bel{NV-12}
\|\Prm(a_1,A_1)-\Prm(a_2,A_2)\|\le (2\lambda+\tlm)\,\rho(x,y).
\ee
\par Let us note that this inequality is immediate from \rf{A-12} provided $a_i\in A_i$, $i=1,2$, (because in this case $\Prm(a_i,A_i)=a_i$).
\par Suppose that
$$
\text{either}~~~a_1\notin A_1~~~\text{or}~~~
a_2\notin A_2.
$$
\par Without loss of generality, we may assume that $a_1\notin A_1$. Fix $\ve>0$ and prove that there exists a half-plane $H_1\in\HPL$ such that
\bel{H1-A1}
H_1\supset A_1,~~~~H_1+\tlm\,\rho(x,y)Q_0\supset A_2,
\ee
and
\bel{D-A1}
\|\Prm(a_1,A_1)-\Prm(a_1,H_1)\|<\ve.
\ee
\par We construct the half-plane $H_1$ as follows: Because $a_1\notin A_1$, we have $\Prm(a_1,A_1)\ne a_1$ so that
$(\Prm(a_1,A_1),a_1]$ is a nonempty semi-open line segment in $\RT$. Let us pick a point
$$
a^{(\ve)}\in(\Prm(a_1,A_1),a_1]
$$
such that
\bel{H-EP1}
\|a^{(\ve)}-\Prm(a_1,A_1)\|<\ve.
\ee
See Fig. 29.

\begin{figure}[H]
\hspace{9mm}
\includegraphics[scale=0.7]{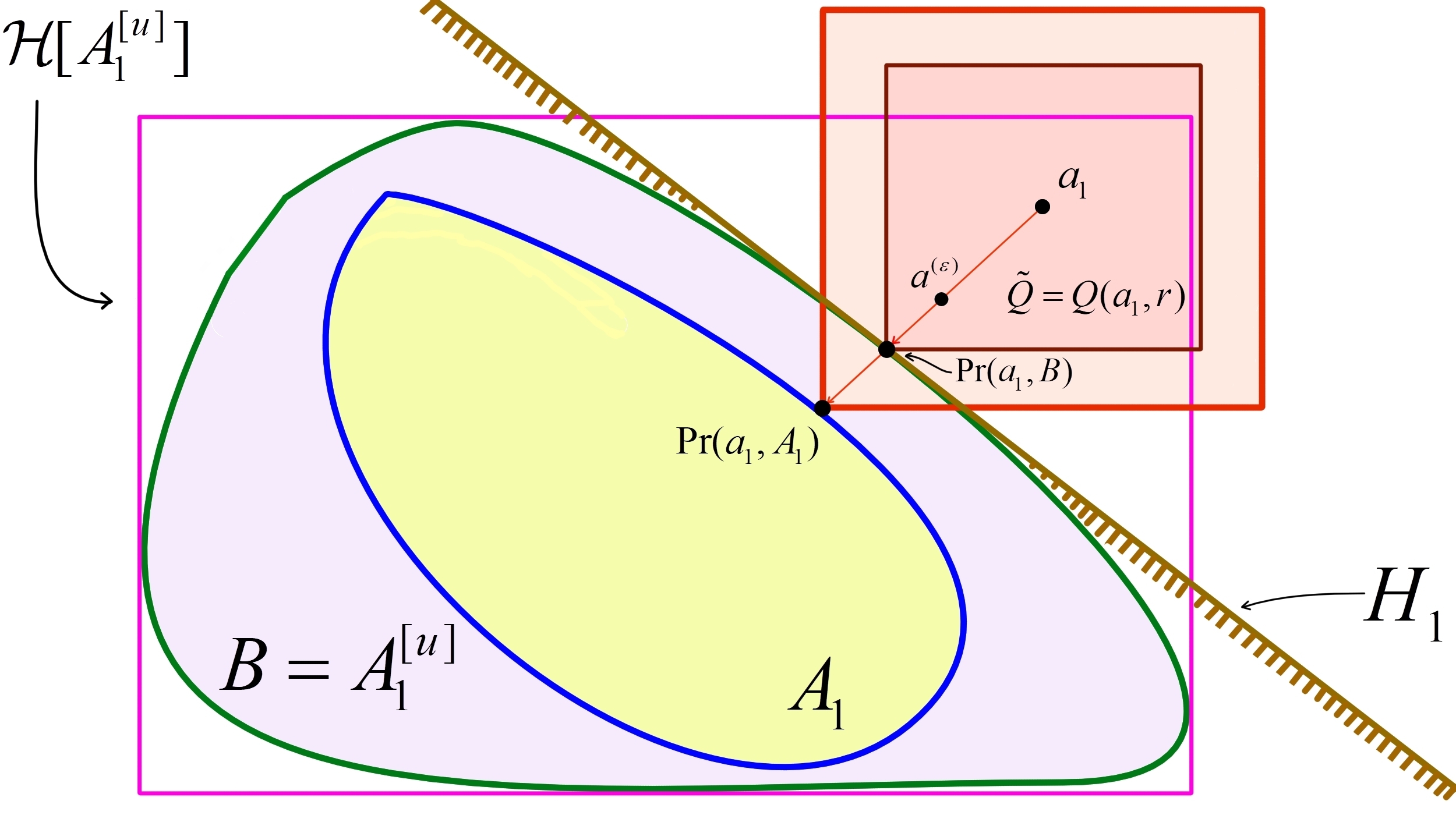}
\hspace*{10mm}
\caption{Metric projections of $a_1$ onto $A_1$ and $B$.}
\end{figure}
\msk
\par Because $\Prm(a_1,A_1)$ is the nearest to $a_1$ point on $A_1$, we have
$$
(\Prm(a_1,A_1),a_1]\cap A_1=\emp.
$$
Therefore,
$$
a^{(\ve)}\notin A_1=\bigcap_{u\in\Mc}\,A_1^{[u]}\,.
$$
\par See \rf{A-R12}. This implies the existence of an element $u\in\Mc$ such that $a^{(\ve)}\notin A_1^{[u]}$.
\smsk

\par  Let $B=A_1^{[u]}$. Thus,
\bel{AE-B}
a^{(\ve)}\notin B=A_1^{[u]}=F(u)+\tlm\,\rho(u,x) Q_0.
\ee
See \rf{AU-12}. Thanks to \rf{HA-12} and \rf{A-R12},
\bel{B-A1}
a_1\in\HR[A_1]~~~~\text{and}~~~~A_1\subset B.
\ee
Therefore, thanks to Lemma \reff{AB-PR}, the metric projections $\Prm(a_1,A_1)$ and $\Prm(a_1,B)$ are singletons such that
$$
\Prm(a_1,B)\in[\Prm(a_1,A_1),a_1].
$$

\par Moreover, $\Prm(a_1,B)\in [\Prm(a_1,A_1),a^{(\ve)}]$; indeed, otherwise $a^{(\ve)}\in [\Prm(a_1,A_1),\Prm(a_1,B)]\subset B$, a contradiction. See \rf{AE-B}. Hence, thanks to \rf{H-EP1},
\bel{PR-E1}
\|\Prm(a_1,A_1)-\Prm(a_1,B)\|\le \|a^{(\ve)}-\Prm(a_1,A_1)\|<\ve.
\ee
\par Let
$$
\tQ=Q(a_1,r)~~~~~\text{where}~~~~r=\dist(a_1,B).
$$
Thus,
$$
\tQ\cap B=\Prm(a_1,B).
$$
Therefore, thanks to the separation theorem, there exists a half-plane
$$
H_1\in\HPL~~\text{which contains}~~B~~
\text{and separates (not strictly)}~~\tQ~~\text{and}~~ B.
$$
\par Thus, $B\subset H_1$  and $\tQ\cap H_1=\Prm(a_1,B)$
as it shown on Fig. 29. In particular, these properties imply the equality
\bel{P-HB}
\Prm(a_1,H_1)=\Prm(a_1,B).
\ee
\par Let us see that inclusions \rf{H1-A1} and inequality \rf{D-A1} hold for the half-plane $H_1$, i.e.,
\bel{H1-A1N}
H_1\supset A_1,~~~~H_1+\tlm\,\rho(x,y)Q_0\supset A_2,
\ee
and
\bel{D-A1N}
\|\Prm(a_1,A_1)-\Prm(a_1,H_1)\|<\ve.
\ee
\par In fact, \rf{D-A1N} is immediate from \rf{PR-E1} and \rf{P-HB}.
\par Prove \rf{H1-A1N}. We know that $A_1\subset B$, see \rf{B-A1}, so that $A_1\subset B\subset H_1$. We also recall that
$$
B=F(u)+\tlm\,\rho(u,x) Q_0
$$
(see \rf{AE-B}). Therefore,
$$
H_1+\tlm\,\rho(x,y)Q_0\supset B+\tlm\,\rho(x,y)Q_0=
F(u)+\tlm\,\rho(u,x) Q_0+\tlm\,\rho(x,y)Q_0=
F(u)+\tlm\,(\rho(u,x)+\rho(x,y))\,Q_0.
$$
Therefore, thanks to the triangle inequality, \rf{A-R12} and \rf{AU-12}, we have
$$
H_1+\tlm\,\rho(x,y)Q_0\supset
F(u)+\tlm\,\rho(u,y) \,Q_0=A_2^{[u]}\supset A_2
$$
proving \rf{H1-A1N}.
\msk
\par Next, let us construct a half-plane $H_2\in\HPL$ having the following properties:
\bel{H2-A2}
H_2\supset A_2,~~~~H_2+\tlm\,\rho(x,y)Q_0\supset A_1,
\ee
and
\bel{D-A2}
\|\Prm(a_2,A_2)-\Prm(a_2,H_2)\|<\ve.
\ee
\par If $a_2\notin A_2$, we define $H_2$ in the same way as we have defined $H_1$ for $a_1$. In this case, properties  \rf{H2-A2}, \rf{D-A2} are a complete analog of properties
\rf{H1-A1N} and \rf{D-A1N} obtained for the point $a_1$.
\par If $a_2\in A_2$, we set
$$
H_2=H_1+\tlm\,\rho(x,y)Q_0.
$$
\par Clearly, $H_2$ is a half-plane. Let us see that inclusions \rf{H2-A2} and inequality \rf{D-A2} hold for this choice of $H_2$. Indeed, thanks to the second inclusion in \rf{H1-A1N}, we have $H_2\supset A_2$. In turn, thanks to the first inclusion,
$$
H_2+\tlm\,\rho(x,y)Q_0=(H_1+\tlm\,\rho(x,y)Q_0)
+\tlm\,\rho(x,y)Q_0\supset H_1\supset A_1,
$$
proving \rf{H2-A2}. Finally, inequality \rf{D-A2} is trivial because $\Prm(a_2,A_2)=\Prm(a_2,H_2)(=a_2)$. (Recall that $a_2\in A_2\subset H_2$.)
\msk
\par Now, we set
\bel{DL}
\delta=\tlm\,\rho(x,y)+\ve.
\ee
\par Let us prove that the points $a_1,a_2$ and the half-planes $H_1$ and $H_2$ satisfy conditions of the hypothesis of Proposition \reff{TWO-HP}, i.e., conditions \rf{H12-DF}, \rf{A-I} and \rf{A-II}.
\msk
\par Thanks to \rf{H1-A1}, $H_1\supset A_1$, and, thanks to \rf{H2-A2}, $H_2+\tlm\,\rho(x,y)Q_0\supset A_1$. Hence,
$$
H_1\cap(H_2+\tlm\,\rho(x,y)Q_0)\supset A_1.
$$
\par Note that $\tlm\,\rho(x,y)\le\delta$; see \rf{DL}. Therefore,
\bel{WW}
H_1\cap(H_2+\delta\,Q_0)\supset A_1.
\ee
\par We also know that $A_1\ne\emp$ so that  $H_1\cap(H_2+\delta\,Q_0)\ne\emp$ as well. This proves that  
$$
\dist(H_1,H_2)\le\delta
$$
so that condition \rf{H12-DF} is satisfied.
\smsk
\par Let us prove that the points $a_1$ and $a_2$ satisfy condition \rf{A-I}. Indeed, inclusion \rf{WW} tells us that
$$
\HR[A_1]\subset \HR[H_1\cap(H_2+\delta\,Q_0)].
$$
But, thanks to \rf{HA-12}, $a_1\in \HR[A_1]$ proving that $a_1\in \HR[H_1\cap(H_2+\delta\,Q_0)]$.
\par In the same fashion we show that $a_2\in \HR[H_2\cap(H_1+\delta\,Q_0)]$ proving that condition \rf{A-I} holds.
\smsk
\par Let us show that condition \rf{A-II} is satisfied as well.
\par Thanks to \rf{H2-A2}, $A_1\subset H_2+\tlm\,\rho(x,y)\,Q_0$, so that
$$
\Prm(a_1,A_1)\in A_1\subset H_2+\tlm\,\rho(x,y)\,Q_0.
$$
Thanks to \rf{D-A1}, $\|\Prm(a_1,A_1)-\Prm(a_1,H_1)\|<\ve$
so that
$$
\Prm(a_1,H_1)\in \Prm(a_1,A_1)+\ve\,Q_0\subset
A_1+\ve\,Q_0.
$$
Recall that $\delta=\tlm\,\rho(x,y)+\ve$. Therefore,
$$
\Prm(a_1,H_1)\in A_1+\ve\,Q_0\subset (H_2+\tlm\,\rho(x,y)\,Q_0)+\ve\,Q_0=
H_2+\delta\,Q_0.
$$
\par In the same way we show that $\Prm(a_2,H_2)\in
H_1+\delta\,Q_0$ completing the proof of \rf{A-II}.
\msk
\par Therefore, thanks to Proposition \reff{TWO-HP}, inequality \rf{PRA-D} holds. This inequality together with \rf{A-12} and \rf{DL} imply the following:
$$
\|\Prm(a_1,H_1)-\Prm(a_2,H_2)\|\le 2\|a_1-a_2\|+\delta
\le 2\lambda\,\rho(x,y)+ \tlm\,\rho(x,y)+\ve.
$$
From this, \rf{D-A1} and \rf{D-A2}, we have
\be
\|\Prm(a_1,A_1)-\Prm(a_2,A_2)\|&\le&
\|\Prm(a_1,A_1)-\Prm(a_1,H_1)\|+
\|\Prm(a_1,H_1)-\Prm(a_2,H_2)\|
\nn\\
&+&\|\Prm(a_2,H_2)-\Prm(a_2,A_2)\|
\le (2\lambda+\tlm)\,\rho(x,y)+3\ve.
\nn
\ee
\par Since $\ve>0$ is arbitrary, this implies \rf{NV-12} proving the required inequality \rf{LIP-FL} and completing the proof of Proposition \reff{LSEL-F}.\bx
\smsk
\par Finally, combining part (i) of Proposition \reff{WLS-T} with Proposition \reff{LSEL-F}, we conclude that the mapping $f$ defined by formula \rf{FX-PRN}, is
a Lipschitz selection of the set-valued mapping $F$ with Lipschitz seminorm at most $2\lambda+\tlm$.
\par The proof of Theorem \reff{W-CR} is complete.\bx

\SECT{5. Main properties of the Projection Algorithm: proofs.}{5}

\addtocontents{toc}{5. The $\lmv$-Projection Algorithm and its versions.\hfill\thepage\par \VST}

\indent

\par {\bf 5.1 Main properties of the Projection Algorithm: proofs.}
\addtocontents{toc}{~~~~5.1 Main properties of the Projection Algorithm: proofs.\hfill \thepage\par\VST}

\msk
\par \underline{\sc Proof of Theorem \reff{ALG-T}.} {\it (i)} Let us assume that the $\lmv$-Projection Algorithm produces the outcome {\bf ``No go''} and prove that there does not exist a Lipschitz selection of $F$ with Lipschitz seminorm at most $\lambda^{\min}=\min\{\lambda_1,\lambda_2\}$. We recall that, thanks to \rf{NG-O}, in this case there exists a point $\tx\in\Mc$ such that
$$
\text{either} ~~F^{[1]}[\tx:\lambda_1]=\emp~~\text{or}~~
\Tc^{[1]}_{F,\lambda_1}[\tx:\lambda_2]=\emp.
$$
\par Let us suppose that $F$ has a Lipschitz selection $f:\Mc\to\RT$ with $\|f\|_{\Lip(\Mc)}\le\lambda^{\min}$, and prove that
$F^{[1]}[\tx:\lambda_1]\ne\emp$ and     $\Tc^{[1]}_{F,\lambda_1}[\tx:\lambda_2]\ne\emp$.
\par Indeed, we know that $f(x)\in F(x)$ and
$$
\|f(x)-f(y)\|\le\lambda^{\min}\,\rho(x,y)\le
\lambda_1\,\rho(x,y)
~~~\text{for all }~~~x,y\in\Mc.
$$
Therefore,
$$
f(x)\in F(y)+\lambda_1\,\rho(x,y)Q_0~~~~\text{for all}~~~~ x,y\in\Mc.
$$
Hence,
$$
f(x)\in F^{[1]}[x:\lambda_1]~~~~\text{for every}~~~~ x\in\Mc.
$$
See  definition \rf{F-1}. In particular,
$f(\tx)\in F^{[1]}[\tx:\lambda_1]$ proving that $F^{[1]}[\tx:\lambda_1]\ne\emp$.
\msk
\par Let us show that $\Tc^{[1]}_{F,\lambda_1}[\tx:\lambda_2]\ne\emp$ as well. Indeed,
$$
f(x)\in F^{[1]}[x:\lambda_1]\subset \Hc[F^{[1]}[x:\lambda_1]]=\Tc_{F,\lambda_1}(x)~~~
\text{for every}~~~x\in\Mc.
$$
See \rf{TC-DF}. Therefore, $f(\tx)\in\Tc_{F,\lambda_1}(\tx)$,
$f(y)\in \Tc_{F,\lambda_1}(y)$ and
$$
\|f(\tx)-f(y)\|\le\lambda^{\min}\,\rho(\tx,y)\le
\lambda_2\,\rho(\tx,y)~~~\text{for every}~~~y\in\Mc.
$$
Hence,
$$
f(\tx)\in \Tc_{F,\lambda_1}(y)+\lambda_2\,\rho(\tx,y)Q_0
~~~\text{for every}~~~y\in\Mc,
$$
so that, thanks to \rf{RT-1},
$f(\tx)\in \Tc^{[1]}_{F,\lambda_1}[\tx:\lambda_2]$. Hence, $\Tc^{[1]}_{F,\lambda_1}[\tx:\lambda_2]\ne\emp$, a contradiction.
\msk
\par {\it (ii)}  Let us assume that the  $(\lambda_1,\lambda_2)$-PA produces the outcome {\bf ``Success''}. Then, thanks to \rf{SC-P}, for every $x\in\Mc$, we have
\bel{S-P5}
F^{[1]}[x:\lambda_1]\neq\emp~~~\text{and}~~~
\Tc^{[1]}_{F,\lambda_1}[x:\lambda_2]\neq\emp.
\ee
\par We recall that $F^{[1]}[x:\lambda_1]$ is
the $\lambda_1$-metric refinement of $F$,
$$
F^{[1]}[x:\lambda_1]=\bigcap_{y\in \Mc}\,
\left[F(y)+\lambda_1\,\rho(x,y)\,Q_0\right],
$$
and $\Tc_{F,\lambda_1}(x)=\Hc[F^{[1]}[x:\lambda_1]]$. We also recall that
$$
\Tc^{[1]}_{F,\lambda_1}[x:\lambda_2]=\bigcap_{y\in \Mc}\,
\left[\Tc_{F,\lambda_1}(y)+\lambda_2\,\rho(x,y)\,Q_0\right]
$$
is the $\lambda_2$-metric refinement of the mapping $\Tc_{F,\lambda_1}$.
\smsk
\par Let us prove property ($\bigstar A$). Clearly, part {\it (a)} of this property is immediate from \rf{S-P5}.
\msk
\par Prove part {\it (b)}. We recall that
\bel{G-CNTN}
g_F(x)=\cent\left
(\Prm\left(O,\Tc^{[1]}_{F,\lambda_1}[x:\lambda_2]\right) \right), ~~~~~~x\in\Mc.
\ee
\par Note that for every $x\in\Mc$, the metric projection of the origin $O$ onto the rectangle $\Tc^{[1]}_{F,\lambda_1}[x:\lambda_2]$, is either a singleton or a compact line segment parallel to one of the coordinate axes. Therefore, the point $g_F(x)$, the center of the line segment $\Prm(O,\Tc^{[1]}_{F,\lambda_1}[x:\lambda_2])$, see \rf{G-CNT}, is well defined.
\msk
\par Let us prove part {\it (c)}. Recall that
$$
f_{\lmv;F}(x)=\Prm(g_F(x),F^{[1]}[x:\lambda_1]), ~~~~~~x\in\Mc.
$$
\par Thanks to \rf{TC-DF}, \rf{RT-1} and \rf{G-CNTN},
\bel{G-IN-T}
g_F(x)\in\Tc^{[1]}_{F,\lambda_1}[x:\lambda_2]
\subset\Tc_{F,\lambda_1}(x)
=\Hc[F^{[1]}[x:\lambda_1]].
\ee
Therefore, thanks to Lemma \reff{WMP-S}, $f_{\lmv;F}(x)$, the metric projection of $g_F(x)$ onto $F^{[1]}[x:\lambda_1]$, see \rf{F-LF}, is a singleton, and the proof of property ($\bigstar A$) is complete.
\smsk
\par Let us prove property ($\bigstar B$) of the theorem.
We will follow the scheme of the proof of the key Theorem \reff{W-CR} for the special case
\bel{LM-12}
\tlm=\lambda_1~~~~\text{and}~~~~\lambda=\lambda_2.
\ee
\par Thanks to \rf{S-P5},
$$
\Tc^{[1]}_{F,\tlm}[x:\lambda]\ne\emp~~~~\text{for every}~~~~x\in\Mc.
$$
\par Part {\it (b)} of Proposition \reff{X2-C} tells us that in this case the mapping $\Tc^{[1]}_{F,\tlm}[\cdot:\lambda]$ is Lipschitz with respect to the Hausdorff distance. See \rf{TAM-1}. Thus,
\bel{T-HD}
\dhf\left(\Tc^{[1]}_{F,\tlm}[x:\lambda],
\Tc^{[1]}_{F,\tlm}[y:\lambda]\right)
\le \lambda\,\rho(x,y)~~~\text{for all}~~~x,y\in\Mc.
\ee
\par Let
$$
\delta(x)=\dist(O,\Tc^{[1]}_{F,\tlm}[x:\lambda]).
$$
(Recall that $O=(0,0)$ is the origin.) It is known that the distance from a fixed point to a set is a Lipschitz function of the set (with respect to the Hausdorff distance). Therefore,
$$
|\delta(x)-\delta(y)|=
|\dist(O,\Tc^{[1]}_{F,\tlm}[x:\lambda])-
\dist(O,\Tc^{[1]}_{F,\tlm}[y:\lambda])|
\le \dhf\left(\Tc^{[1]}_{F,\tlm}[x:\lambda],
\Tc^{[1]}_{F,\tlm}[y:\lambda]\right)
$$
so that, thanks to \rf{T-HD},
\bel{L-DLT}
|\delta(x)-\delta(y)|\le \lambda\,\rho(x,y)~~~
\text{for every}~~~x,y\in\Mc.
\ee
\par We note that
$$
\Prm(O,\Tc^{[1]}_{F,\tlm}[x:\lambda])=
Q(O,\delta(x))\cap \Tc^{[1]}_{F,\tlm}[x:\lambda].
$$
From this and Lemma \reff{H-RB}, we have
\be
\Prm(O,\Tc^{[1]}_{F,\tlm}[x:\lambda])+\lambda\,\rho(x,y)Q_0
&=&
\{Q(O,\delta(x))+\lambda\,\rho(x,y)Q_0\}\cap \{\Tc^{[1]}_{F,\tlm}[x:\lambda]+\lambda\,\rho(x,y)Q_0\}
\nn\\
&=&
Q(O,\delta(x)+\lambda\,\rho(x,y))\cap \{\Tc^{[1]}_{F,\tlm}[x:\lambda]+\lambda\,\rho(x,y)Q_0\}.
\nn
\ee
\par Note that, thanks to \rf{L-DLT},
$$
\delta(y)\le\delta(x)+\lambda\,\rho(x,y),
$$
and, thanks to \rf{T-HD},
$$
\Tc^{[1]}_{F,\tlm}[y:\lambda]\subset\
\Tc^{[1]}_{F,\tlm}[x:\lambda]+\lambda\,\rho(x,y)Q_0.
$$
Hence,
$$
\Prm(O,\Tc^{[1]}_{F,\tlm}[x:\lambda])
+\lambda\,\rho(x,y)Q_0
\supset
Q(O,\delta(y))\cap \Tc^{[1]}_{F,\tlm}[y:\lambda]=
\Prm(O,\Tc^{[1]}_{F,\tlm}[y:\lambda]).
$$
\par By interchanging the roles of $x$ and $y$, we obtain also
$$
\Prm(O,\Tc^{[1]}_{F,\tlm}[y:\lambda])
+\lambda\,\rho(x,y)Q_0
\supset\Prm(O,\Tc^{[1]}_{F,\tlm}[x:\lambda]).
$$
These two inclusions imply the following inequality:
\bel{TH-PR}
\dhf\left(\Prm(O,\Tc^{[1]}_{F,\tlm}[x:\lambda]),
\Prm(O,\Tc^{[1]}_{F,\tlm}[y:\lambda])\right)
\le \lambda\,\rho(x,y),~~~~x,y\in\Mc.
\ee
\par As we have noted above, the set $\Prm(O,\Tc^{[1]}_{F,\tlm}[x:\lambda])$ is either a singleton or a compact line segment. Thus $\Prm(O,\Tc^{[1]}_{F,\tlm}[x:\lambda])$ is a {\it bounded rectangle}, so that, thanks to part (ii) of Claim \reff{TWO} and \rf{TH-PR}, the mapping $g_F$ defined by formula \rf{G-CNTN} is Lipschitz with Lipschitz seminorm at most $\lambda$.
\smsk
\par Let us also not that, thanks to \rf{G-IN-T} and \rf{LM-12}, for every $x\in\Mc$, we have
$$
g_F(x)\in\Tc_{F,\tlm}(x)=\Hc[F^{[1]}[x:\tlm]].
$$
These properties of $g_F$ show that the mapping $g=g_F$ satisfies conditions \rf{WG-HF} and \rf{G-LIP}.
\smsk
\par Thus, all conditions of the hypothesis of Proposition \reff{LSEL-F} are satisfied for $g$ and $f=f_{\lmv;F}$.
This proposition tells us that $f$ is a selection of $F$ with Lipschitz seminorm
$$
\|f\|_{\Lip(\Mc)}=\|f_{\lmv;F}\|_{\Lip(\Mc)}\le \tlm +2\lambda=\lambda_1+2\lambda_2.
$$
\par The proof of Theorem \reff{ALG-T} is complete.\bx
\msk
\par \underline{\sc Proof of Theorem \reff{PA-W}.} Let $\Mf=\MR$ be a finite pseudometric space, and let $F:\Mc\to\CRT$ be a set-valued mapping. We know that, given $\tlm,\lambda\ge 0$, condition \rf{WNEW} holds for \mbox{every} $x,x',x'',y,y',y''\in\Mc$. Lemma \reff{WF1-NE} tells us that, in this case the set $F^{[1]}[x:\tlm]$ is nonempty for every $x\in\Mc$. Furthermore, part {\it (ii)} of
Proposition \reff{WLS-T} tells us that, for every $x\in\Mc$, the set $\Tc_{F,\tlm}^{[1]}[x:\lambda]$ is nonempty as well. See \rf{TC-NEMP}.
\par Therefore, thanks to \rf{SC-P}, the $\lmv$-PA with $\lmv=(\tlm,\lambda)$ produces the outcome {\bf ``Success''}. Furthermore, part {\it (ii)} of Theorem \reff{ALG-T} tells us that, in this case the $\lmv$-PA
returns the mapping $f_{\lmv;F}:\Mc\to\RT$ having the following properties:
$f_{\lmv;F}$ is a Lipschitz selection of $F$ with $$\|f_{\lmv;F}\|_{\Lip(\Mc)}\le 2\lambda+\tlm.$$ Thus, \rf{FL-LS} holds and the proof of Theorem \reff{PA-W} is complete.\bx
\bsk

\par {\bf 5.2 Other versions of the Projection Algorithm.}
\addtocontents{toc}{~~~~5.2 Other versions of the Projection Algorithm.\hfill \thepage\par\VST}

\msk
\par Theorems \reff{ALG-T} and \reff{PA-W} hold for various versions of the $\lmv$-Projection Algorithm relating to the particular choice of the mapping $g_F$ in {\bf STEP 4} of the algorithm. In particular, the proof of Theorem \reff{ALG-T} presented above shows that the only requirement for $g_F$ is that it be a Lipschitz selection of the set-valued mapping
$$
\Tc_{F,\lambda_1}=
\HR[F^{[1]}[\cdot:\lambda_1]]~~~~
\text{with}~~~~\|g_F\|_{\Lip(\Mc)}\le\lambda_2.
$$
In other words, $g_F$ have to satisfy \rf{WG-HF} and inequality \rf{G-LIP} with constants $\tlm=\lambda_1$ and $\lambda=\lambda_2$.
\smsk
\par Let us indicate some of these versions.
\smsk
\par ($\bigstar 1$) We can define $g_F$ by formula \rf{G-CNT} with replacing the origin $O$ by an {\it arbitrary} point in $\RT$;
\msk
\par  ($\bigstar 2$) Suppose that
\bel{TC-BO}
\text{for every}~~x\in\Mc~~\text{the set}~~ F^{[1]}[x:\lambda_1]~~\text{is bounded}.
\ee
Then the rectangle $\Tc^{[1]}_{F,\lambda_1}[x:\lambda_2]$ (see \rf{RT-1}) defined at {\bf STEP 3} of the algorithm is also bounded for all $x\in\Mc$. In this case, we can define $g_F$ by the formula
\bel{GF-CNT}
g_F(x)=\cent
\left(\Tc^{[1]}_{F,\lambda_1}[x:\lambda_2]\right),
~~~~~x\in\Mc.
\ee
See Fig. 30.
\msk
%
%
\begin{figure}[h!]
\hspace{30mm}
\includegraphics[scale=0.65]{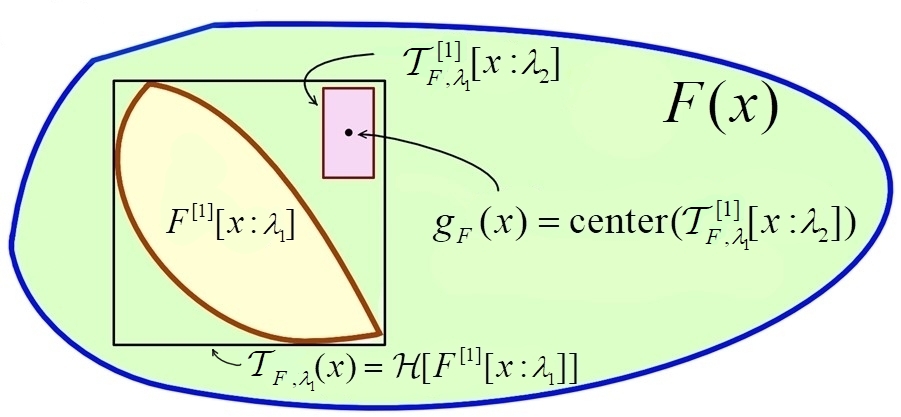}
\caption{{\bf STEP 4} of the Projection Algorithm: the case of bounded sets $F^{[1]}[x:\lambda_1]$.}
\end{figure}
\msk
\par Then, thanks to part (iii) of Proposition \reff{WLS-T}, $g_F$ is Lipschitz with $\|g_F\|_{\Lip(\Mc)}\le\lambda_2$. Thus, for this choice of $g_F$ both property \rf{WG-HF} and inequality \rf{G-LIP} hold.
\par We note that, thanks to Lemma \reff{WF1-NE}, property \rf{TC-BO} holds provided the set $\Mc$ is {\it infinite},
so that in this case we can define $g_F$ by formula \rf{GF-CNT}.
\par Of course, in general we cannot guarantee that property \rf{TC-BO} holds. However, in this case, the following property of the rectangles $\{\Tc^{[1]}_{F,\lambda_1}[x:\lambda_2]:x\in\Mc\}$ maybe useful: if the pseudometric $\rho$ is finite and one of the rectangles of this family is a {\it bounded set}, then {\it all rectangles} from this family are bounded as well, i.e., \rf{TC-BO} holds. This property is immediate  from the fact that at {\bf STEP 4} the mapping $\Tc^{[1]}_{F,\lambda_1}[\cdot:\lambda_2]$ is Lipschitz with respect to the Hausdorff distance. See inequality \rf{T-HD}. (Recall that $\tlm=\lambda_1$ and $\lambda=\lambda_2$.)
\msk
\par ($\bigstar 3$) Because $\Mc$ is finite, there exists $R>0$ such that
$$
\Tc^{[1]}_{F,\lambda_1}[x:\lambda_2]\cap Q(O,R)\ne\emp
~~~~\text{for every}~~~x\in\Mc.
$$
\par This enables us to define $g_F$ by the formula
$$
g_F(x)=\cent
\left(\Tc^{[1]}_{F,\lambda_1}[x:\lambda_2]\cap Q(O,R)\right),~~~~~x\in\Mc.
$$
\par Clearly,
$$
g_F(x)\in\Tc^{[1]}_{F,\lambda_1}[x:\lambda_2] \subset\Tc_{F,\lambda_1}(x)~~~\text{on}~~\Mc
$$
so that property \rf{WG-HF} holds. The proof of \rf{G-LIP} in this case follows the same scheme as the proof of this inequality for $g_F$ defined by formula \rf{G-CNT}. We leave the details to the interested reader.

\SECT{6. From geometric description to analytical: five basic algorithms.}{6}

\addtocontents{toc}{6. From geometric description to analytical: five basic algorithms. \hfill\thepage\par \VST}


\indent
\par Let us proceed to a detailed analytical description of the $\lmv$-Projection Algorithm.
\par First, let us describe {\it the input data}.
\begin{statement}\lbl{INP-S} The inputs of the $\lmv$-PA are the following objects:
\par ($\blbig 1$) A positive integer $N$;
\smsk
\par ($\blbig 2$) A vector $\lmv=(\lambda_1,\lambda_2)$ with non-negative coordinates $\lambda_1$ and $\lambda_2$;
\smsk
\par ($\blbig 3$) An $N\times N$ distance matrix $\Df=\{r_{ij}:i,j=1,...,N\}$ with real entries $r_{ij}\ge 0$ satisfying the standard (pseudo)metric axioms:
\smsk
\par (i) $r_{ii}=0$~ for all~ $1\le i\le N$;
\smsk
\par (ii) $r_{ij}=r_{ji}$~ for all~ $1\le i,j\le N$, i.e., $\Df$ is a symmetric matrix;
\smsk
\par (iii) $r_{ij}\le r_{ik}+r_{kj}$~ for all~ $i,j,k=1,...,N$ (the triangle inequality).
\smsk
\par ($\blbig 4$) $N$ two dimensional non-zero vectors $\lh_i=(a_i,b_i)\in\RT$, and $N$ real numbers $\alpha_i$, $i=1,...,N$.
\end{statement}
\par These inputs generate a finite pseudometric space
\bel{MSP}
\Mf=\MR~~\text{with elements}~~x_i=i,~~i=1,...,N,
\ee
a pseudometric
\bel{RHO}
\rho(x_i,x_j)=r_{ij},~~~~i,j=1,...,N,
\ee
and two mappings, $\lh:\Mc\to\RT\setminus\{0\}$ and $\alpha:\Mc\to\R$ defined by the formulae
\bel{H-ALF}
\lh(x)=\lh_i~~~\text{and}~~~\alpha(x)=\alpha_i~~~\text{for}~~~ x=i\in\Mc,~~i=1,...,N.
\ee
We also define a set-valued mapping $F:\Mc\to\HPL$ by letting
\bel{ID-F}
F(x)=\{u\in\RT:\ip{\lh(x),u}+\alpha(x)\le 0\},~~~x\in\Mc.
\ee
Thus, each $F(x)$ is a half-plane in $\RT$ and $\lh(x)$
is a vector directed outside of $F(x)$ and orthogonal its boundary, the straight line $$\ell(x)=\{u\in\RT:\ip{\lh(x),u}+\alpha(x)=0\}.$$
\par Given these inputs, the $\lmv$-Projection Algorithm produces the outcomes described in Theorem \reff{ALG-T}. This procedure follows the approach outlined in Section 1.2 and includes five main algorithms, which we present in this section.
\par As we noted in Section 1.2, our approach to the corresponding algorithms at {\bf STEPS 1-2} of the $\lmv$-PA relies on a classical result on low-dimensional linear programming which we cite below.
\par Consider a linear program in $\R^d$ expressed in the inequality form
\bel{LP-A}
\text{maximize}~~~~z=\ip{c,x}~~~~\text{subject to}~~~
Ax\le b
\ee
where $c\in\R^d$, $b\in\RM$ and the matrix $A\in\R^{m\times d}$ are the input data, and $x\in\R^d$. (The vector inequality $Ax\le b$ in \rf{LP-A} is with respect to the componentwise partial order in $\RM$.)
\par M. E. Dyer \cite{Dy-1984} and N. Megiddo \cite{M-1983-2} found, independently, an efficient algorithm for the linear programming problem \rf{LP-A} in the case $d = 2$.
\begin{theorem}\lbl{MD-L} The algorithms presented in \cite{Dy-1984} and \cite{M-1983-2} give a solution (or say it is unbounded or unfeasible) to the linear program \rf{LP-A} in two variables and $m$ constraints.
\par The total work and storage required by this algorithm are at most $Cm$ where $C$ is an absolute constant.
\end{theorem}
\begin{remark} {\em As usual, the notion of ``feasibility'' that we use in the statement of Theorem \reff{MD-L} means the following: in $O(m)$ running time and $O(m)$ storage, the above algorithm determines {\it whether the feasible set of a linear program} (i.e., the set of points $x\in\R^d$ satisfying  the constraints $Ax\le b$) {\it is empty or not}.\rbx}
\end{remark}

\begin{remark}\lbl{FEAS} {\em Let us recall something of the history of algorithms for the linear programming problem, where the running time is analyzed as a function of the number of constraints $m$ and the number of variables $d$. See \cite{Ch-2018} for more details.
\smsk
\par First, we note that it is also shown in \cite{M-1983-2} that Theorem \reff{MD-L} holds in the three dimensional case.
\par N. Megiddo \cite{M-1984} generalized Theorem  \reff{MD-L} to an arbitrary $d$, proving that
a linear program in $d$ variables and $m$ constraints can be solved in $2^{O(2^d)}m$ time. This result has become a cornerstone of computational geometry, with many applications in low-dimensional settings. See also \cite{DGMW-2018}.
\smsk
\par It was improved by K. L. Clarkson \cite{Cl-1986} and M. E. Dyer \cite{Dy-1986} to $O((5/9)^d 3^{d^2}m)$ and $O(3^{(d+2)^2}m)$ respectively. Further improvements were obtained by P. K. Agarwal, M. Sharir and S. Toledo \cite{AST-1993} to $O(d)^{10d}(\log d)^{2d} m$, by
B. Chazelle and J. Matou\v{s}ek \cite{ChM-1996} to $O(d)^{7d}(\log d)^d m$ and H. Br\"{o}nnimann, B. Chazelle and J. Matou\v{s}ek \cite{BCM-1999} to $O(d)^{7d}(\log d)^d m$. The fastest algorithm known to the moment is due to
T. M. Chan \cite{Ch-2018} with the running time $O(d)^{d/2}(\log d)^{3d} m$.\rbx}
\end{remark}
\bsk

\par {\bf 6.1 The $\lambda$-metric refinement of the set-valued mapping - STEP 1.}

\addtocontents{toc}{~~~~6.1 The $\lambda$-metric refinement of the set-valued mapping - STEP 1.\hfill \thepage\par}

\msk
\par The goal of this section is to provide an efficient algorithm for recognizing whether, for each $x\in\Mc$, the set $F^{[1]}[x:\lambda_1]$ defined by formula \rf{F-LM1} (the $\lambda_1$-metric refinement of $F$) is nonempty or not.
\begin{alg}\lbl{AL-1}  There exists an algorithm with the inputs given in Statement \reff{INP-S}, which
produces one of the following outcomes:\smsk
\par \underline{\sc Outcome 1:} We guarantee that there exists $x\in\Mc$ such that the set $F^{[1]}[x:\lambda_1]$ is empty.
\smsk
\par \underline{\sc Outcome 2:} We guarantee that for each  $x\in\Mc$, the set $F^{[1]}[x:\lambda_1]\ne\emp$.
\smsk
\par The total work and storage required by this algorithm are at most $CN^2$ where $C$ is an absolute constant.
\end{alg}
\par {\it Explanation:} We recall that, for every $x\in\Mc$,
\bel{NF-I}
F^{[1]}[x:\lambda_1]=\bigcap_{y\in\Mc}\,H(x;y)
\ee
where
$$
H(x;y)=F(y)+\lambda_1\,\rho(x,y)\,Q_0,~~~y\in\Mc.
$$
See \rf{F-1}. Thanks to \rf{ID-F},
$$
H(x;y)=\{u\in\RT:\ip{\lh(y),u}+\alpha(y)\le 0\} +\lambda_1\,\rho(x,y)\,Q_0.
$$
Thus, $H(x,y)$ is also a half-plane having the following representation:
\bel{HI-A}
H(x;y)=\{u\in\RT:\ip{\lh(y),u}+\alpha(y)\le  \lambda_1\,\rho(x,y)\,\|\lh(y)\|_1\},
~~~y\in\Mc.
\ee
Recall that
$$
\|\lh\|_1=|h_1|+|h_2|~~~~\text{for}~~~~ \lh=(h_1,h_2)\in\RT
$$
is the $\ell^2_1$-norm in $\RT$.
\smsk
\par Fix $x\in\Mc$ and consider a system of linear inequalities in two variables $(u_1,u_2)=u$:
\bel{SY-1}
\ip{\lh(y),u}+\alpha(y)\le\lambda_1\,\rho(x,y)\,\|\lh(y)\|_1,
~~~~y\in\Mc.
\ee
Thanks to \rf{NF-I} and \rf{HI-A}, this system has a solution if and only if the set $F^{[1]}[x:\lambda_1]\ne\emp$.
\smsk
\par We apply Theorem \reff{MD-L} to the system \rf{SY-1} of linear inequalities with respect to two variables.
This theorem tells us that in $O(N)$ running time with $O(N)$ storage we can check the feasibility of this system, i.e., to say whether it has a solution or not.
\smsk
\par We apply this procedure to each $x=i\in\Mc$, $i=1,...,N$, starting from $x=1$. If for a certain element $\hx=i_0, 1\le i_0\le N$, Theorem \reff{MD-L} will tell us that the system has no a feasible solution, then $F^{[1]}[\hx:\lambda_1]=\emp$. In this case we produce {\sc Outcome 1}.
\par Recall that in this case the $\lmv$-Projection Algorithm also produces the outcome {\bf ``No Go''} (see   Theorem \reff{PA-PG-I}) and terminates.
\smsk
\par If Theorem \reff{MD-L} tells us that {\it for every} $x\in\Mc$ the system \rf{SY-1} has a solution, we produce {\sc Outcome 2}. In this case, the set $F^{[1]}[x:\lambda_1]\ne\emp$~ for every $x\in\Mc$ which corresponds to the outcome {\bf ``Success''} in Theorem \reff{PA-PG-I}.
\par Since we apply the algorithm of Theorem \reff{MD-L} at most $N$ times, and each time this algorithm produces an answer in $O(N)$ running time with $O(N)$ storage, the work and storage required to perform all steps of this procedure are $O(N^2)$.\bx
\bsk\msk

\par {\bf 6.2 The rectangular hulls of metric refinements - STEP 2.}

\addtocontents{toc}{~~~~6.2 The rectangular hulls of metric refinements - STEP 2.\hfill \thepage\par\VST}

\msk
\par In this section we present an efficient algorithm for computing the rectangle $\Tc_{F,\lambda_1}(x)$, see \rf{TC-DF}, i.e., the rectangular hull of the set $F^{[1]}[x:\lambda_1]$, $x\in\Mc$.
\begin{alg}\lbl{AL-2} There exists an algorithm which for any inputs from Statement \reff{INP-S} such that the set
$F^{[1]}[x:\lambda_1]\ne\emp$~ for every $x\in\Mc$, produces $4N$ numbers
$$
s_1(x),s_2(x),t_1(x),t_2(x)\in\R
\cup\{\pm\infty\},
~~s_1(x)\le s_2(x),~t_1(x)\le t_2(x),~~~x\in\Mc,
$$
which determine the rectangle $\Tc_{F,\lambda_1}(x)=\Hc[F^{[1]}[x:\lambda_1]]$~ for every $x\in\Mc$.
\par This means the following: $\Tc_{F,\lambda_1}(x)=I_1(x)\times I_2(x)$ where
$I_1(x)\in\Ic(Ox_1)$ and $I_2(x)\in\Ic(Ox_2)$ are nonempty closed intervals on the axes $Ox_1$ and $Ox_2$.
These intervals have the following properties:
\par The numbers $s_1(x)$ and $s_2(x)$ are the left and right ends of $I_1(x)$, see \rf{LR1-END}, i.e.,
$$
s_1(x)=\inf\{s\in\R:(s,0)\in I_1(x)\}~~~\text{and}~~~
s_2(x)=\sup\{s\in\R:(s,0)\in I_1(x)\},
$$
and $t_1(x)$ and $t_2(x)$ are the left and right ends of $I_2(x)$, see \rf{LR2-END}, i.e.,
$$
t_1(x)=\inf\{t\in\R:(0,t)\in I_2(x)\}~~~\text{and}~~~
t_2(x)=\sup\{t\in\R:(0,t)\in I_2(x)\}.
$$
\par The work and storage required by this algorithm are at most $CN^2$ where $C$ is an absolute constant.
\end{alg}
\par {\it Explanation:} Note that, thanks to \rf{P1XP2} and
\rf{HS-U},
$$
I_i(x)=\Prj_i[\Tc_{F,\lambda_1}(x)]=
\Prj_i[F^{[1]}[x:\lambda_1]],~~~i=1,2.
$$
\par Let $S=F^{[1]}[x:\lambda_1]$, and let $e_1=(1,0)$ and $e_2=(0,1)$. As we have noted in Remark \reff{N-ST},
\bel{ST-P1}
s_1(x)=\inf\{\ip{e_1,u}:u\in S\},~~~
s_2(x)=\sup\{\ip{e_1,u}:u\in S\},
\ee
and
\bel{ST-P2}
t_1(x)=\inf\{\ip{e_2,u}:u\in S\},~~~
t_2(x)=\sup\{\ip{e_2,u}:u\in S\}.
\ee
\par We recall, that, thanks to \rf{NF-I} and \rf{HI-A}, given $x\in\Mc$, the set $S=F^{[1]}[x:\lambda_1]$ is a polygon determined by $N$ linear constrains given in \rf{SY-1}.
\par Given $\xi\in\RT$, let us consider a linear program in $\RT$
\bel{LP-S2}
\text{maximize}~~~~z=\ip{\xi,u}~~~~\text{subject to}~~~
u\in S.
\ee
Then, thanks to \rf{ST-P1} and \rf{ST-P2}, the numbers
$s_2(x)$ and $t_2(x)$ are solutions to this problem where $\xi$ equals to $e_1$ or $e_2$ respectively. In turn, the numbers $s_1(x)$ and $t_1(x)$ are solutions to the linear program
\bel{LP-S1}
\text{minimize}~~~~z=\ip{\xi,u}~~~~\text{subject to}~~~
u\in S,
\ee
with the same vectors $\xi=e_1$ or $\xi=e_2$. 
\par Theorem \reff{MD-L} tells us that each of these four linear programs (i.e., the linear programs \rf{LP-S2} and \rf{LP-S1} for computing the numbers $s_1(x),s_2(x),t_1(x)$ and $t_2(x)$) can be solved in $O(N)$ running time with $O(N)$ storage.
\par Then, we apply this procedure to every $x_i\in\Mc$, $i=1,...,N$, starting from $i=1$. After $N$ steps, we will determine all $N$ rectangles $\Tc_{F,\lambda_1}(x_i)$, $i=1,...,N$, so that the total work and storage required to proceed all steps of this procedure are $O(N^2)$.\bx
\bsk\bsk


\par {\bf 6.3 The metric refinements of the rectangular hulls - STEP 3.}

\addtocontents{toc}{~~~~6.3 The metric refinements of the rectangular hulls - STEP 3.\hfill \thepage\par\VST}

\msk
\par At this step of the $\lmv$-Projection Algorithm our inputs are  $4N$ numbers
$$
s_1(x),s_2(x),t_1(x),t_2(x)\in\R
\cup\{\pm\infty\},
~~s_1(x)\le s_2(x),~t_1(x)\le t_2(x),~~~x\in\Mc,
$$
obtained in Algorithm \reff{AL-2}. These numbers determine the rectangles $\{\Tc_{F,\lambda_1}(x):x\in\Mc\}$, which means that  for each $x\in\Mc$ the numbers $s_1(x)$ and $s_2(x)$ are the left and right ends of $\Prj_1[\Tc_{F,\lambda_1}(x)]$, and $t_1(x)$ and $t_2(x)$ are the left and right ends of $\Prj_2[\Tc_{F,\lambda_1}(x)]$. Let us also note the following properties of these numbers:
$$
s_1(x),~t_1(x)<+\infty~~~\text{and}~~~
s_2(x),~t_2(x)>-\infty.
$$
\par Our aim is to give an efficient algorithm to recognize whether for every $x\in\Mc$ the $\lambda_2$-metric refinement of the rectangle $\Tc_{F,\lambda_1}$, i.e., the set
$$
\Tc^{[1]}_{F,\lambda_1}[x:\lambda_2]=\bigcap_{y\in \Mc}\,
\left[\Tc_{F,\lambda_1}(y)+\lambda_2\,\rho(x,y)\,Q_0\right],
$$
is nonempty or not. Recall that, according to {\bf STEP 3}
of the $\lmv$-Projection Algorithm, if this set is empty for some $x\in\Mc$, this algorithm produces the outcome {\bf ``No Go''} and terminates. If the set $\Tc^{[1]}_{F,\lambda_1}[x:\lambda_2]\ne\emp$ for all $x\in\Mc$, we produce the ends of the projections of this set onto the coordinate axes and turn to {\bf STEP 4} of the $\lmv$-PA.
\begin{alg}\lbl{AL-3} There exists an algorithm with the inputs $\{s_1(x),s_2(x),t_1(x),t_2(x):x\in\Mc\}$ which produces one of the following outcomes:\smsk
\par \underline{\sc Outcome 1:} We guarantee that there exists $x\in\Mc$ such that the set $\Tc^{[1]}_{F,\lambda_1}[x:\lambda_2]$ is empty.
\par \underline{\sc Outcome 2:} We guarantee that, for every $x\in\Mc$, the set $\Tc^{[1]}_{F,\lambda_1}[x:\lambda_2]\ne\emp$.
\smsk
\par Furthermore, at this stage the algorithm produces $4N$ numbers
$$
\s_1(x),\,\s_2(x),\,\st_1(x),\,\st_2(x)\in\R
\cup\{\pm\infty\},
~~~\s_1(x)\le \s_2(x),~~\st_1(x)\le \st_2(x),~~~x\in\Mc,
$$
which for each $x\in\Mc$ determine the rectangle $\Tc^{[1]}_{F,\lambda_1}[x:\lambda_2]$.
\smsk
\par The total work and storage required by this algorithm are at most $CN^2$ where $C$ is a universal constant.
\end{alg}
\par {\it Explanation:} Let us note that the numbers $\s_i(x)$ and $\st_i(x)$, $i=1,2$, determine the rectangle $\Tc^{[1]}_{F,\lambda_1}[x:\lambda_2]$ in the same sense as
in the formulation of Algorithm \reff{AL-2}: the rectangle
$$
\Tc^{[1]}_{F,\lambda_1}[x:\lambda_2]=\tI_1(x)\times \tI_2(x)~~~\text{where}~~~
\tI_i(x)\in\Ic(Ox_i),~~~i=1,2,
$$
are nonempty closed intervals on the axes $Ox_1$ and $Ox_2$. Thus,
$$
\tI_i(x)=\Prj_i[\Tc^{[1]}_{F,\lambda_1}[x:\lambda_2]],
~~~i=1,2.
$$
These intervals have the following properties: $\s_1(x)$ and $\s_2(x)$ are the left and right ends of $\tI_1(x)$, see \rf{LR1-END}, and  $\st_1(x)$ and $\st_2(x)$ are the left and right ends of $\tI_2(x)$, see \rf{LR2-END}.
\msk
\par Our algorithm includes the following two steps:
\smsk

\par ($\bigstar 1$)~ If there exist $x,y\in\Mc$ such that
\bel{CR-CL-RE}
\max\{s_1(x)-s_2(y),s_1(y)-s_2(x),
t_1(x)-t_2(y),t_1(y)-t_2(x)\}\,>\,\lambda_2\,\rho(x,y)
\ee
then we produce {\sc Outcome 1}. We state that in this case
the set $\Tc^{[1]}_{F,\lambda_1}[x:\lambda_2]=\emp$.
\smsk
\par ($\bigstar 2$)~ If condition \rf{EM-TN} does not hold, i.e.,
$$
\max\{s_1(x)-s_2(y),s_1(y)-s_2(x),
t_1(x)-t_2(y),t_1(y)-t_2(x)\}\,\le\,\lambda_2\,\rho(x,y)
$$
for every $x,y\in\Mc$, we produce {\sc Outcome 2}. In this case we guarantee that
$\Tc^{[1]}_{F,\lambda_1}[x:\lambda_2]\ne\emp$ for every $x\in\Mc$.
\smsk
\par Furthermore, Lemma \reff{RF-ENDS} tells us that in this case the rectangle $\Tc^{[1]}_{F,\lambda_1}[x:\lambda_2]$ is determined by the
numbers $\s_i(x),\st_i(x)$, $i=1,2$, defined as follows:
\bel{S12-NF}
\s_1(x)=\max_{y\in\Mc}\,\{s_1(y)-\lambda_2\,\rho(x,y)\},
~~~~~
\s_2(x)=\min_{y\in\Mc}\,\{s_2(y)+\lambda_2\,\rho(x,y)\},
\ee
and
\bel{T12-NF}
\st_1(x)=\max_{y\in\Mc}\,\{t_1(y)-\lambda_2\,\rho(x,y)\},
~~~~~
\st_2(x)=\min_{y\in\Mc}\,\{t_2(y)+\lambda_2\,\rho(x,y)\}.
\ee
\par Let us explain parts ($\bigstar 1$) and ($\bigstar 2$) of our algorithm.
\par Part {\it (a)} of Proposition \reff{X2-C} tells us that $\Tc^{[1]}_{F,\lambda_1}[x:\lambda_2]\ne\emp$
on $\Mc$ if and only if
\bel{EM-TN}
\Tc_{F,\lambda_1}(x)\cap\left[\Tc_{F,\lambda_1}(y)
+\lambda_2\,\rho(x,y)\,Q_0\right]\ne\emp
\ee
for every $x,y\in\Mc$.
\par In turn, Claim \reff{CL-M} tells us that, given $x,y\in\Mc$, property \rf{EM-TN} holds if and only if
\bel{CL-RE}
\max\{s_1(x)-s_2(y),s_1(y)-s_2(x),
t_1(x)-t_2(y),t_1(y)-t_2(x)\}\,\le\,\lambda_2\,\rho(x,y).
\ee
\msk
\par These two statements, criterion \rf{EM-TN} and
inequality \rf{CL-RE}, justify our algorithm introduced in items ($\bigstar 1$) and ($\bigstar 2$).
\msk
\par Let us estimate the work required to execute this algorithm.
\par Clearly, for each pair $x,y$ of elements of $\Mc$ we can decide whether inequality \rf{CR-CL-RE} holds
or not using at most $M$ computer operations where $M>0$ is an absolute constant.
\smsk
\par Thus, to check criteria ($\bigstar 1$) we need to consider at most $N(N+1)/2$ pairs of elements $x,y\in\Mc$.
Because each check requires at most $M$ computer operations, the total work of the algorithm suggested in ($\bigstar 1$) does not exceed $O(N^2)$.
\smsk
\par Furthermore, formulae \rf{S12-NF} and \rf{T12-NF} show that, given $x\in\Mc$, the algorithm requires $O(N)$ running time to compute all numbers $\s_i(x)$, $\st_i(x)$, $i=1,2$. Therefore, the total running time for computing these numbers for all $N$-element set is $O(N^2)$.
\par Let us also note that this algorithm requires at most $C N^2$ units of the computer memory where $C$ is a universal constant.\bx
\bsk\msk

\par {\bf 6.4 The centers of the metric projections onto rectangles - STEP 4.}

\addtocontents{toc}{~~~~6.4 The centers of the metric projections onto rectangles - STEP 4.\hfill \thepage\par\VST}

\msk
\par At this stage of the Projection Algorithm our inputs are $4N$ numbers
\bel{NS-31}
\s_1(x),\,\s_2(x),\,\st_1(x),\,\st_2(x)\in\R
\cup\{\pm\infty\},
~~~\s_1(x)\le \s_2(x),~~\st_1(x)\le \st_2(x),~~~x\in\Mc,
\ee
obtained with the help of Algorithm \reff{AL-3}. These numbers determine the rectangles $\Tc^{[1]}_{F,\lambda_1}[x:\lambda_2]$ for every $x\in\Mc$. Our goal is to present an algorithm that, based on these inputs, constructs the mapping $g_F:\Mc\to\RT$ defined by formula \rf{G-CNT}.
\begin{alg}\lbl{AL-4} There exists an algorithm which for the inputs \rf{NS-31} produces $N$ points $g_F(x)\in\RT$
defined by the formula
$$
g_F(x)=\cent\left
(\Prm\left(O,\Tc^{[1]}_{F,\lambda_1}[x:\lambda_2]\right) \right), ~~~~~~x\in\Mc.
$$
\par The work and storage required by this algorithm are at most $CN$ where $C$ is an absolute constant.
\end{alg}
\par {\it Explanation:} We compute the coordinates of the points $g_F(x)=(g_1(x),g_2(x))\in\RT$, $x\in\Mc$, using the explicit formulae for these coordinates given in Claim \reff{MPR-C}. Indeed, we know that for each $x\in\Mc$, the rectangle
$$
\Tc^{[1]}_{F,\lambda_1}[x:\lambda_2]=I_1(x)\times I_2(x)
$$
where $I_i(x)\in\Ic(Ox_i)$ are closed intervals on the axes $Ox_i$ determined by the numbers $\s_i(x),\,\st_i(x)$, $i=1,2$, from the inputs \rf{NS-31}. More specifically,
$\s_1(x)$ and $\s_2(x)$ are the left and right ends of $I_1(x)$ while $\st_1(x)$ and $\st_2(x)$ are the left and right ends of $I_2(x)$.
\msk
\par Claim \reff{MPR-C} tells us that we can compute the required coordinates $g_1(x)$ and $g_2(x)$ by letting
$$
g_1(x)=g_1~~~~\text{and}~~~~g_2(x)=g_2
$$
where the numbers $g_1$ and $g_2$ are determined by formulae \rf{DS-Q}, \rf{G-2T} and \rf{C-MP} with the initial parameters
$$
a_1=\s_1(x),~~b_1=\s_2(x),~~
a_2=\st_1(x)~~\text{and}~~b_2=\st_2(x).
$$
\par These formulae tell us that, for every $x\in\Mc$, we need at most $C$ arithmetic operations and $C$ units of computer memory to calculate the coordinates of the point $g_F(x)$. Here $C$ is a universal constant.
\par Therefore, the total work and storage required by this algorithm for computing the mapping $g_F$ on $\Mc$ are at most $CN$.\bx
\bsk\msk

\par {\bf 6.5 The final step: metric projections onto the metric refinements - STEP 5.}

\addtocontents{toc}{~~~~6.5 The final step: metric projections onto the metric refinements - STEP 5.\hfill \thepage\par\VST}

\msk
\par We move on to the final part of the Projection Algorithm. At this stage our input data are the inputs given in Statement \reff{INP-S} and the mapping $g_F$ constructed in Algorithm \reff{AL-4}. Below we give an efficient algorithm which produces the mapping $f_{\lmv;F}:\Mc\to\RT$ defined by formula \rf{F-LF}, i.e., the metric projection of $g_F(x)$ onto the set $F^{[1]}[x:\lambda_1]$, ~$x\in\Mc$. The mapping $f_{\lmv;F}$ is the finite goal of the $\lmv$-Projection Algorithm. Theorem \reff{ALG-T} tells us that $f_{\lmv;F}$ is a Lipschitz selection of $F$ with Lipschitz seminorm at most $\lambda_1+2\lambda_2$.
\begin{alg}\lbl{AL-5} There exists an algorithm that, for the mapping $g_F$ and the inputs given in Statement \reff{INP-S}, produces the mapping $f_{\lmv;F}$ defined by the formula
\bel{PRF}
f_{\lmv;F}(x)=\Prm\left(g_F(x),F^{[1]}[x:\lambda_1]\right), ~~~~~~x\in\Mc.
\ee
\par The work and storage required by this algorithm are at most $CN^2$ where $C$ is a universal constant.
\end{alg}
\par {\it Explanation:} We recall that, for every $x\in\Mc$, the set $F^{[1]}[x:\lambda_1]$ is given by the formula
$$
F^{[1]}[x:\lambda_1]=\bigcap_{y\in\Mc}\,H(x;y)
$$
where
\bel{HI-A-L}
H(x;y)=\{u\in\RT:\ip{\lh(y),u}+\alpha(y)-  \lambda_1\,\rho(x,y)\,\|\lh(y)\|_1\le 0\},
~~~y\in\Mc.
\ee
See \rf{NF-I} and \rf{HI-A}.
\msk
\par Given $x\in\Mc$, we construct the point $f(x)=f_{\lmv;F}(x)$, defined by the formula \rf{PRF}, in two steps.
\par We recall that $g_F(x)\in\Hc[F^{[1]}[x:\lambda_1]]$. Part {\it (i)} of Lemma \reff{D-23} tells us that in this case
\bel{D-FY}
\dist\left(g_F(x),F^{[1]}[x:\lambda_1]\right)=
\max_{y\in\Mc}\dist\left(g_F(x),H(x;y)\right).
\ee
\par At the first step of our algorithm, we fix a half-plane $H_0=H(x;y_0)$ for which the maximum in the right hand side of the equality \rf{D-FY} is attained.
\par We can find $y_0$ in $N$-steps because, thanks to formula \rf{DA-HN} and representation \rf{HI-A-L}, for each $y\in\Mc$, we have an explicit formula for the distance from $g_F(x)$ to $H(x;y)$:
$$
\dist(g_F,H(x;y))=[\ip{\lh(y),g_F(x)}+\alpha(y)-  \lambda_1\,\rho(x,y)\,\|\lh(y)\|_1]_+\,/\,\|\lh(y)\|_1.
$$
\par Thus, in $N$ steps we compute the maximum in the right hand side of \rf{D-FY} using at most $C$ operation at each step. This enables us in at most $CN$ running time to determine the half-plane $H_0=H(x;y_0)$ with the property
$$
\dist\left(g_F(x),F^{[1]}[x:\lambda_1]\right)=
\dist\left(g_F(x),H(x;y_0)\right).
$$
\par We know that the point $g_F(x)$ belongs to the rectangular hull of the set $F^{[1]}[x:\lambda_1]$, i.e., to the set $\HR[F^{[1]}[x:\lambda_1]]$. Because $F^{[1]}[x:\lambda_1]\subset H_0$, we have $g_F(x)\in\HR[H_0]$, so that, thanks to Lemma \reff{WMP-S},
the metric projection $\Prm\left(g_F(x),H_0\right)$ of $g_F(x)$ onto $H_0$ is a {\it singleton}.
\par Furthermore, let
$$
R=\dist\left(g_F(x),F^{[1]}[x:\lambda_1]\right)(=
\dist\left(g_F(x),H_0)\right)).
$$
Then
$$
\Prm\left(g_F(x),F^{[1]}[x:\lambda_1]\right)=
F^{[1]}[x:\lambda_1]\cap Q(g_F(x),R)\subset
H_0\cap Q(g_F(x),R)=
\Prm\left(g_F(x),H_0\right).
$$
\par Because $\Prm\left(g_F(x),H_0\right)$ is a singleton, we conclude that
$$
\Prm\left(g_F(x),F^{[1]}[x:\lambda_1]\right) =\Prm\left(g_F(x),H_0\right).
$$
Hence, thanks to \rf{PRF},
$$
f(x)=f_{\lmv;F}(x)=\Prm\left(g_F(x),H_0\right).
$$
\par It remains to compute the coordinates of the point $f(x)=(f_1(x),f_2(x))$ as the nearest point to $g_F(x)$ on the half-plane $H_0=H(x;y_0)$. Recall that, thanks to \rf{HI-A-L},
$$
H_0=H(x;y_0)=\{u\in\RT:\ip{\lh(y_0),u}+\alpha(y_0)-  \lambda_1\,\rho(x,y_0)\,\|\lh(y_0)\|_1\le 0\}.
$$
\par Clearly, $f(x)=g_F(x)$ provided
\bel{G-H0}
\ip{\lh(y_0),g_F(x)}+\alpha(y_0)\le  \lambda_1\,\rho(x,y_0)\,\|\lh(y_0)\|_1
\ee
because in this case $g_F(x)\in H(x;y_0)$.
\smsk
\par If inequality \rf{G-H0} is not satisfied, i.e., $g_F(x)\notin H(x;y_0)$, we can make use of formulae \rf{NP-HPN} and \rf{NP-HN} for the coordinates of the metric projection of a point onto a half-plane in $\RT$.  We set in these formulas
$$
g=g_F(x),~~\lh=\lh(y_0)~~~\text{and}~~~
\alpha=\alpha(y_0)-\lambda_1\,\rho(x,y_0)\,\|\lh(y_0)\|_1,
$$
and get the required point $f(x)=f_{\lmv;F}(x)$. Clearly, these formulas require at most $C$ computer operations to calculate $f(x)=f_{\lmv;F}(x)$. Here $C$ is an absolute constant.
\smsk
\par Summarizing, we conclude that, following the above procedure, for every $x\in\Mc$ we need $O(N)$ running time and $O(N)$ storage to compute the coordinates of the point $f_{\lmv;F}(x)$.
\par Therefore, the total work and storage required by this algorithm for computing the mapping $f_{\lmv;F}$ on $\Mc$ are at most $CN^2$.\bx
\msk
\par We are in a position to prove the main Theorem \reff{PA-PG-I}.
\msk
\par \underline{\sc Proof of Theorem \reff{PA-PG-I}.} Let
$$
\lambda_1=\lambda_2=\lambda~~~\text{and let}~~~ \lmv=(\lambda_1,\lambda_2)=(\lambda,\lambda).
$$
Without loss of generality we can assume that our pseudometric space $\MR$ and the set-valued mapping $F$ are generated by the inputs given in Statement \reff{INP-S} and formulas \rf{MSP} - \rf{ID-F}.
\smsk
\par We successively apply Algorithms \reff{AL-1} - \reff{AL-5} to these initial data. These algorithms correspond to {\bf STEP 1} - {\bf STEP 5} of the $\lmv$-Projection Algorithm and produce either the outcome {\bf(``No Go'')} or the outcome {\bf(``Success'')} described in Theorem \reff{ALG-T}. This theorem tells us that in the case of the outcome {\bf ``No go''} we guarantee that there does not exist a Lipschitz selection of $F$ with Lipschitz seminorm at most
$$
\min\{\lambda_1,\lambda_2\}=\min\{\lambda,\lambda\}
=\lambda.
$$
Thus, this outcome coincides with {\sc Outcome 1} {\bf(``No Go'')} of Theorem \reff{PA-PG-I}.
\par Theorem \reff{ALG-T} also tells us that in the case of
the outcome {\bf ``Success''}, the $\lmv$-Projection Algorithm returns the mapping $f_{\lmv;F}:\Mc\to\RT$, which is a Lipschitz selection of $F$ on $\Mc$ with Lipschitz seminorm
$$
\|f_{\lmv;F}\|_{\Lip(\Mc)}\le \lambda_1+2\lambda_2=\lambda+2\lambda=3\lambda.
$$
Thus, this outcome coincides with {\sc Outcome 1} {\bf(``Success'')} of Theorem \reff{PA-PG-I}.
\smsk
\par Let us also note that the work and memory space needed for each of Algorithms \reff{AL-1} - \reff{AL-5} are at most $O(N^2)$. Therefore, the total work and storage used by $\lmv$-Projection Algorithm are at most $CN^2$,
where $C>0$ is a universal constant.
\par The proof of Theorem \reff{PA-PG-I} is complete.\bx
\bsk\msk

\par {\bf 6.6 Lipschitz selections of polygon-set valued mappings.}

\addtocontents{toc}{~~~~6.6 Lipschitz selections of polygon-set valued mappings.\hfill \thepage\par\VST}

\msk
\par As promised in Section 1.1, in this section we give an analogue of Theorem \reff{FP-LSA} for the case $\tau=0$ and $D=2$.
\par Let $L$ be a positive integer, and let $\PLG$ be the family of all polygons in $\RT$ (not necessarily bounded), each of which is defined by at most $L$ linear constraints.
\begin{theorem}\lbl{PLG} There exists an algorithm which, given an $N$-element pseudometric space $\MR$, a set-valued mapping $F:\Mc\to\PLG$, and a constant $\lambda>0$, produces the following two outcomes:

\vspace*{3mm}
\par \underline{\sc Outcome 1} {\bf(``No-Go'')}: The algorithm guarantees that there does not exist a Lipschitz mapping $f:\Mc\to\RT$ with Lipschitz constant at most $\lambda$ such that $f(x)\in F(x)$ for all $x\in \Mc$.
\smsk
\par \underline{\sc Outcome 2} {\bf(``Success'')}: The algorithm returns a Lipschitz mapping $f:\Mc\to\RT$, with Lipschitz constant at most $3\lambda$ such that $f(x)\in F(x)$ for all $x\in \Mc$.

\vspace*{3mm}
\par This Selection Algorithm requires at most
$CL^2\,N^2$ computer operations and at most $CL^2\,N^2$ units of computer memory. Here $C>0$ is a universal constant.
\end{theorem}
\par {\it Proof.} We know that for every $x\in\Mc$ the  polygon $F(x)$ can be represented as an intersection of at most $L$ half-planes. We denote this family of half-planes by $\Hcr_F(x)$. Thus, $\#\Hcr_F(x)\le L$ and
\bel{FH-CR}
F(x)=\bigcap_{H\in\,\Hcr_F(x)}\hspace*{-2mm}H
~~~~~~~\text{for every}~~~ x\in\Mc.
\ee
\par Let us introduce a new pseudometric space $\tMf=(\tMc,\trh)$ whose elements are all couples $u=(x,H)$ where $x\in\Mc$ and $H\in\Hcr_F(x)$. We define a pseudometric $\trh$ on $\tMc$ by letting
\bel{D-RTL}
\trh(u,u')=\rho(x,x')~~~\text{for every} ~~~u=(x,H),~u'=(x',H')\in\tMc.
\ee
\par Finally, we define a new set-valued mapping $\tF:\tMc\to\HPL$ by letting
\bel{TF-4}
\tF((x,H))=H~~~\text{provided}~~~x\in \Mc~~\text{and}~~ H\in\Hcr_F(x).
\ee
\par In particular,
$$
\trh((x,H),(x,H'))=0~~~\text{for every}~~~ x\in\Mc~~~\text{and}~~~H,H'\in\Hcr_F(x).
$$
It is also clear that
$$
\#\tMc=\sum_{x\in\Mc}\,\#\Hcr_F(x)\,\le L\cdot N.
$$
\par Let us apply Theorem \reff{PA-PG-I} to the constant $\lambda$, the pseudometric space $\tMf=(\tMc,\trh)$ and  the set-valued mapping $\tF:\tMc\to\HPL$. Theorem \reff{PA-PG-I} tells us that the Projection Algorithm with these inputs produces one of the two following outcomes:
\smsk
\par ($\blbig\tilde{1}$) \underline{\sc Outcome} {\bf(``No Go''):} We guarantee that there does not exist $\tf\in\Lip(\tMc,\trh)$ with Lipschitz seminorm $\|\tf\|_{\Lip(\tMc,\trh)}\le \lambda$ such that $\tf(u)\in \tF(u)$ for all $u\in\tMc$.
\smsk
\par $(\blbig \tilde{2})$ \underline{\sc Outcome} {\bf(``Success''):} The algorithm produces a mapping $\tf:\tMc\to\RT$ with Lipschitz seminorm
$\|\tf\|_{\Lip(\tMc,\trh)}\le 3\lambda$ such that
$\tf(u)\in \tF(u)$ for all $u\in\tMc$.
\par Furthermore, this algorithm requires at most $C(\#\tMc)^2$ computer operations and at most $C(\#\tMc)^2$ units of computer memory. Here $C>0$ is an absolute constant.
\smsk
\par Our algorithm produces
\smsk
\par {\it (i)} {\sc Outcome 1} {\bf(``No-Go'')} if {\it the outcome ($\blbig\tilde{1}$) holds} and
\smsk
\par {\it (ii)} {\sc Outcome 2} {\bf(``Success'')} provided {\it the outcome ($\blbig\tilde{2}$) holds}.
\msk
\par Prove that these outcomes satisfy the requirements of Theorem \reff{PLG}. Indeed, let us show that in the case of the outcome ($\blbig\tilde{1}$) the algorithm guarantees that there does not exist a Lipschitz mapping $f:\Mc\to\RT$ with Lipschitz constant at most $\lambda$ such that $f(x)\in F(x)$ for all $x\in \Mc$.
\smsk
\par Suppose that it is not true, i.e., there exist a mapping $f\in\Lip(\Mc)$ with
\bel{NF-LP3}
\|f\|_{\Lip(\Mc)}\le \lambda
\ee
such that $f(x)\in F(x)$ for all $x\in \Mc$. This mapping generates a new mapping $\tf:\tMc\to\RT$ defined as follows:
$$
\tf((x,H))=f(x)~~~\text{for every}~~~x\in\Mc ~~~\text{and}~~~H\in\Hcr_F(x).
$$
\par Then, thanks to \rf{FH-CR}, \rf{D-RTL} and \rf{TF-4}, for every $x\in\Mc$ and $H\in\Hcr_F(x)$, we have
$$
\tf((x,H))=f(x)\in F(x)=\cap
\{\tH:\tH\in\Hcr_F(x)\}\subset H=\tF((x,H))
$$
which proves that $\tf(u)\in\tF(u)$ for every $u\in\tMc$.
\par Furthermore, let
$$
u=(x,H),u'=(x',H')\in\tMc,
$$
i.e., $H\in\Hcr_F(x)$ and $H'\in\Hcr_F(x')$. Then, thanks to \rf{NF-LP3},
$$
\|\tf(u)-\tf(u')\|=\|\tf((x,H))-\tf((x,H'))\|=
\|f(x)-f(x')\|\le \lambda\,\rho(x,x')=\lambda\,\trh(u,u').
$$
Hence, $\|\tf\|_{\Lip(\tMc,\trh)}\le\lambda$.
\par The existence of the mapping $\tf$ with these properties contradicts the requirement of the outcome ($\blbig\tilde{1}$) proving that the statement of {\sc Outcome 1} {\bf(``No-Go'')} of Theorem \reff{PLG} holds.
\smsk
\par Let us see that if the outcome ($\blbig\tilde{2}$) holds, then {\sc Outcome 2} {\bf(``Success'')} of Theorem \reff{PLG} holds as well. In fact, let us define a mapping $f:\Mc\to\RT$ as follows: given $x\in\Mc$ let us fix some  $H\in\Hcr_F(x)$ and set
\bel{DF-LT}
f(x)=\tf((x,H)).
\ee
\par Prove that $f$ is well defined on $\Mc$, i.e., the value of $f(x)$ does not depend on the choice of the half-plane $H\in\Hcr_F(x)$ in definition \rf{DF-LT}.
\par Indeed, thanks to \rf{D-RTL}, for every $H,H'\in\Hcr_F(x)$, we have
$$
\|\tf((x,H))-\tf((x,H'))\|\le 3\lambda\,\trh(u,u')=3\lambda\,\rho(x,x)=0.
$$
Thus,
$$
\tf((x,H))=\tf((x,H'))
$$
for every $H,H'\in\Hcr_F(x)$
proving that the value of $f(x)$ in \rf{DF-LT} is well defined.
\par Furthermore, given $x,x'\in\Mc$ and $H\in\Hcr_F(x)$ and $H'\in\Hcr_F(x')$, we have
$$
\|f(x)-f(x')\|=\|\tf((x,H))-\tf((x,H'))\|
\le 3\lambda\,\trh((x,H),(x,H'))=3\lambda\,\rho(x,x')
$$
proving that $\|f\|_{\Lip(\Mc)}\le 3\lambda$.
\par Finally, let $x\in\Mc$ and $H\in\Hcr_F(x)$. Then,
$$
f(x)=\tf((x,H))\in\tF((x,H))=H.
$$
From this and \rf{FH-CR}, we have
$$
f(x)\in\cap \{H:H\in\Hcr_F(x)\}=F(x)
$$
proving that $f$ is a selection of $F$.
\par Thus, $f$ satisfies all the requirements of  {\sc Outcome 2} {\bf(``Success'')} of Theorem \reff{PLG} proving that this outcome holds.
\par It remains to note that the Projection Algorithm produces the outcomes ($\blbig\tilde{1}$) and ($\blbig\tilde{2}$) using the work and storage at most $C(\#\tMc)^2$ where $C>0$ is an absolute constant.
Therefore, the Selection Algorithm proposed in this proof produces {\sc Outcome 1} {\bf(``No-Go'')} and {\sc Outcome 2} {\bf(``Success'')} using the work and storage at most
$CL^2\,N^2$.
\par The proof of the theorem is complete.\bx
\msk
\par Let $\PLGA$ be the family of all (not necessarily bounded) polygons in $\RT$. Each of these polygons is determined by a finite number of linear constraints. Given a polygon $P\in\PLGA$ we let $L(P)$ denote the number of linear constraints determining this polygon.
\par The method of proof of Theorem \reff{PLG} allows us to show that the following, somewhat more general version of this theorem is true.
\begin{theorem} There exists an algorithm which, given a finite pseudometric space $\MR$, a set-valued mapping $F:\Mc\to\PLGA$, and a constant $\lambda>0$, produces {\sc Outcome 1} {\bf(``No-Go'')} and
{\sc Outcome 2} {\bf(``Success'')} described in Theorem \reff{PLG}.
\smsk
\par This Selection Algorithm requires at most
$CM^2$ work and storage where $C>0$ is a universal constant and
$$
M=\sum_{x\in\Mc}\, L(F(x))
$$
is the total number of linear constraints determining the polygons $\{F(x):x\in\Mc\}$.
\end{theorem}
\SECT{7. Further results and comments}{7}

\addtocontents{toc}{7. Further results and comments. \hfill \thepage\par\VST}

\indent
\par {\bf 7.1 Generalizations of the PA: sets determined by polynomial inequalities and more.}

\addtocontents{toc}{~~~~7.1 Generalizations of the PA: sets determined by polynomial inequalities and more.\hfill \thepage\par\VST}

\msk
\par Theorem \reff{PA-PG-I} can be generalized to wider families of convex subsets in $\RT$ than the family $\HPL$ of all closed half-planes. (We study this problem in the forthcoming paper \cite{S-2025}.) Let us discuss some of these generalizations.
\par Let $\Tf$ be a family of convex closed subsets of $\RT$. Analyzing the projection algorithm for family $\Tf$, we conclude that {\bf STEP 1} and {\bf STEP 2} are the most complicated parts of this algorithm.
\smsk
\par For instance, at {\bf STEP 2} we construct the rectangular hull $\Hc[F^{[1]}[x:\lambda]]$ of the set $F^{[1]}[x:\lambda]$ defined by formula \rf{F-LM1}. The problem of computing of this rectangular is equivalent to the following problem of convex optimization:
\bel{C-OPT}
\text{maximize}~~~~z=u_1~~~~\text{subject to}~~~
u=(u_1,u_2)\in \bigcap_{y\in \Mc}\,
\left[F(y)+\lambda\,\rho(x,y)\,Q_0\right].
\ee
Here $\lambda>0$, $x\in\Mc$, and the set-valued mapping
$F:\Mc\to\Tf$ defined on the pseudometric space $\MR$
are the input data.
\smsk
\par As we have noted in Section 6.2, if $\Tf=\HPL$ is the family of all closed half-planes in $\RT$, then, thanks to  the results of  N. Megiddo \cite{M-1983-2} and M. E. Dyer \cite{Dy-1984}, the program \rf{C-OPT} can be solved in $O(m)$ time provided $\#\Mc=m$. See Theorem \reff{MD-L}. This allows us, for each $x\in\Mc$, to construct the set $\Hc[F^{[1]}[x:\lambda]]$ in $O(m)$ time. Applying this algorithm to all elements $x\in\Mc$, we build the family of rectangles $\{\Hc[F^{[1]}[x:\lambda]]:x\in\Mc\}$ using at most $O(m^2)$ computer operations.
\smsk
\par Thus, if we want to have similar estimates at {\bf STEP 2} of the PA for a set $\Tf$ of convex closed subsets in $\RT$, we need an analog of Theorem \reff{MD-L} for the convex program \rf{C-OPT}.
\smsk
\par A family of sets $\Tf$ with this property was described by J. Puerto, A. M. Rodr\'{\i}guez-Ch\'{\i}a and A. Tamir in \cite{PRT-2010}. Let $k$ be a positive integer. Let us consider the following convex programming problem, defined on $\RT$:
\bel{GCV-O}
\text{maximize}~~~~z=u_1~~~~\text{subject to}~~~
u=(u_1,u_2)\in\bigcap_{i=1}^m S_i.
\ee
Here
$$
S_i=\{x\in\RT: P_i(x)\le 0\}, ~~i=1,...,m,
$$
where each $P_i=P_i(x_1,x_2)$, $i=1,...,m$, is a {\it convex} polynomial on $\RT$ of degree at most $k$. It is shown in \cite[Section 3]{PRT-2010} that this problem can be solved in $O(m)$ time. The proof of this statement relies on the results of the works \cite{Am-1994} and \cite{ChM-1996}.
\par Basing on this property, we can show that the convex program \rf{C-OPT} also can be solved in $O(m)$ time provided each set $S\in\Tf$ can be represented in the form
$$
S=\{x\in\RT: P(x)\le 0\}
$$
where $P$ is a {\it convex} polynomial on $\RT$ of degree at most $k$.
\par We let $\Pc\Cc_k(\RT)$ denote the family of all such sets $S$. For instance, if $k=2$, as a family $\Tf$ we can take the family of all disks or all ellipses in $\RT$.
\par  Another approach to the problem \rf{GCV-O} with the family of convex sets even wider than $\Pc\Cc_k(\RT)$ was suggested by M. E. Dyer in \cite{Dy-1992}.
\msk
\par T. M. Chan \cite{Ch-1998} proposed an approach that solves the problem \rf{GCV-O} for a very large family $\Sc=\{S_1,...,S_m\}$ of convex closed subsets of $\RT$, but with somewhat more computer operations.
\par Let us introduce two operations with the sets from $\Sc$:
\smsk
\par (D1)~~ Given $S,S'\in\Sc$ find the ends the points $v_{min},v_{max}\in S\cap S'$ such that the projections of the line segment $[v_{min},v_{max}]$ and the intersection $S\cap S'$ onto the axis $Ox_1$ coincide;
\smsk
\par (D2)~~ Given $S\in\Sc$ and a straight line $\ell\subset\RT$, compute the line segment $S\cap\ell$.
We refer to (D1)and (D2) as {\it primitive operations}.
\smsk
\par The main result of \cite{Ch-1998} states that {\it the convex program \rf{GCV-O} in $\RT$ with $m$ constraints  can be solved using $O(m\log m)$ primitive operations}.
\smsk
\par Applying this result to the optimization problem \rf{C-OPT} we can essentially extend the class of sets in the Projection Algorithm. But in this version, this algorithm will require at most $C N^2\log N$ computer operations (and $C N^2$ units of computer memory). Here $C>0$ is an absolute constant.
\bsk\msk

\par {\bf 7.2 A polynomial algorithm for the sharp selection problem.}

\addtocontents{toc}{~~~~7.2 A polynomial algorithm for the sharp selection problem.\hfill \thepage\par\VST}

\msk
\par Let $\lambda>0$, an $N$-point pseudometric space $\MR$, and a set-valued mapping $F:\Mc\to\HPL$ be the same as in Theorem \reff{PA-PG-I}. We ask the following question: {\it Does there exist a polynomial algorithm which provides the same outcomes as in this theorem, but with the estimate $\|f\|_{\Lip(\Mc)}\le \lambda$ instead of inequality \rf{G-3}?}
\smsk
\par We refer to this problem as {\it the sharp selection problem}. The next theorem gives an affirmative answer to this problem.
\begin{theorem} There exists a polynomial algorithm which receives as input a real number $\lambda>0$, a finite pseudometric space $\MR$ with $\#\Mc=N$, and a set-valued mapping $F:\Mc\to\HPL$, and produces the same outcomes as in Theorem \reff{PA-PG-I} (i.e., {\bf (``Success'')} and {\bf (``No Go'')}) with the following refinement of the outcome {\bf (``Success'')}: in this case the algorithm produces a Lipschitz selection $f$ of $F$ with Lipschitz seminorm $\|f\|_{\Lip(\Mc)}\le \lambda$.
\par This algorithm requires at most $CN^4\log N$ computer operations and at most $CN^3 $ units of computer memory. Here $C>0$ is an absolute constant.
\end{theorem}
\par {\it Proof.} This result is immediate from polynomial algorithms for solving systems of linear inequa\-li\-ties of order $m\times n$ {\it with at most two variables per inequality}. N. Megiddo \cite{M-1983-1} was the first to discover such an algorithm requiring $O(mn^3\log m)$ computer operations. Then E. Cohen and N. Megiddo \cite{CM-1996} obtained a new $O(mn^2(\log m+\log^2 n))$ running time algorithm. The fastest algorithm known so far for this problem is by D. S. Hochbaum and J. Naor \cite{HN-1994} with $O(mn^2\log m)$ running time. Note also that all these algorithms require $O(mn)$ units of computer memory.
\smsk
\par Let us apply these algorithms to our problem. This problem can be formulated as follows: construct an algorithm which, in polynomial time, given $\lambda>0$,  $\MR$ with $\#\Mc=N$ and $F:\Mc\to\HPL$, either tells us that no exists a Lipschitz selection of $F$ with Lipschitz seminorm at most $\lambda$ (outcome {\bf (``No Go'')}) or produces such a selection (outcome {\bf (``Success'')}).
\par In turn, this problem is equivalent the feasibility of the following linear programming problem: Let $\Mc=\{x_1,...,x_N\}$ and let
$$
F(x_i)=\{(u,v)\in\RT: a_i u+b_i v\le c_i\}, ~~i=1,...,N.
$$
\par Let us consider a system of linear inequalities with respect to the $2N$ variables $u_1,v_1,...,u_N,v_N$ defined as follows:
\be
\text{Linear Program:} && a_i u_i+b_i v_i\le c_i,~~~i=1,...,N,\lbl{LPR-L}\\
&& u_i-u_j\le \lambda\,\rho(x_i,x_j), ~~~
v_i-v_j\le \lambda\,\rho(x_i,x_j),\nn\\
&& u_j-u_i\le \lambda\,\rho(x_i,x_j), ~~~
v_j-v_i\le \lambda\,\rho(x_i,x_j),\nn\\
&& i,j=1,...,N.
\nn
\ee
\par This is a system of $m=N+N(N+1)=N(N+2)$ linear inequalities with respect to $n=2N$ variables $u_1,v_1,...,u_N,v_N$ where each inequality involves at most two variables.
\par Our goal is either to find a point satisfying all inequalities (the feasible solution of the system) or to conclude that no such point exists (i.e., Linear Program \rf{LPR-L} is infeasible).
\par Applying to this Linear Program the algorithm suggested by D. S. Hochbaum and J. Naor in \cite{HN-1994}, we conclude that this algorithm solves this problem in
$$
O(mn^2\log m)=O(N(N+2)(2N)^2\log N(N+2))=O(N^4\log N)
$$
running time. Also, this algorithm requires
$$
O(mn)=O(N(N+2)(2N))=O(N^3)
$$
units of computer memory.
\par The proof of the theorem is complete.\bx
\begin{remark} {\em Let us note that applying to Linear Program \rf{LPR-L} general-purpose linear programming we obtain algorithms with running time depending on $N$ exponentially.
\par For instance, as we noted in Remark \reff{FEAS}, the fastest of them known to the moment is the algorithm presented by M. Chan in \cite{Ch-2018}. For solving linear programs with $m$ constraints and $n$ variables this algorithm requires $O(n)^{n/2}(\log n)^{3n}m$ running time. In the case of Linear Program \rf{LPR-L}, we have $n=2N$ and $m=N(N+2)$, which leads to an algorithm with the exponential grough in $N$ the number  of computer operations.\rbx}
\end{remark}
\bsk
\par {\bf 7.3 Distance oracle and intrinsic metrics.}

\addtocontents{toc}{~~~~7.3 Distance oracle and intrinsic metrics.\hfill \thepage\par\VST}
\msk
\par In this section we give some remarks related to the following version of Problem \reff{LS-PRB} formulated by
C. Fefferman \cite[Chapter 8.7]{FI-2020}.
\begin{problem}\lbl{PR-F-B}{\em Let $D,L$ be positive integers. Suppose that $\MR$ is an $N$-point metric space. Given $x,y\in\Mc$, we suppose that {\it a distance oracle} tells us $\rho(x, y)$ (perhaps only up to a universal constant factor). For each $x\in\Mc$, suppose we are given a convex polytope $F(x)\subset\RD$, deﬁned by at most $L$ linear constraints.
\par We would like to compute a function $f:\Mc\to\RD$ satisfying $f(x)\in F(x)$ for all $x\in\Mc$, with Lipschitz seminorm as small as possible up to a factor $C(D)$.
How many computer operations does it take?}
\end{problem}
\smsk

\par Our remarks concern the concept of {\it the distance oracle} mentioned in Problem \reff{PR-F-B}. In many cases, in applications, we do not know the true value of the distance $\rho(x,y)$ between any two elements $x,y\in\Mc$ (even up to an absolute constant factor). Instead of this exact value of $\rho(x,y)$, we often know only a certain upper bound for this distance. In other words, our input is a non-negative weight function $W:\Mc\times\Mc\to[0,+\infty]$ defined on pairs of elements of $\Mc$. We assume that
\bel{PR-MC}
W(x,x)=0~~~\text{and}~~~W(x,y)=W(y,x)~~~\text{for all}~~~ x,y\in\Mc.
\ee
\par Note that, in general, {\it $W$ does not necessarily satisfy the triangle inequality}. We refer to the function $W$ as a {\it ``pre-metric''} on $\Mc$.
\par The following selection problem arises in various applications.
\begin{problem}\lbl{W-PR}{\em Let $W$ be non-negative weight function on $\Mc$ satisfying conditions \rf{PR-MC}. For each $x\in\Mc$, suppose we are given a convex closed set $F(x)\subset\RD$. Find an efficient algorithm for constructing a nearly-optimal selection $f$ of $F$ which is ``Lipschitz'' with respect to $W$, i.e., $f(x)\in F(x)$ on $\Mc$, and for some  constant $\lambda\ge 0$ we have
\bel{PR-D}
\|f(x)-f(y)\|\le \lambda\,W(x,y)~~~\text{for all}~~~x,y\in\Mc.
\ee
\par Here ``nearly-optimal'' means ``with $\lambda$ as small as possible'' (up to a universal constant factor).}
\end{problem}
\par Let us note that Problem \reff{W-PR} easily reduces to an equivalent problem where the pre-metric $W$ is replaced with a certain {\it pseudometric $\DI$ generated by $W$}. The pseudometric $\DI$ is the classical object of the graph theory known in the literature as {\it the intrinsic pseudometric} induced by the pre-metric $W$. We recall that $\DI$ is defined as follows: given $x,y\in\Mc$, we set
\bel{GEO}
\DI(x,y)=\inf_{\{u_0,...,u_k\}}\, \sum_{i=0}^{k-1} W(u_i,u_{i+1}),
\ee
where the infimum is taken over all subsets $\{u_0,...,u_k\}\subset\Mc$ such that $u_0=x$ and $u_k=y$.
\par Clearly, $\DI$ is a pseudometric. It is also clear that a function $f:\Mc\to\RD$ satisfies inequality \rf{PR-D} if and only if $f$ is $\lambda$-Lipschitz with respect to $\DI$, i.e.,
$$
\|f(x)-f(y)\|\le \lambda\,\DI(x,y)~~~\text{for every}~~~ x,y\in\Mc.
$$
\par This observation shows that the problems of constructing a nearly optimal Lipschitz selection with respect to $W$ and with respect to $\DI$ are equivalent.
Therefore, at the preparatory stage of our algorithm we can construct $\DI$, and then, at the following stages of the algorithm, deal with the pseudometric $\DI$ rather than  the pre-metric $W$.
\msk
\par The problem of computing the intrinsic metric $\DI$ is one of the classic problems in the algorithmic graph theory. Let us recall something of the history of this problem.
\par Let $G=(V,E)$ be a graph with the sets of the nodes $V$ and the sets of the edges $E$. Let $n=\#V$ and $m=\#E$ be the number of nodes and the number of edges respectively. Finally, let $W:E\to[0,\infty]$ be a length function. (In particular, in these settings, the intrinsic pseudometric $\DI(x,y)$ induced by $W$ (see \rf{GEO}) is the infimum of the lengths of all paths from $x$ to $y$ in the complete graph with the nodes in $\Mc$.)
\smsk
\par Let us note several results related to two problems known in the literature as the {\it Single Source Shortest Path problem}, SSSP in short, and the {\it All Pairs Shortest Path problem}, APSP in short.
\smsk
\par The SSSP problem is the problem of finding the length of the shortest path from a single-source node $x\in V$ to all other nodes $y\in V$. The APSP is to find, for each pair of nodes, $x,y\in V$, the length of the shortest path from $x$ to $y$. (A length of a path is defined as the sum of the length of its edges.) We refer the reader to the papers \cite{B-1958,Ch-2008,Ch-2010,Ch-2012,CLR-2009,
CW-2021,Di-1959,Fl-1962,
FT-1987,FW-1993,HT-2016,Jo-1977,KNS-2018,Th-1999,Wi-2021} and references therein for numerous results on the SSSP and
APSP algorithms.
\smsk
\par Among the multitude of results devoted to the SSSP we mention the classical Dijkstra’s algorithm \cite{Di-1959}
which solves the SSSP in $O(n^2)$ time; see also Bellman's algorithm \cite{B-1958} with running time $O(mn)$. Let us note various improvements of these result: $O(m)$ \cite{Th-1999} ($W(x)$, $x\in V$, are non-negative integers), $O((m+n)\log n)$ \cite{Jo-1977}, $O(m+n\log n)$ \cite{FT-1987} and
$O(m +n\log n/\log\log n)$ \cite{FW-1993}.
\smsk
\par There is also an extensive literature devoted to different aspects of the APSP problem. We refer the reader to the papers \cite{Ch-2008,Wi-2021} for the detailed reviews of the results obtained in this directions. Here we mention the algorithms with the following running time: $O(n^3)$ \cite{Fl-1962,Wa-1962} (the classical Floyd - Warshall algorithm), $O(n^3(\log\log n/\log n)^{1/3})$ \cite{Fr-1976}, $O(n^3/\log n)$ \cite{Ch-2008},
$O(n^3(\log\log n)^3/\log^2 n)$ \cite{Ch-2010},
$O(n^3\log\log n/\log^2 n)$ \cite{HT-2016}. Finally, the algorithm suggested by T. M. Chan and R. Williams \cite{CW-2021} improves the running time to $O(n^3/2^{30(\log n)^{1/2}})$.
\par For the case of the {\it sparse} graphs (i.e., the graphs in which the number of edges is much less than the maximal number of edges) we mention the algorithms of the papers \cite{Jo-1977,KKP-1993} with the running time
$O(mn + n^2\log n)$, \cite{Pe-2004} with $O(mn\log\alpha(m,n))$ (here $\alpha$ is the inverse-Ackermann function), and \cite{PeR-2005} with the running time $O(mn + n^2\log\log n)$. See also \cite{Ch-2012} for a faster algorithm which provides running time $o(mn)$.
\begin{remark}{\em In Problem \reff{W-PR}, sparse graphs appear in a natural way for a weight $W$ that takes the value $+\infty$. Indeed, we introduce a graph structure on $\Mc$, allowing $x,y\in\Mc$ to be joined by an edge if and only if $W(x,y)<+\infty$. Clearly, in definition \rf{GEO} of the intrinsic pseudometric $\DI$ it suffices to consider only paths $\{u_0,...u_k\}$ from $x$ to $y$ such that $u_i$ and $u_{i+1}$ are joined by an edge
for all $i=0,...,k-1$. In other words, $\DI$ is the intrinsic pseudometric induced by the graph $G=(\Mc,E)$ and the weight function $W|_E$ where $E=\{(x,y)\in\Mc\times\Mc:
W(x,y)<+\infty\}$.
\smsk
\par It may happen that $G$ is a sparse graph, i.e.,
$\#E<< n(n+1)/2$.
\par This will allow us to use the above-mentioned fast APSP algorithms for sparse graphs when constructing the pseudometric  $\DI$.\rbx}
\end{remark}

\fontsize{12}{13.5}\selectfont


\begin{thebibliography}{A}


\bibitem {AST-1993} P. K. Agarwal, M. Sharir, S. Toledo, An Efficient Multi-Dimensional Searching Technique and Its
Applications. Technical Report CS-1993-20. Department of Computer Science, Duke University, 1993.
\bibitem {Am-1994} N. Amenta, Helly-type theorems and generalized linear programming. ACM Symposium on Computational Geometry (San Diego, CA, 1993). Discrete Comput. Geom. 12 (1994), no. 3, 241--261.
\bibitem {AF-1990} J.-P. Aubin, H. Frankowska,
Set-valued analysis, Systems \& Control:
Foundations \& Applications, 2. Birkhauser Boston, 1990.
\bibitem {B-1958} R. Bellman, On a routing problem.
Quart. Appl. Math. 16 (1958), 87–90.
\bibitem {BL-2000}  Y. Benyamini,  J. Lindenstrauss,
Geometric nonlinear functional analysis, Vol. 1,
in: American Mathematical Society Colloquium Publications, 48. American Mathematical Society, Pro\-vidence, RI, 2000. xii+488 pp.
\bibitem {BCM-1999} H. Br\"{o}nnimann, B. Chazelle, J. Matou\v{s}ek, Product range spaces, sensitive sampling, and derandomization. SIAM J. Comput. 28 (1999), no. 5, 1552--1575.
\bibitem {BS-1994} Yu. Brudnyi, P. Shvartsman, Generalizations of Whitney's extension theorem.  Internat. Math. Res. Notices (1994), no. 3, 129--139.
\bibitem {BS-2001} Yu. Brudnyi, P. Shvartsman,
Whitney Extension Problem for Multivariate
$C^{1,\omega}$-functions. Trans. Amer. Math. Soc.
353 (2001), no. 6, 2487--2512.
\bibitem {Ch-1998} T. M. Chan, Deterministic algorithms for 2-d convex programming and 3-d online linear programming. J. Algorithms 27 (1998), no. 1, 147--166.
\bibitem {Ch-2008} T. M. Chan, All-pairs shortest paths with real weights in $O(n^3/\log n)$ time. Algorithmica 50 (2008), no. 2, 236--243.
\bibitem {Ch-2010} T. M. Chan, More algorithms for all-pairs shortest paths in weighted graphs. SIAM J. Comput. 39 (2010), no. 5, 2075--2089.
\bibitem {Ch-2012} T. M. Chan, All-pairs shortest paths for unweighted undirected graphs in $o(mn)$ time. ACM Trans. Algorithms 8 (2012), no. 4, Art. 34, 17 pp.
\bibitem {Ch-2018} T. M. Chan, Improved deterministic algorithms for linear programming in low dimensions. ACM Trans. Algorithms 14 (2018), no. 3, Art. 30, 10 pp.
\bibitem {CW-2021} T. M. Chan, R. Williams, Deterministic APSP, orthogonal vectors, and more: quickly derandomizing Razborov-Smolensky. ACM Trans. Algorithms 17 (2021), no. 1, Art. 2, 14 pp.
\bibitem {ChM-1996}  B. Chazelle, J. Matou\v{s}ek, On linear-time deterministic algorithms for optimization problems in fixed dimension. J. Algorithms 21 (1996), no. 3, 579--597.
\bibitem {CM-1996} E. Cohen, N. Megiddo, Improved algorithms for linear inequalities with two variables per inequality. SIAM J. Comput. 23 (1994), no. 6, 1313--1347.
\bibitem {Cl-1986} K. L. Clarkson, Linear programming in $O(n\times3^{d^2})$ time. Inform. Process. Lett. 22 (1986), no. 1, 21--24.
\bibitem {CLR-2009} T. H. Cormen, C. E. Leiserson, R. L. Ronald, C. Stein, Introduction to algorithms. Third edition. MIT Press, Cambridge, MA, 2009.
\bibitem {DGK-1963} L. Danzer, B. Gr\"unbaum, V. Klee, Helly's Theorem and its relatives. in: AMS Symposium on Convexity, Seattle, Proceedings of Symposium on Pure Mathematics, Vol. 7, Amer. Math. Soc., Providence, RI, 1963, pp. 101--180.
\bibitem {Di-1959}  E. W. Dijkstra, A note on two problems in connexion with graphs. Numer. Math. 1 (1959), 269--271.
\bibitem {Dy-1984} M. E. Dyer, Linear time algorithms for two- and three-variable linear programs. SIAM J. Comput. 13 (1984), no. 1, 31--45.
\bibitem {Dy-1986} M. E. Dyer, On a multidimensional search problem and its application to the Euclidean
one-centre problem, SIAM J. Comput. 15 (1986), no. 3, 725--738.
\bibitem {Dy-1992} M. E. Dyer, A Class of Convex Programs with Applications to Computational Geometry. 8th Annual Computational Geometry, 6/92, Berlin, Germany, 9--15.
\bibitem {DGMW-2018} M. Dyer, B. G\"artner, N. Megiddo, E. Welzl, Linear programming. Handbook of discrete and computational geometry, 1291--1309, CRC Press, Boca Raton, FL, 2018.
\bibitem {F-2005} C. Fefferman, A sharp form of Whitney's extension theorem. Ann. of Math. (2) 161 (2005), no. 1, 509–577.
\bibitem {F-2009} C. Fefferman, Whitney extension problems and interpolation of data. Bull. Amer. Math. Soc. 46 (2009), no. 2, 207--220.
\bibitem {F-2009-2} C. Fefferman, Fitting a $C^m$-smooth function to data. III. Ann. of Math. (2) 170 (2009), no. 1, 427--441.
\bibitem {F-2013} C. Fefferman, Smooth interpolation of data by efficient algorithms. Excursions in harmonic analysis. Vol. 1, 71--84, Appl. Numer. Harmon. Anal., Birkhäuser/Springer, New York, 2013.
\bibitem {F-2019} C. Fefferman, Unsolved Problems, 11th Whitney Extension Problems Workshop, Trinity College Dublin, August 13-17, 2019,
https://cms-math.net.technion.ac.il/open-problems-whitney/
\bibitem {FI-2020} C. Fefferman, A. Israel, Fitting Smooth Functions to Data. CBMS Regional Conference Series in Mathematics, 135. American Mathematical Society, Providence, RI, 2020. xi+160 pp.
\bibitem {FIL-2016} C. Fefferman, A. Israel, G. K. Luli, Finiteness principles for smooth selection.  Geom. Funct. Anal. 26 (2016), no. 2, 422--477.
\bibitem {FJL-2023} C. Fefferman, F. Jiang, G. K. Luli,
$C^2$-interpolation with range restriction. Rev. Mat. Iberoam. 39 (2023), no. 2, 649--710.
\bibitem {FK-2009} C. Fefferman, B. Klartag, Fitting a $C^m$-smooth function to data. I., Ann. of Math. (2) 169 (2009), no. 1, 315--346.
\bibitem {FK-2009-2} C. Fefferman, B. Klartag, Fitting a $C^m$-smooth function to data. II., Rev. Mat. Iberoam.  25 (2009), no. 1, 49--273.
\bibitem {FP-2019} C. Fefferman, B. Pegueroles, Efficient Algorithms for Approximate Smooth Selection. J. Geom. Anal. 31 (2021), no. 7, 6530--6600.
\bibitem {FS-2018} C. Fefferman, P. Shvartsman, Sharp finiteness principles for Lipschitz selections. Geom. Funct. Anal. 28 (2018), no. 6, 1641--1705.
\bibitem {Fl-1962} R. W. Floyd, Algorithm 97: Shortest path. Comm. ACM, 5 (1962), 345.
\bibitem {Fr-1976} M. L. Fredman, New bounds on the complexity of the shortest path problem. SIAM J. Comput. 5 (1976), no. 1, 83--89.
\bibitem {FT-1987} M. L. Fredman, R. E. Tarjan, Fibonacci heaps and their uses in improved network optimization algorithms. J. Assoc. Comput. Mach. 34 (1987), no. 3, 596--615.
\bibitem {FW-1993}  M. L. Fredman, D. E. Willard, Surpassing the information-theoretic bound with fusion trees. Proceedings of the 22nd Annual ACM Symposium on Theory of Computing (Baltimore, MD, 1990). J. Comput. System Sci. 47 (1993), no. 3, 424--436.
\bibitem {HT-2016} Y. Han, T. Takaoka, An
$O(n^3\,\log\log n/\log^2n)$ time algorithm for all pairs shortest paths. J. Discrete Algorithms 38/41 (2016), 9--19.
\bibitem {HN-1994} D. S. Hochbaum, J. Naor, Simple and fast algorithms for linear and integer programs with two variables per inequality. SIAM J. Comput. 23 (1994), no. 6, 1179--1192.
\bibitem {JL-2021} F. Jiang, G. K. Luli, Algorithms for nonnegative $C^2({\mathbb R}^2)$ interpolation. Adv. Math. 385 (2021), Paper No. 107756, 43 pp.
\bibitem {JLO-2022-1} F. Jiang, G. K. Luli, K. O'Neill, Smooth selection for infinite sets. Adv. Math. 407 (2022), Paper No. 108566, 62 pp.
\bibitem {JLO-2022-2} F. Jiang, G. K. Luli, K. O'Neill, On the Shape Fields Finiteness Principle. IMRN 2022 (2022), no. 23, 18895--18918.
\bibitem {JLLL-2023} F. Jiang, C. Lian, Y. Liang, G. K. Luli, Univariate range-restricted $C^2$ interpolation algorithms. J. Comput. Appl. Math. 425 (2023), Paper No. 115040, 19 pp.
\bibitem {Jo-1977} D. B. Johnson, Efficient algorithms for shortest paths in sparse networks. J. Assoc. Comput. Mach. 24 (1977), no. 1, 1--13.
\bibitem {KKP-1993} D. R. Karger, D. Koller, S. J. Phillips, Finding the hidden path: time bounds for all-pairs shortest paths. SIAM J. Comput. 22 (1993), no. 6, 1199--1217.
\bibitem {KNS-2018} M. J. Kashyop, T. Nagayama, K. Sadakane, Faster algorithms for shortest path and network flow based on graph decomposition. J. Graph Algorithms Appl. 23 (2019), no. 5, 781--813.
\bibitem{M-1983-1} N. Megiddo, Towards a genuinely polynomial algorithm for linear programming. SIAM J. Comput. 12 (1983), no. 2, 347--353.
\bibitem{M-1983-2} N. Megiddo, Linear-time algorithms for linear programming in $\R^3$ and related problems. SIAM J. Comput. 12 (1983), no. 4, 759--776.
\bibitem{M-1984} N. Megiddo, Linear programming in linear time when dimension is fixed. J. Assoc. Comput. Mach. 31 (1984), no. 1, 114--127.
\bibitem{Pe-2004} S. Pettie, A new approach to all-pairs shortest paths on real-weighted graphs. Automata, languages and programming. Theoret. Comput. Sci. 312 (2004), no. 1, 47--74.
\bibitem{PeR-2005} S. Pettie, V. A. Ramachandran, A shortest path algorithm for real-weighted undirected graphs. SIAM J. Comput. 34 (2005), no. 6, 1398--1431.
\bibitem {PR-1992}  K. Przes{\l}awski, L. E. Rybinski, Concepts of lower semicontinuity and continuous selections for convex valued multifunctions. J. Approx. Theory 68 (1992), no. 3, 262--282.
\bibitem {PY-1989} K. Przes{\l}awski, D. Yost, Continuity properties of selectors and Michael's theorem. Mich. Math. J. 36 (1989), no. 1, 113--134.
\bibitem {PY-1995} K. Przes{\l}awski, D. Yost, Lipschitz Retracts, Selectors and Extensions. Mich. Math. J. 42 (1995), no. 3, 555--571.
\bibitem {PRT-2010} J. Puerto, A. M. Rodríguez-Chía, A. Tamir, On the planar piecewise quadratic $1$-center problem. Algorithmica 57 (2010), no. 2, 252--283.
\bibitem {S-1987} P. A. Shvartsman, Traces of functions of Zygmund class. (Russian) Sibirsk. Mat. Zh. 28 (1987), no. 5, 203–215; English transl. in: Sib. Math. J. 28 (1987), 853--863.
\bibitem {S-2001} P. Shvartsman, On Lipschitz selections of affine-set valued mappings. Geom. Funct. Anal.  11  (2001), no. 4, 840--868.
\bibitem {S-2002} P. Shvartsman, Lipschitz selections of set-valued mappings and Helly's theorem. J. Geom. Anal. 12 (2002), no. 2, 289--324.
\bibitem {S-2004} P. Shvartsman, Barycentric selectors and a Steiner-type point of a convex body in a Banach space. J. Funct. Anal. 210 (2004), no. 1, 1--42.
\bibitem {S-2008} P. Shvartsman, The Whitney extension problem and Lipschitz selections of set-valued mappings in jet-spaces. Trans. Amer. Math. Soc. 360 (2008),
    no. 10, 5529--5550.
\bibitem {S-2021-L} P. Shvartsman, The Core of a 2-Dimensional Set-Valued Mapping. Existence Criteria and Efficient Algorithms for Lipschitz Selections of Low Dimensional Set-Valued Mappings.\\
    arXiv: 2010.04540v2
\bibitem {S-2023} P. Shvartsman, Existence Criteria for Lipschitz Selections of Set-Valued Mappings in $\RT$.
 arXiv:2306.14042v1
\bibitem {S-2025} P. Shvartsman, Efficient Algorithms for Lipschitz Selections of Low Dimensional Set-Valued Mappings. (to appear)
\bibitem {Th-1999} M. Thorup, Undirected single-source shortest paths with positive integer weights in linear time. J. ACM 46 (1999), no. 3, 362--394.
\bibitem {Wa-1962} S. Warshall, A theorem on boolean matrices. J. Assoc. Comput. Mach. 9 (1962), 11--12.
\bibitem {Wh-1934} H. Whitney, Analytic extension of differentiable functions defined in closed sets.  Trans. Amer. Math. Soc.  36 (1934), no. 1, 63--89.
\bibitem {Wi-2021} R. Williams, From circuit complexity to faster all-pairs shortest paths. SIAM Rev. 63 (2021), no. 3, 559--582.


\end{thebibliography}
\end{document}